\documentclass[10pt]{article}
\usepackage[english]{babel}
\usepackage{amsfonts}
\usepackage{amsthm}
\usepackage{amsmath}    
\usepackage{amssymb}
\usepackage{bm}         
\usepackage[utf8]{inputenc}
\usepackage{libertine}
\usepackage{graphicx}   
\usepackage{verbatim}   
\usepackage{color}      
\usepackage{subfigure}  
\usepackage{epsfig}
\usepackage{tikz}       
\usepackage{stmaryrd}
\usepackage{hyperref}   
\numberwithin{equation}{section} 
\usepackage{cleveref} 
\usepackage{fancyhdr} 
\usepackage{rotating} 
\usepackage{mathrsfs}
\usepackage{thmtools}
\usepackage{caption}
\usepackage{pgfplots,xparse}
\usepackage{cite}
\usepackage{gensymb}   
\usepackage{standalone} 

\usepackage[square,numbers]{natbib}

\pgfplotsset{
    colormap={plasma}{
      rgb=(0.050383, 0.029803, 0.527975)
      rgb=(0.063536, 0.028426, 0.533124)
      rgb=(0.075353, 0.027206, 0.538007)
      rgb=(0.086222, 0.026125, 0.542658)
      rgb=(0.096379, 0.025165, 0.547103)
      rgb=(0.105980, 0.024309, 0.551368)
      rgb=(0.115124, 0.023556, 0.555468)
      rgb=(0.123903, 0.022878, 0.559423)
      rgb=(0.132381, 0.022258, 0.563250)
      rgb=(0.140603, 0.021687, 0.566959)
      rgb=(0.148607, 0.021154, 0.570562)
      rgb=(0.156421, 0.020651, 0.574065)
      rgb=(0.164070, 0.020171, 0.577478)
      rgb=(0.171574, 0.019706, 0.580806)
      rgb=(0.178950, 0.019252, 0.584054)
      rgb=(0.186213, 0.018803, 0.587228)
      rgb=(0.193374, 0.018354, 0.590330)
      rgb=(0.200445, 0.017902, 0.593364)
      rgb=(0.207435, 0.017442, 0.596333)
      rgb=(0.214350, 0.016973, 0.599239)
      rgb=(0.221197, 0.016497, 0.602083)
      rgb=(0.227983, 0.016007, 0.604867)
      rgb=(0.234715, 0.015502, 0.607592)
      rgb=(0.241396, 0.014979, 0.610259)
      rgb=(0.248032, 0.014439, 0.612868)
      rgb=(0.254627, 0.013882, 0.615419)
      rgb=(0.261183, 0.013308, 0.617911)
      rgb=(0.267703, 0.012716, 0.620346)
      rgb=(0.274191, 0.012109, 0.622722)
      rgb=(0.280648, 0.011488, 0.625038)
      rgb=(0.287076, 0.010855, 0.627295)
      rgb=(0.293478, 0.010213, 0.629490)
      rgb=(0.299855, 0.009561, 0.631624)
      rgb=(0.306210, 0.008902, 0.633694)
      rgb=(0.312543, 0.008239, 0.635700)
      rgb=(0.318856, 0.007576, 0.637640)
      rgb=(0.325150, 0.006915, 0.639512)
      rgb=(0.331426, 0.006261, 0.641316)
      rgb=(0.337683, 0.005618, 0.643049)
      rgb=(0.343925, 0.004991, 0.644710)
      rgb=(0.350150, 0.004382, 0.646298)
      rgb=(0.356359, 0.003798, 0.647810)
      rgb=(0.362553, 0.003243, 0.649245)
      rgb=(0.368733, 0.002724, 0.650601)
      rgb=(0.374897, 0.002245, 0.651876)
      rgb=(0.381047, 0.001814, 0.653068)
      rgb=(0.387183, 0.001434, 0.654177)
      rgb=(0.393304, 0.001114, 0.655199)
      rgb=(0.399411, 0.000859, 0.656133)
      rgb=(0.405503, 0.000678, 0.656977)
      rgb=(0.411580, 0.000577, 0.657730)
      rgb=(0.417642, 0.000564, 0.658390)
      rgb=(0.423689, 0.000646, 0.658956)
      rgb=(0.429719, 0.000831, 0.659425)
      rgb=(0.435734, 0.001127, 0.659797)
      rgb=(0.441732, 0.001540, 0.660069)
      rgb=(0.447714, 0.002080, 0.660240)
      rgb=(0.453677, 0.002755, 0.660310)
      rgb=(0.459623, 0.003574, 0.660277)
      rgb=(0.465550, 0.004545, 0.660139)
      rgb=(0.471457, 0.005678, 0.659897)
      rgb=(0.477344, 0.006980, 0.659549)
      rgb=(0.483210, 0.008460, 0.659095)
      rgb=(0.489055, 0.010127, 0.658534)
      rgb=(0.494877, 0.011990, 0.657865)
      rgb=(0.500678, 0.014055, 0.657088)
      rgb=(0.506454, 0.016333, 0.656202)
      rgb=(0.512206, 0.018833, 0.655209)
      rgb=(0.517933, 0.021563, 0.654109)
      rgb=(0.523633, 0.024532, 0.652901)
      rgb=(0.529306, 0.027747, 0.651586)
      rgb=(0.534952, 0.031217, 0.650165)
      rgb=(0.540570, 0.034950, 0.648640)
      rgb=(0.546157, 0.038954, 0.647010)
      rgb=(0.551715, 0.043136, 0.645277)
      rgb=(0.557243, 0.047331, 0.643443)
      rgb=(0.562738, 0.051545, 0.641509)
      rgb=(0.568201, 0.055778, 0.639477)
      rgb=(0.573632, 0.060028, 0.637349)
      rgb=(0.579029, 0.064296, 0.635126)
      rgb=(0.584391, 0.068579, 0.632812)
      rgb=(0.589719, 0.072878, 0.630408)
      rgb=(0.595011, 0.077190, 0.627917)
      rgb=(0.600266, 0.081516, 0.625342)
      rgb=(0.605485, 0.085854, 0.622686)
      rgb=(0.610667, 0.090204, 0.619951)
      rgb=(0.615812, 0.094564, 0.617140)
      rgb=(0.620919, 0.098934, 0.614257)
      rgb=(0.625987, 0.103312, 0.611305)
      rgb=(0.631017, 0.107699, 0.608287)
      rgb=(0.636008, 0.112092, 0.605205)
      rgb=(0.640959, 0.116492, 0.602065)
      rgb=(0.645872, 0.120898, 0.598867)
      rgb=(0.650746, 0.125309, 0.595617)
      rgb=(0.655580, 0.129725, 0.592317)
      rgb=(0.660374, 0.134144, 0.588971)
      rgb=(0.665129, 0.138566, 0.585582)
      rgb=(0.669845, 0.142992, 0.582154)
      rgb=(0.674522, 0.147419, 0.578688)
      rgb=(0.679160, 0.151848, 0.575189)
      rgb=(0.683758, 0.156278, 0.571660)
      rgb=(0.688318, 0.160709, 0.568103)
      rgb=(0.692840, 0.165141, 0.564522)
      rgb=(0.697324, 0.169573, 0.560919)
      rgb=(0.701769, 0.174005, 0.557296)
      rgb=(0.706178, 0.178437, 0.553657)
      rgb=(0.710549, 0.182868, 0.550004)
      rgb=(0.714883, 0.187299, 0.546338)
      rgb=(0.719181, 0.191729, 0.542663)
      rgb=(0.723444, 0.196158, 0.538981)
      rgb=(0.727670, 0.200586, 0.535293)
      rgb=(0.731862, 0.205013, 0.531601)
      rgb=(0.736019, 0.209439, 0.527908)
      rgb=(0.740143, 0.213864, 0.524216)
      rgb=(0.744232, 0.218288, 0.520524)
      rgb=(0.748289, 0.222711, 0.516834)
      rgb=(0.752312, 0.227133, 0.513149)
      rgb=(0.756304, 0.231555, 0.509468)
      rgb=(0.760264, 0.235976, 0.505794)
      rgb=(0.764193, 0.240396, 0.502126)
      rgb=(0.768090, 0.244817, 0.498465)
      rgb=(0.771958, 0.249237, 0.494813)
      rgb=(0.775796, 0.253658, 0.491171)
      rgb=(0.779604, 0.258078, 0.487539)
      rgb=(0.783383, 0.262500, 0.483918)
      rgb=(0.787133, 0.266922, 0.480307)
      rgb=(0.790855, 0.271345, 0.476706)
      rgb=(0.794549, 0.275770, 0.473117)
      rgb=(0.798216, 0.280197, 0.469538)
      rgb=(0.801855, 0.284626, 0.465971)
      rgb=(0.805467, 0.289057, 0.462415)
      rgb=(0.809052, 0.293491, 0.458870)
      rgb=(0.812612, 0.297928, 0.455338)
      rgb=(0.816144, 0.302368, 0.451816)
      rgb=(0.819651, 0.306812, 0.448306)
      rgb=(0.823132, 0.311261, 0.444806)
      rgb=(0.826588, 0.315714, 0.441316)
      rgb=(0.830018, 0.320172, 0.437836)
      rgb=(0.833422, 0.324635, 0.434366)
      rgb=(0.836801, 0.329105, 0.430905)
      rgb=(0.840155, 0.333580, 0.427455)
      rgb=(0.843484, 0.338062, 0.424013)
      rgb=(0.846788, 0.342551, 0.420579)
      rgb=(0.850066, 0.347048, 0.417153)
      rgb=(0.853319, 0.351553, 0.413734)
      rgb=(0.856547, 0.356066, 0.410322)
      rgb=(0.859750, 0.360588, 0.406917)
      rgb=(0.862927, 0.365119, 0.403519)
      rgb=(0.866078, 0.369660, 0.400126)
      rgb=(0.869203, 0.374212, 0.396738)
      rgb=(0.872303, 0.378774, 0.393355)
      rgb=(0.875376, 0.383347, 0.389976)
      rgb=(0.878423, 0.387932, 0.386600)
      rgb=(0.881443, 0.392529, 0.383229)
      rgb=(0.884436, 0.397139, 0.379860)
      rgb=(0.887402, 0.401762, 0.376494)
      rgb=(0.890340, 0.406398, 0.373130)
      rgb=(0.893250, 0.411048, 0.369768)
      rgb=(0.896131, 0.415712, 0.366407)
      rgb=(0.898984, 0.420392, 0.363047)
      rgb=(0.901807, 0.425087, 0.359688)
      rgb=(0.904601, 0.429797, 0.356329)
      rgb=(0.907365, 0.434524, 0.352970)
      rgb=(0.910098, 0.439268, 0.349610)
      rgb=(0.912800, 0.444029, 0.346251)
      rgb=(0.915471, 0.448807, 0.342890)
      rgb=(0.918109, 0.453603, 0.339529)
      rgb=(0.920714, 0.458417, 0.336166)
      rgb=(0.923287, 0.463251, 0.332801)
      rgb=(0.925825, 0.468103, 0.329435)
      rgb=(0.928329, 0.472975, 0.326067)
      rgb=(0.930798, 0.477867, 0.322697)
      rgb=(0.933232, 0.482780, 0.319325)
      rgb=(0.935630, 0.487712, 0.315952)
      rgb=(0.937990, 0.492667, 0.312575)
      rgb=(0.940313, 0.497642, 0.309197)
      rgb=(0.942598, 0.502639, 0.305816)
      rgb=(0.944844, 0.507658, 0.302433)
      rgb=(0.947051, 0.512699, 0.299049)
      rgb=(0.949217, 0.517763, 0.295662)
      rgb=(0.951344, 0.522850, 0.292275)
      rgb=(0.953428, 0.527960, 0.288883)
      rgb=(0.955470, 0.533093, 0.285490)
      rgb=(0.957469, 0.538250, 0.282096)
      rgb=(0.959424, 0.543431, 0.278701)
      rgb=(0.961336, 0.548636, 0.275305)
      rgb=(0.963203, 0.553865, 0.271909)
      rgb=(0.965024, 0.559118, 0.268513)
      rgb=(0.966798, 0.564396, 0.265118)
      rgb=(0.968526, 0.569700, 0.261721)
      rgb=(0.970205, 0.575028, 0.258325)
      rgb=(0.971835, 0.580382, 0.254931)
      rgb=(0.973416, 0.585761, 0.251540)
      rgb=(0.974947, 0.591165, 0.248151)
      rgb=(0.976428, 0.596595, 0.244767)
      rgb=(0.977856, 0.602051, 0.241387)
      rgb=(0.979233, 0.607532, 0.238013)
      rgb=(0.980556, 0.613039, 0.234646)
      rgb=(0.981826, 0.618572, 0.231287)
      rgb=(0.983041, 0.624131, 0.227937)
      rgb=(0.984199, 0.629718, 0.224595)
      rgb=(0.985301, 0.635330, 0.221265)
      rgb=(0.986345, 0.640969, 0.217948)
      rgb=(0.987332, 0.646633, 0.214648)
      rgb=(0.988260, 0.652325, 0.211364)
      rgb=(0.989128, 0.658043, 0.208100)
      rgb=(0.989935, 0.663787, 0.204859)
      rgb=(0.990681, 0.669558, 0.201642)
      rgb=(0.991365, 0.675355, 0.198453)
      rgb=(0.991985, 0.681179, 0.195295)
      rgb=(0.992541, 0.687030, 0.192170)
      rgb=(0.993032, 0.692907, 0.189084)
      rgb=(0.993456, 0.698810, 0.186041)
      rgb=(0.993814, 0.704741, 0.183043)
      rgb=(0.994103, 0.710698, 0.180097)
      rgb=(0.994324, 0.716681, 0.177208)
      rgb=(0.994474, 0.722691, 0.174381)
      rgb=(0.994553, 0.728728, 0.171622)
      rgb=(0.994561, 0.734791, 0.168938)
      rgb=(0.994495, 0.740880, 0.166335)
      rgb=(0.994355, 0.746995, 0.163821)
      rgb=(0.994141, 0.753137, 0.161404)
      rgb=(0.993851, 0.759304, 0.159092)
      rgb=(0.993482, 0.765499, 0.156891)
      rgb=(0.993033, 0.771720, 0.154808)
      rgb=(0.992505, 0.777967, 0.152855)
      rgb=(0.991897, 0.784239, 0.151042)
      rgb=(0.991209, 0.790537, 0.149377)
      rgb=(0.990439, 0.796859, 0.147870)
      rgb=(0.989587, 0.803205, 0.146529)
      rgb=(0.988648, 0.809579, 0.145357)
      rgb=(0.987621, 0.815978, 0.144363)
      rgb=(0.986509, 0.822401, 0.143557)
      rgb=(0.985314, 0.828846, 0.142945)
      rgb=(0.984031, 0.835315, 0.142528)
      rgb=(0.982653, 0.841812, 0.142303)
      rgb=(0.981190, 0.848329, 0.142279)
      rgb=(0.979644, 0.854866, 0.142453)
      rgb=(0.977995, 0.861432, 0.142808)
      rgb=(0.976265, 0.868016, 0.143351)
      rgb=(0.974443, 0.874622, 0.144061)
      rgb=(0.972530, 0.881250, 0.144923)
      rgb=(0.970533, 0.887896, 0.145919)
      rgb=(0.968443, 0.894564, 0.147014)
      rgb=(0.966271, 0.901249, 0.148180)
      rgb=(0.964021, 0.907950, 0.149370)
      rgb=(0.961681, 0.914672, 0.150520)
      rgb=(0.959276, 0.921407, 0.151566)
      rgb=(0.956808, 0.928152, 0.152409)
      rgb=(0.954287, 0.934908, 0.152921)
      rgb=(0.951726, 0.941671, 0.152925)
      rgb=(0.949151, 0.948435, 0.152178)
      rgb=(0.946602, 0.955190, 0.150328)
      rgb=(0.944152, 0.961916, 0.146861)
      rgb=(0.941896, 0.968590, 0.140956)
      rgb=(0.940015, 0.975158, 0.131326)
    }
}

\newcommand{\logLogSlopeTriangle}[5]
{

    \pgfplotsextra
    {
        \pgfkeysgetvalue{/pgfplots/xmin}{\xmin}
        \pgfkeysgetvalue{/pgfplots/xmax}{\xmax}
        \pgfkeysgetvalue{/pgfplots/ymin}{\ymin}
        \pgfkeysgetvalue{/pgfplots/ymax}{\ymax}

        \pgfmathsetmacro{\xArel}{#1}
        \pgfmathsetmacro{\yArel}{#3}
        \pgfmathsetmacro{\xBrel}{#1-#2}
        \pgfmathsetmacro{\yBrel}{\yArel}
        \pgfmathsetmacro{\xCrel}{\xArel}

        \pgfmathsetmacro{\lnxB}{\xmin*(1-(#1-#2))+\xmax*(#1-#2)} 
        \pgfmathsetmacro{\lnxA}{\xmin*(1-#1)+\xmax*#1} 
        \pgfmathsetmacro{\lnyA}{\ymin*(1-#3)+\ymax*#3} 
        \pgfmathsetmacro{\lnyC}{\lnyA+#4*(\lnxA-\lnxB)}
        \pgfmathsetmacro{\yCrel}{\lnyC-\ymin)/(\ymax-\ymin)} 

        \coordinate (A) at (rel axis cs:\xArel,\yArel);
        \coordinate (B) at (rel axis cs:\xBrel,\yBrel);
        \coordinate (C) at (rel axis cs:\xCrel,\yCrel);

        \draw[#5]   (A)-- node[pos=0.5,anchor=north] {1}
                    (B)--
                    (C)-- node[pos=0.5,anchor=west] {#4}
                    cycle;
    }
}

\pgfplotsset{compat=newest}

\declaretheoremstyle[notefont=\bfseries,notebraces={}{},%
    headpunct={},postheadspace=1em]{mystyle}
\declaretheorem[style=mystyle,numbered=no,name=Theorem]{thm-hand}

\declaretheoremstyle[notefont=\bfseries,notebraces={}{},%
    headpunct={},postheadspace=1em]{mystyle}
\declaretheorem[style=mystyle,numbered=no,name=Proposition]{prop-hand}

\declaretheoremstyle[notefont=\bfseries,notebraces={}{},%
    headpunct={},postheadspace=1em]{mystyle}
\declaretheorem[style=mystyle,numbered=no,name=Corollary]{Cor-hand} 

\declaretheoremstyle[notefont=\bfseries,notebraces={}{},%
    headpunct={},postheadspace=1em]{mystyle}
\declaretheorem[style=mystyle,numbered=no,name=Lemma]{lem-hand}

\declaretheoremstyle[notefont=\bfseries,notebraces={}{},%
    headpunct={},postheadspace=1em]{mystyle}
\declaretheorem[style=mystyle,numbered=no,name=Remark]{rem-hand}

\usepackage{longtable}
\usepackage{booktabs}

\makeatletter
\renewcommand*\env@matrix[1][\arraystretch]{%
  \edef\arraystretch{#1}%
  \hskip -\arraycolsep
  \let\@ifnextchar\new@ifnextchar
  \array{*\c@MaxMatrixCols c}}
\makeatother

\renewcommand{\thefigure}{\arabic{section}.\arabic{figure}}

\allowdisplaybreaks

\newcommand{\avg}[1]{\ensuremath{\{\!\!\{#1\}\!\!\}} }
\newcommand{\jump}[1]{\ensuremath{[\![#1]\!]} }

\newcommand {\R}{\mathbb{R}}

\newcommand{\matx}[1]{\mathcal{\mathcal{#1}}} 
\newcommand{\blockmatx}[1]{\boldsymbol{\mathcal{#1}}}
\newcommand{\blockvec}[1]{\overset{\leftrightarrow}{#1}}
\newcommand{\blockveccon}[1]{\overset{\leftrightarrow}{\tilde{#1}}}

\newcommand{\interpolation}[2]{\mathbb{I}_{#1}\!#2}
\newcommand{\Dprojection}[2]{\vec{\mathbb{D}}_{#1}\!#2}

\newcommand{\pderivative}[2]{\frac{\partial #1}{\partial #2}}

\newcommand{\etal}{et al.~}

\begin{document}


\title{\Large Split form ALE discontinuous Galerkin methods with applications to under-resolved turbulent low-Mach number flows}

\setcounter{footnote}{0}
\font\myfont=cmr12 at 10pt 
\author{{\myfont Nico Krais\footnote{\ Inst. of Aerodynamics and Gas Dynamics (IAG), University of Stuttgart, Email: krais@iag.uni-stuttgart.de}, \ 
Gero Schnücke\footnote{\ Dep. for Mathematics and Computer Science, University of Cologne, Email: gschnuec@uni-koeln.de},  
\ Thomas Bolemann\footnote{\ Inst. of Aerodynamics and Gas Dynamics (IAG), University of Stuttgart, Email: bolemann@iag.uni-stuttgart.de} \ and     
Gregor J. Gassner \footnote{\ Dep. for Mathematics/Computer Science; Center for Data and Simulation Science, University of Cologne, Email: ggassner@uni-koeln.de}}} 
\pagestyle{fancy}
\lhead{N. Krais, G. Schn\"ucke,  T. Bolemann and G. J. Gassner}
\rhead{Split form ALE DG methods for turbulent flows}
\date {\myfont \today}
\maketitle
{\hrule height 1.5pt} 
\begin{abstract}
\noindent
The construction of discontinuous Galerkin (DG) methods for the compressible 
Euler or Navier-Stokes equations (NSE) includes the approximation of non-linear 
flux terms in the volume integrals. The terms can lead to aliasing and stability 
issues in turbulence simulations with moderate Mach numbers ($\textrm{Ma}
\lesssim 0.3$), e.g. due to under-resolution of vortical dominated structures 
typical in large eddy simulations (LES). The kinetic energy or entropy  are 
elevated in smooth, but under-resolved parts of the solution which are affected 
by aliasing. It is known that the kinetic energy is not a conserved quantity for 
compressible flows, but for small Mach numbers minor deviations from a conserved 
evolution can be expected. While it is formally possible to construct kinetic energy  
preserving (KEP) and entropy conserving (EC) DG methods for the Euler equations, due 
to the viscous terms in case of the NSE, we aim to construct kinetic energy 
dissipative (KED) or entropy stable (ES) DG methods on moving curved hexahedral 
meshes. The Arbitrary Lagrangian-Eulerian (ALE) approach is used to include the 
effect of mesh motion in the split form DG methods. First, we use the three 
dimensional Taylor-Green vortex to investigate and analyze our theoretical findings 
and the behavior of the novel split form ALE DG schemes for a turbulent vortical 
dominated flow. Second, we apply the framework to a complex aerodynamics application. 
An implicit LES split form ALE DG approach is used to simulate the transitional flow 
around a plunging SD7003 airfoil at Reynolds number $\textrm{Re}=40,000$ and Mach 
number $\textrm{Ma}=0.1$. We compare the standard nodal ALE DG scheme, the ALE DG 
variant with consistent overintegration of the non-linear terms and the novel KED and 
ES split form ALE DG methods in terms of robustness, accuracy and computational 
efficiency. 
{\hrule height 1.5pt}
\end{abstract}
\newpage

\newpage
\section{Introduction}\label{sec:Int}
The numerical simulation of under-resolved turbulent flows in the regime of low or 
moderate Mach numbers (the local Mach number is $\textrm{Ma}\lesssim 0.3$) 
requires accurate dispersion and dissipation behavior 
\cite{gassner2013accuracy, moura2015linear,nigro2019low}. It is desirable that the 
dissipation errors are very low for well resolved scales and are very high for 
scales close to the Nyquist cutoff, to get rid of small scale noise. In addition, 
the dispersion error should be small for a wide range of scales. This motivates the 
application of high order methods, e.g. discontinuous Galerkin (DG) methods 
\cite{Cockburn2001,hesthaven2007nodal,Kopriva2009}, compact finite 
difference methods \cite{lele1992compact} or summation-by-parts finite 
difference methods \cite{nordstrom1999boundary}, to simulate the transition to 
turbulence and fully turbulent compressible flows. 

DG methods are a class of finite element methods using piecewise polynomials as 
basis functions and  Riemann solver based interface numerical flux functions 
along the element interfaces. The numerical dissipation caused by the interface 
Riemann solver acts as a filter of high frequency solution components. 
Furthermore, high order DG methods have excellent dispersion behavior 
\cite{gassner2011comparison, hesthaven2007nodal}. These observations motivate 
the application of DG methods in implicit large eddy simulations (iLES). For 
instance in \cite{de2012dns, flad2016simulation, flad2017use, moura2017eddy,  
uranga2011implicit, stoellinger2019dynamic}, DG iLES were used to simulate 
flows at moderate Reynolds numbers. In all these simulations no additional explicit 
sub-grid-scale models were used, since the high frequency type dissipation of the 
numerical interface flux functions is interpreted as a dissipative implicit 
sub-grid-scale model.

In many cases, the iLES DG simulations are negatively affected by aliasing errors 
due to the strong non-linearity of the flux functions \cite{blaisdell1996effect, 
kirby2003aliasing}. These errors are generated in the volume integrals and cannot be 
sufficiently controlled with the dissipation of the numerical surface fluxes alone. 
A careful treatment of aliasing is however necessary, as these errors may even cause 
fatal failure of the simulation. Standard de-aliasing techniques in the DG framework 
are based on either projection of the non-linear flux functions, e.g., 
\cite{flad2017use, winters2018comparative} or exact evaluations of the integrals in 
the variational formulation, sometimes termed "overintegration" or "consistent 
integration", e.g., \cite{gassner2013accuracy, 
kirby2003aliasing,kopriva2018stability,mengaldo2015dealiasing}.

The nodal DG spectral element method (DGSEM), e.g.,  \cite{Kopriva2009}, is 
constructed with local tensor-product Lagrange polynomial basis functions computed 
from Legendre-Gauss-Lobatto (LGL) points. The collocation of interpolation and 
quadrature nodes is used in the spatial discretization. This approach ensures 
that the derivative matrix in the DGSEM  provides a summation-by-parts (SBP)  
operator \cite{Gassner2013}. A SBP operator gives a discrete analogue of the 
integration-by-parts formula \cite{Fernandez2014,Gassner2013,kreiss1} and thus 
ideas from the continuous stability analysis can be mimicked at the discrete 
level. The SBP property is a powerful tool to construct the numerical 
approximation in a way that aliasing issues are avoided. 

The equations of gas dynamics, e.g. the Euler and Navier-Stokes equations
(NSE), are equipped with a mathematical entropy equal to the scaled negative
thermodynamic entropy. Hence, the mathematical model correctly captures the
second law of thermodynamics \cite{Barth2018}. A numerical scheme for these
equations should reflect the properties of the entropy on the discrete level.
In \cite{Tadmor2003} Tadmor gave a discrete entropy criteria to construct
entropy conservative (EC) Riemann solver based interface two point flux
functions for first oder finite volume (FV) and finite difference (FD) methods
to solve systems of conservation laws. This criteria was used to construct EC
two point flux functions for the Euler equations in \cite{Chandrashekar2013,
ismail2009affordable, Ranocha2018}. These low order EC methods can be modified
by adding dissipation to the numerical fluxes such that the entropy is
decreasing for all times. Then the method becomes entropy stable (ES) such that
the entropy of the solution is bounded. The kinetic energy is another important
physical quantity for the dynamics of turbulent flows. Therefore, a numerical scheme 
should also accurately capture its evolution and the transfer of kinetic and inner 
energy. Jameson \cite{Jameson2008} gave a discrete criteria to construct kinetic 
energy preserving (KEP) two point flux functions for first oder FV/FD methods to 
solve the Euler equations. However, in under-resolved turbulent parts of the 
simulation, it might be necessary to increase the robustness of the
discretization by adding dissipation to the KEP numerical fluxes. Here a
suitable dissipation term should be chosen such that the kinetic energy is
guaranteed to be dissipated by the scheme (no pile up of energy in higher
frequencies). In the current work, methods with this property will be denoted
kinetic energy dissipative (KED) methods. It should be noted that a method for
the compressible NSE can only be KDE or ES and not KEP or EC, since the model
contains viscous and dissipative terms in the momentum and energy evolution. 

Gassner \cite{gassner2014kinetic} used the SBP property to 
construct a KEP DGSEM for the one dimensional Euler equations. This method is a 
discrete analogue of a quasilinear skew-symmetric formulation for the Euler 
equations and can be written as the standard DGSEM with an additional source 
term that acts as a de-aliasing mechanism. It is important to note that the 
skew-symmetric formulation is still fully conservative and satisfies the Lax-
Wendroff theorem. For under-resolved computations, the source term is active and 
accounts for the error of the discrete product rule. On the other hand in smooth 
parts of the solution without large variation the term almost vanishes to zero, 
since the 
approximation is consistent. The stabilizing effect of high order conservative 
skew-symmetric or skew-symmetric-like schemes in under-resolved computations 
was also observed in other publications, e.g. \cite{ducros2000high, 
Kennedy2008, kuya2018kinetic,  morinishi2010skew, Pirozzoli2010}. However, 
since the thermodynamic entropy is a logarithmic function depending on density 
and pressure for the Euler equations and NSE, it is difficult to find an explicit 
discrete skew-symmetric or skew-symmetric-like formulation. 

A framework to construct high order EC schemes in periodic domains has been 
given by LeFloch et al. \cite{lefloch2002fully}. Fisher and Carpenter 
\cite{Fisher2013} combined this approach with SBP operators and proved that 
two-point EC fluxes can be used to construct high order schemes when the 
derivative approximations in space are SBP operators. Gassner et al. 
\cite{Gassner2016, Gassner2017} showed that skew-symmetric-like (split form) 
DGSEM formulations for the Euler equations can be discretely recovered when  
specific numerical volume fluxes in the flux form volume integral of Fisher and 
Carpenter are chosen. In particular, the three dimensional version of the skew-symmetric
KEP DGSEM in \cite{gassner2014kinetic} can be recovered by using the 
two point flux from Morinishi et al. \cite{morinishi2010skew} as volume flux. 
Hence, ES DGSEM can be constructed when the derivative matrix satisfies the SBP 
property and a two-point EC flux function is used as volume flux in the flux 
form volume integral of Fisher and Carpenter. The construction of KED DGSEM 
requires more attention, since it has been shown numerically that there are two 
point flux functions, e.g. \cite{Chandrashekar2013}, which are KEP in the sense 
of Jameson \cite{Jameson2008}, but when they are used as volume flux in the 
DGSEM the kinetic energy is necessarily not preserved (see  \cite{Gassner2016, 
Ranocha2018} or Figure~\ref{fig:TGV_EC_KEP_Euler} in Section \ref{EC_KEP_Test} of 
this work). 

It should be mentioned that the described approach and the associated analysis  
is for semi-discrete ES or KED high order DGSEM. Moreover, the construction of 
these methods rely on the assumption that the discrete density and pressure are 
positive. In general, this migh not always be the case for high order methods, 
since the positivity of the discrete density and pressure can be affected by 
spurious oscillations in the numerical solution. This phenomenon is for 
instance observed for strong shock waves which appear in the regime of high 
Mach numbers, or severe under-resolution of strong shear-layers.

Another important component for the simulation of turbulent flows are adaptive 
discretizations, where resolution is increased in regions with large spatial 
variations. The r-adaptive method involves the re-distribution of the mesh 
nodes in regions of rapid variation of the solution \cite{tang2005moving}. In 
comparison with h-adaptive discretizations, where the mesh is refined and coarsened 
by changing the number of elements in the tessellation, the r-adaptive method has 
some advantages, e.g. no hanging nodes appear and the number of elements does not 
change. On the other hand a r-adaptive method can be only used when the effect of 
mesh movement is appropriately accounted for the discretization of the system of 
conservation laws. This can be done by an Arbitrary Lagrangian-Eulerian (ALE) 
approach \cite{donea2017arbitrary}. In the last decades several ALE DG methods have 
been developed, e.g. \cite{kopriva2016provably, lomtev1999discontinuous, 
minoli2011discontinuous, nguyen2010arbitrary, persson2009discontinuous, 
winters2014ale}. Furthermore, in \cite{Schnuecke2018} a provably ES moving mesh 
ALE DGSEM for the three dimensional Euler equations on curved hexahedral 
elements and in \cite{yamaleev2019entropy} a moving mesh ES spectral collocation 
scheme for the three dimensional NSE were constructed.   

The structure of the present work is as follows: First, the approaches 
\cite{Gassner2017, Schnuecke2018} are used to construct a provably ES ALE DGSEM 
for the three dimensional NSE. Then a provable KEP ALE DGSEM for the Euler 
equations on curved moving elements is constructed. The discrete energy balance of 
the discrete kinetic, internal and total energy is analyzed for 
the split form ALE DGSEM  when different well known split forms \cite{Chandrashekar2013, Kennedy2008, kuya2018kinetic, 
Pirozzoli2010, Ranocha2018} are used. Afterward, the constructed KEP DGSEM is 
combined with the approach in \cite{Gassner2017} to construct a KED ALE DGSEM for 
the three dimensional NSE. 
The theoretical findings and the constructed methods are numerically 
investigated for the Taylor-Green Vortex (TGV) problem \cite{shu2005numerical}. 
Finally, the constructed ALE DGSEM is applied in an iLES for a complex aerodynamics 
application. An interesting area of application for problems with moving boundaries 
in aerodynamics is the moderate Reynolds number flow around pitching 
airfoils. These flows can serve as a model for the wing motion of natural 
flyers or the situation in micro unmanned vehicles. As a specific example, the 
transitional flow around a plunging SD7003 airfoil at Reynolds number $
\textrm{Re}=40,000$ and Mach number $\textrm{Ma}=0.1$ is simulated. This setup, 
among others, was investigated by Visbal \cite{visbal2009high}, which gives the 
possibility for a comparison to available experimental measurements and 
numerical simulation results from literature \cite{galbraith2010implicit, 
mcgowan2008computation, ol2009computation, radespiel2007numerical, 
yuan2005investigation}.

\section{The Navier Stokes equations (NSE) in three dimensions}\label{NSE} 
We apply the block vector notation as in \cite{Gassner2017} to present the three dimensional NSE. In Appendix \ref{Nomenclature:BlockVector} the key elements of this notation are briefly summarized. The three dimensional NSE in compact block vector form are given by 
\begin{equation}\label{eq:NSELawBlock}
\pderivative{\textbf{u}}{t}+\vec{\nabla}_{x}\cdot\blockvec{\textbf{f}}\left(\textbf{u}\right)=\vec{\nabla}_{x}\cdot\blockvec{\textbf{f}}{}^{\text{v}}\left(\textbf{u},\vec{\nabla}_{x}\textbf{u}\right).
\end{equation}
These equations are considered on a time-dependent domain $\Omega\left(t\right)\subseteq \R^3$ with suitable initial data and boundary conditions. The state vector is given by    
\begin{equation}\label{Statevector}
\mathbf{u}=\left[\rho,\rho\vec{u},\rho e+\frac{1}{2}\rho\left|\vec{u}\right|^{2}\right]^{T},
\end{equation}
where $\rho$ is the mass density, $\vec{u}=\left[u_{1},u_{2},u_{3}\right]^{T}$ the velocity, $p$ is the pressure and the specific internal energy is     
\begin{equation}\label{SpecificInnerEnergy}
e=\frac{p}{\left(\gamma-1\right)\rho}, \qquad \gamma >1. 
\end{equation}
The advective block flux $\blockvec{\textbf{f}}=\left[\textbf{f}_{1}^{T},\textbf{f}_{2}^{T},\textbf{f}_{3}^{T}\right]^{T}$ consists of the flux vectors for $\iota=1,2,3$ given by  
\begin{align}\label{AdvectiveFluxes}
\begin{split}
\mathbf{f}_{\iota}^{1\phantom{+1}}=& \phantom{+} \rho u_{\iota}, \\
\mathbf{f}_{\iota}^{\upsilon+1}=& \phantom{+} \rho u_{\iota}u_{\upsilon}+p \delta_{\iota \upsilon},\qquad\upsilon=1,2,3, \\
\mathbf{f}_{\iota}^{5\phantom{+1}}=& \phantom{+} \rho u_{\iota}\left[e+\frac{1}{2}\left|\vec{u}\right|^{2}\right]+pu_{\iota},
\end{split}
\end{align}
where $\delta_{\iota \upsilon}$ is the Kronecker delta.
The viscous block flux $\blockvec{\textbf{f}}{}^{\text{v}}=\left[(\textbf{f}^{\,\text{v}}_{1})^{T},(\textbf{f}^{\,\text{v}}_{2})^{T},(\textbf{f}^{\,\text{v}}_{3})^{T}\right]^{T}$ consists of the flux vectors 
\begin{align}\label{viscousFlux}
\begin{split}
\textbf{f}_{1}^{\,\text{v}}=&\left[0,\sigma_{11},\sigma_{12},\sigma_{13},\sum_{j=1}^{3}u_{j}\sigma_{1j}-q_{1}\right]^{T}, \\
\textbf{f}_{2}^{\,\text{v}}=&\left[0,\sigma_{21},\sigma_{22},\sigma_{23},\sum_{j=1}^{3}u_{j}\sigma_{2j}-q_{2}\right]^{T}, \\
\textbf{f}_{3}^{\,\text{v}}=&\left[0,\sigma_{31},\sigma_{32},\sigma_{33},\sum_{j=1}^{3}u_{j}\sigma_{3j}-q_{3}\right]^{T},
\end{split}
\end{align}
where 
\begin{equation}\label{viscousTensor}
\sigma_{ij}=\mu\left(\pderivative{u_{j}}{x_{i}}+\pderivative{u_{i}}{x_{j}}\right)-\frac{2}{3}\mu\left(\vec{\nabla}_{x}\cdot\vec{u}\right)\delta_{ij}, 
\qquad 
q_{i}=-\kappa\pderivative{T}{x_{i}},
\qquad i,j=1,2,3,  
\end{equation}
and $T=\frac{p}{R\rho}$ is the temperature.\footnote{Here the Greek letter $\sigma$ is used to characterize the shear stress tensor.
The common notation $\tau$ is not used, since $\tau$ will be used as time variable on the reference element in the following
sections.} The introduced constants $\mu, \kappa, R>0$ describe the dynamic viscosity, thermal conductivity and the universal gas constant. Instead of defining the thermal conductivity and universal gas constant, we state the Prandtl number, which is determined by these quantities, i.e.   
\begin{equation}\label{PrandtlNumber}
\textrm{Pr}=\frac{\gamma\mu R}{\left(\gamma-1\right)\kappa}. 
\end{equation}
In order to classify the viscous effects of the system, we define the Reynolds number 
\begin{equation}\label{ReynoldsNumber}
\textrm{Re}=\frac{\rho \left|\vec{u}\right|L}{\mu},
\end{equation}
where $L$ is a characteristic length. 

\subsection{Entropy representation of the viscous fluxes}\label{ViscousFluxSymmetrization} 
The NSE are equipped with the entropy/entropy flux pair   
\begin{equation}\label{EntropyNSE}
s=-\frac{\rho\varsigma}{\gamma-1},\qquad 
f_{\iota}^{s}=-\frac{\rho\varsigma u_{\iota}}{\gamma-1},\qquad \iota=1,2,3,
\end{equation} 
where $\varsigma=\log\left(p\rho^{-\gamma}\right)$ is the thermodynamic
entropy of the fluid (cf. Harten \cite{Harten1983}). The entropy variables are given by  
\begin{align}\label{EntropyVariables}
\pderivative{s}{\textbf{u}}=\left[\frac{\gamma-\varsigma}{\gamma-1}-\frac{\rho}{2p}\left|\vec{u}\right|^{2},\frac{\rho\vec{u}}{p},-\frac{\rho}{p}\right]^{T}=:\mathbf{w}.
\end{align}
Hughes et al. \cite{Hughes1986} proved that the viscous flux in the NSE \eqref{eq:NSELawBlock} can be written as  
\begin{equation}\label{ViscousBlockFlux1}
\blockvec{\textbf{f}}{}^{\text{v}}\left(\textbf{u},\vec{\nabla}_{x}\textbf{u}\right)=\blockmatx{B}^{\,\text{v}}\left(\textbf{w}\right)\vec{\nabla}_{x}\textbf{w}\qquad\text{or}\qquad\textbf{f}_{\iota}^{\,\text{v}}=\sum_{\vartheta=1}^{3}\matx{B}_{\iota \vartheta}^{\,\text{v}}\left(\textbf{w}\right)\pderivative{\textbf{w}}{x_{j}},\qquad \iota=1,2,3,
\end{equation}
where $\blockmatx{B}^{\,\text{v}}\left(\textbf{w}\right)$ is a $15 \times 15$ block matrix which blocks satisfy $\matx{B}_{\iota \vartheta}^{v}\left(\textbf{w}\right)=\matx{B}_{\vartheta\iota}^{v}\left(\textbf{w}\right)^{T}$ for $\iota,\vartheta=1,2,3$. A representation of these blocks can be found in the appendix of Ray's PhD thesis \cite{Ray2017}. The matrix is positive semidefinite such that
\begin{equation}\label{ViscousBlockFlux3}
\mathbf{p}^{T}\left[\blockmatx{B}^{\,\text{v}}\left(\textbf{w}\right)\right]\mathbf{p}\geq0,\qquad\mathbf{p}\in\R^{15}.
\end{equation}
Note that only entropy functions of the type $-c_{0}\varsigma\rho+c_{1}$, with constants $c_{0}$, $c_{1}$ and $c_{0}>0$, provide a positive semidefinite viscous block matrix (cf. Hughes et al. \cite{Hughes1986}).


\subsection{Transformation of the NSE onto a reference element}\label{Sec:TransNSE}
We divide the domain $\Omega(t)$ in $K$ time-dependent non-overlapping curved elements $e_{\alpha}(t)$, $\alpha =1,\dots,K$. Each element $e_{\alpha}(t)$ is mapped onto a reference element $\mathrm{E}=\left[-1,1\right]^{3}$ by an isoparametric transformation 
\begin{equation}\label{IsoparametricTrans}
\vec{x}\left(t\right)=\vec{\chi}\left(\tau,\vec{\xi}\right),\qquad\vec{\xi}=\left(\xi^{1},\xi^{2},\xi^{3}\right)^{T}\in \mathrm{E},\qquad\tau\in\left[0,T\right].
\end{equation}
This mapping provides the covariant basis vectors   
\begin{equation}
\vec{a}_{i}:=\pderivative{\vec{\chi}}{\xi^{i}},\qquad i=1,2,3,
\end{equation}
and the volume weighted contravariant vectors  
\begin{equation}\label{ContravariantVectors}
J\vec{a}^{\,\iota}=\vec{a}_{\alpha}\times \vec{a}_{\beta}, \quad (\iota, \alpha, \beta) \ \text{cyclic}, \quad J:=\det\,\nabla\vec{\chi}, 
\end{equation}
where $\iota,\alpha,\beta\in\left\{ 1,2,3\right\}$ and $\nabla\vec{\chi}$ is the Jacobian matrix of the isoparametric transformation \eqref{IsoparametricTrans}. The contravariant vectors satisfy the metric identities    
\begin{equation}\label{MetricIdentities}
\sum_{\iota=1}^{3}\pderivative{J\vec{a}^{\,\iota}}{\xi^{\iota}}=0.
\end{equation}
In order to combine the contravariant vectors with the block vector nomenclature, the following block matrix has been introduced in \cite{Gassner2017}  
\begin{equation}\label{TransformationBlock}
\blockmatx{M}:=\left[\begin{matrix}J\vec{a}_{1}^{1}\matx{I_{\textit{5}}} & J\vec{a}_{1}^{2}\matx{I_{\textit{5}}} & J\vec{a}_{1}^{3}\matx{I_{\textit{5}}}\\[0.1cm]
J\vec{a}_{2}^{1}\matx{I_{\textit{5}}} & J\vec{a}_{2}^{2}\matx{I_{\textit{5}}} & J\vec{a}_{2}^{3}\matx{I_{\textit{5}}}\\[0.1cm]
J\vec{a}_{3}^{1}\matx{I_{\textit{5}}} & J\vec{a}_{3}^{1}\matx{I_{\textit{5}}} & J\vec{a}_{3}^{3}\matx{I_{\textit{5}}}
\end{matrix}\right],
\end{equation}
where the matrix $\matx{I_{\textit{5}}}$ is the $5\times 5$ identity matrix and $Ja_{j}^{i}$ is the component of $J \vec{a}^{i}$ in
the $j$-th Cartesian coordinate direction. The contravariant vectors and the block matrix \eqref{TransformationBlock} enable the
transformation of the gradient and the divergence on the reference space. The transformation formulas for these differential
operators are given in Appendix \ref{TransformationDifferentialOperator}. Analogous to \cite{Schnuecke2018}, equation
\eqref{eq:NSELawBlock} in a time-dependent element $e_\alpha(t)$ turns into
\begin{subequations}\label{eq:NSELawBlockRefE}
\begin{align}
 \pderivative{J}{\tau}=&\vec{\nabla}_{\xi}\cdot\vec{\tilde{\nu}}, \label{GCL} \\
 \pderivative{\left(J\textbf{u}\right)}{\tau}
+\vec{\nabla}_{\xi}\cdot\blockvec{\tilde{\textbf{g}}}=&\vec{\nabla}_{\xi}\cdot\left(\blockmatx{\tilde{B}}^{\,\text{v}}\left(\textbf{w}\right)\blockvec{\textbf{q}}\right), \label{eq:NSE} \\
 J\blockvec{\textbf{q}}=&\blockmatx{M}\,\vec{\nabla}_{\xi}\textbf{w},
\label{eq:ViscousPart}
\end{align}
\end{subequations}
in the reference element $\mathrm{E}$, where $\vec{\nu}=\begin{bmatrix} \nu_{1},\nu_{2},\nu_{3} \end{bmatrix}^{T}$ is the grid velocity field and the contravariant flux vectors are given by
\begin{equation}\label{G-flux}
\vec{\tilde{\nu}}=\begin{bmatrix}J\vec{a}^{1}\cdot\vec{\nu}\\
J\vec{a}^{2}\cdot\vec{\nu}\\
J\vec{a}^{3}\cdot\vec{\nu}
\end{bmatrix},\qquad\blockveccon{\textbf{g}}=\blockveccon{\textbf{f}}-\vec{\tilde{\nu}}\textbf{u}=\begin{bmatrix}J\vec{a}^{1}\cdot\blockvec{\textbf{f}}-\left(J\vec{a}^{1}\cdot\vec{\nu}\right)\textbf{u}\\
J\vec{a}^{2}\cdot\blockvec{\textbf{f}}-\left(J\vec{a}^{2}\cdot\vec{\nu}\right)\textbf{u}\\
J\vec{a}^{3}\cdot\blockvec{\textbf{f}}-\left(J\vec{a}^{3}\cdot\vec{\nu}\right)\textbf{u}
\end{bmatrix},\qquad\blockmatx{\tilde{B}}^{\,\text{v}}\left(\textbf{w}\right)=\blockmatx{M}^{T}\blockmatx{B}^{\,\text{v}}\left(\textbf{w}\right),
\end{equation}  
with the viscous block matrix $\blockmatx{B}^{\,\text{v}}\left(\textbf{w}\right)$ introduced in section \ref{ViscousFluxSymmetrization} to represent the viscous fluxes. 

\subsection{Analysis of the kinetic energy}\label{KineticEnergyStability}
The kinetic energy is 
\begin{equation}\label{KineticEnergy}
\Bbbk:=\frac{1}{2}\rho\left|\vec{u}\right|^{2}=-\frac{1}{2}\rho\left|\vec{u}\right|^{2}+\rho\vec{u}^{T}\vec{u}=\textbf{v}^{T}\textbf{u},
\qquad 
\textbf{v}:=\left[-\frac{1}{2}\left|\vec{u}\right|^2,\vec{u},0\right]^T.
\end{equation} 
We note that $\textbf{v}^T=\pderivative{\Bbbk}{\textbf{u}}$. Thus, it follows for smooth solutions of the equation \eqref{eq:NSELawBlock} 
\begin{equation}\label{ContinuousKES1}
\mathbf{v}^{T}\pderivative{\left(J\textbf{u}\right)}{\tau}
=\mathbf{v}^{T}\textbf{u}\pderivative{J}{\tau}
+J\left(\mathbf{v}^{T}\pderivative{\,\mathbf{u}}{\tau}\right)
=\Bbbk\pderivative{J}{\tau}+J\pderivative{\Bbbk}{\tau}
=\pderivative{\left(J\Bbbk\right)}{\tau},
\end{equation}
and for the product with the flux divergence
\begin{align}\label{ContinuousKES2}
\begin{split}
\textbf{v}^{T}\vec{\nabla}_{\xi}\cdot\blockvec{\tilde{\textbf{g}}} 
=& J\mathbf{v}^{T}\vec{\nabla}_{x}\cdot\blockvec{\mathbf{g}}=J\vec{\nabla}_{x}\cdot\left[\left(\vec{u}-\vec{\nu}\right)\Bbbk\right]+J\left[\vec{\nabla}_{x}p\right]\cdot\vec{u}=\vec{\nabla}_{\xi}\cdot\left[\left(\vec{\tilde{u}}-\vec{\tilde{\nu}}\right)\Bbbk\right]+\left[\vec{\nabla}_{\xi}p\right]\cdot\vec{\tilde{u}}.
\end{split}
\end{align}

\subsubsection{Kinetic energy preservation for the Euler equations}\label{KineticEnergyEuler}
These identities provide the kinetic energy balance for the Euler equations. Therefore, we set the viscous fluxes to zero, multiply the equation \eqref{eq:NSE} with $\textbf{v}$, integrate over the reference element $\mathrm{E}$ and obtain by \eqref{ContinuousKES1}, \eqref{ContinuousKES2}    
 \begin{align}\label{ContinuousKESEuler}
\pderivative{}{\tau} \int_{\mathrm{E}}J\Bbbk\,d\vec{\xi}=-\int_{\partial\mathrm{E}}\left[\left(\vec{\tilde{u}}-\vec{\tilde{\nu}}\right)\Bbbk\right]^{T}\hat{n}\,dS-\int_{\mathrm{E}}\left[\vec{\nabla}_{\xi}p\right]\cdot\vec{\tilde{u}}\,d\vec{\xi}.
\end{align}
We note that the pressure work is for compressible flows in general a non-conservative term, but the advection terms can be rewritten as surface flux  integrals. A numerical scheme that mimics this behavior is called kinetic energy preserving. Comparing the total energy conservation with the kinetic energy balance gives a balance equation for the internal energy    
\begin{align}\label{ContinuousInnerEuler}
\begin{split}
\pderivative{}{\tau}
\int_{\mathrm{E}}J\rho e\,d\vec{\xi} 
=& \phantom{-} \pderivative{}{\tau}\int_{\mathrm{E}}J\left(\rho e+\frac{1}{2}\rho\left|\vec{u}\right|^{2}\right)\,d\vec{\xi}-
\pderivative{}{\tau}\int_{\mathrm{E}}J\Bbbk\,d\vec{\xi} \\ 
=& -\int_{\partial\mathrm{E}}e\rho\left(\vec{\tilde{u}}-\vec{\tilde{\nu}}\right)^{T}\hat{n}\,dS-\int_{\mathrm{E}}p\vec{\nabla}_{\xi}\cdot\vec{\tilde{u}}\,d\vec{\xi}.
\end{split}
\end{align}

\subsubsection{Kinetic energy dissipation for the NSE}\label{KineticEnergy StabilityNSE}
For the NSE the temporal evolution of the kinetic energy is complemented with a guaranteed kinetic energy dissipating contribution. We multiply the equation \eqref{eq:NSE} with $\textbf{v}$, integrate over $\mathrm{E}$ and obtain by \eqref{ContinuousKES1}, \eqref{ContinuousKES2}    
 \begin{align}\label{ContinuousKES3}
 \begin{split}
\pderivative{}{\tau}\int_{\mathrm{E}}J\Bbbk\,d\vec{\xi} 
=&-\int_{\partial\mathrm{E}}\left[\left(\vec{\tilde{u}}-\vec{\tilde{\nu}}\right)\Bbbk\right]^{T}\hat{n}\,dS-\int_{\mathrm{E}}\left[\vec{\nabla}_{\xi}p\right]\cdot\vec{\tilde{u}}\,d\vec{\xi} \\
 &-\int_{\mathrm{E}}\left[\vec{\nabla}_{\xi}\mathbf{v}\right]^{T}\left[\blockmatx{\tilde{B}}^{\,\text{v}}\left(\textbf{w}\right)\right]\,\left[\blockvec{\textbf{q}}\right]\,d\vec{\xi}+\int_{\partial\mathrm{E}}\textbf{v}^{T}\left\{ \left[\blockmatx{\tilde{B}}^{\,\text{v}}\left(\textbf{w}\right)\right]\,\blockvec{\textbf{q}}\cdot\hat{n}\right\} \,dS.
\end{split}
\end{align}
Then we multiply the equation \eqref{eq:ViscousPart} by 
\begin{equation}
\frac{1}{J}\left(\blockmatx{\tilde{B}}^{\,\text{v}}\left(\textbf{w}\right)\right)^{T}\vec{\nabla}_{\xi}\textbf{v},
\end{equation}
integrate over $\mathrm{E}$ and use \eqref{ViscousBlockFlux3}. This gives    
\begin{align}\label{ContinuousKES4}
\begin{split} 
\int_{\mathrm{E}}\left(\vec{\nabla}_{\xi}\textbf{v}\right)^{T}\left(\blockmatx{\tilde{B}}^{\,\text{v}}\left(\textbf{w}\right)\right)\,\left(\blockvec{\textbf{q}}\right)\,d\vec{\xi} 
=&\int_{\mathrm{E}}\left(\blockvec{\textbf{q}}\right)^{T}\left(\blockmatx{\tilde{B}}^{\,\text{v}}\left(\textbf{w}\right)\right)^{T}\left(\vec{\nabla}_{\xi}\textbf{v}\right)\,d\vec{\xi} \\
=&\int_{\mathrm{E}}\frac{1}{J}\left(\blockmatx{M}\vec{\nabla}_{\xi}\textbf{w}\right)^{T}\blockmatx{B}^{\,\text{v}}\left(\textbf{w}\right)\left(\blockmatx{M}\vec{\nabla}_{\xi}\textbf{v}\right)\,d\vec{\xi}. 
\end{split} 
\end{align}
Next, it follows 
\begin{align}\label{ContinuousKES5}
\begin{split}
\left(\blockmatx{M}\vec{\nabla}_{\xi}\textbf{w}\right)^{T}\blockmatx{B}^{\,\text{v}}\left(\textbf{w}\right)\left(\blockmatx{M}\vec{\nabla}_{\xi}\textbf{v}\right) 
=& \left(\vec{\nabla}_{x}\textbf{v}\right)^{T}\blockmatx{B}^{\,\text{v}}\left(\textbf{w}\right)\left(\vec{\nabla}_{x}\textbf{w}\right) \\
=& \sum_{j=1}^{3}\left(\pderivative{\mathbf{v}}{x_{j}}\right)\textbf{f}_{j}^{\,\text{v}}\left(\textbf{u},\vec{\nabla}_{x}\textbf{u}\right) \\
=& \sum_{i,j=1}^{3}\left(\pderivative{u_{i}}{x_{j}}\right)^{T}\sigma_{ij} \\ 
=& \sum_{i,j=1}^{3}\frac{\mu}{2}\left(\pderivative{u_{i}}{x_{j}}+\pderivative{u_{j}}{x_{i}}-\frac{2}{3}\delta_{ij}\vec{\nabla}\cdot\vec{u}\right)^{2}\geq0,
\end{split}
\end{align}
where $\delta_{ij}$ is the Kronecker delta and $\sigma_{ij}$ are the elements of the shear stress tensor given by \eqref{viscousTensor}. Note that the equality in the penultimate step follows from   the viscous fluxes \eqref{viscousFlux} and the variables $\textbf{v}$ \eqref{KineticEnergy}. The last identity in \eqref{ContinuousKES5} has been proven by Ray in \cite{Ray2017}. Finally, we combine \eqref{ContinuousKES3}, \eqref{ContinuousKES4}, \eqref{ContinuousKES5} and obtain 
\begin{equation}\label{ContinuousKES6}
\pderivative{}{\tau}\int_{\mathrm{E}}J\Bbbk\,d\vec{\xi}\leq-\int_{\mathrm{E}}\left[\vec{\nabla}_{\xi}p\right]\cdot\vec{\tilde{u}}\,d\vec{\xi}-\int_{\partial\mathrm{E}}\left[\left(\vec{\tilde{u}}-\vec{\tilde{\nu}}\right)\Bbbk\right]^{T}\hat{n}\,dS+\int_{\partial\mathrm{E}}\mathbf{v}^{T}\left(\left[\blockmatx{\tilde{B}}^{\,\text{v}}\left(\mathbf{w}\right)\right]\blockvec{\mathbf{q}}\cdot\hat{n}\right)\,dS.
\end{equation} 
  

\subsection{Continuous entropy analysis}\label{EntropyAnalysis}
Analogous to \cite{Schnuecke2018,yamaleev2019entropy}, we obtain for smooth solutions of the equation \eqref{eq:NSELawBlock} 
\begin{equation}\label{EC:Condition1}
\textbf{w}^{T}\pderivative{\left(J\textbf{u}\right)}{\tau}=
\pderivative{\left(Js\right)}{\tau}+\left(\vec{\nabla}_{\xi}\cdot\vec{\tilde{\nu}}\right)\left(\textbf{w}^{T}\textbf{u}-s\right)=
\pderivative{\left(Js\right)}{\tau}+\left(\vec{\nabla}_{\xi}\cdot\vec{\tilde{\nu}}\right)\rho
\end{equation}
and 
\begin{equation}\label{EC:Condition2}
\textbf{w}^{T}\left(\vec{\nabla}_{\xi}\cdot\blockvec{\tilde{\textbf{g}}}\right)=\vec{\nabla}_{\xi}\cdot\left(\vec{\tilde{f}}^{s}-\vec{\tilde{\nu}}s\right)-\left(\vec{\nabla}_{\xi}\cdot\vec{\tilde{\nu}}\right)\left(\textbf{w}^{T}\textbf{u}-s\right)=\vec{\nabla}_{\xi}\cdot\left(\vec{\tilde{f}}^{s}-\vec{\tilde{\nu}}s\right)-\left(\vec{\nabla}_{\xi}\cdot\vec{\tilde{\nu}}\right)\rho,
\end{equation}
where the entropy/entropy flux pair is given by \eqref{EntropyNSE}, the entropy variables $\textbf{w}$ are given by
\eqref{EntropyVariables} and $\vec{f}^{s}:=\left[f_{1}^{s},f_{2}^{s},f_{3}^{s}\right]^{T}$. Next, we multiply the equation \eqref{eq:NSE} with the entropy variables and integrate over $\mathrm{E}$. Then, we obtain by \eqref{EC:Condition1}, \eqref{EC:Condition2}    
\begin{align}\label{ContinuousEntropy1}
\begin{split}
\pderivative{}{\tau}\int_{\mathrm{E}}Js\,d\vec{\xi}+\int_{\partial \mathrm{E}}\left(\vec{\tilde{f}}^{s}-\vec{\tilde{\nu}}s\right)^{T}\hat{n}\,dS  
=& \phantom{+} \int_{\mathrm{E}}\textbf{w}^{T}\frac{d\left(J\textbf{u}\right)}{dt}\,d\vec{\xi}+\int_{\mathrm{E}}\textbf{w}^{T}\left(\vec{\nabla}_{\xi}\cdot\blockvec{\tilde{\textbf{g}}}\right)\,d\vec{\xi} \\
=& \phantom{+} \int_{\mathrm{E}}\textbf{w}^{T}\vec{\nabla}_{\xi}\cdot\left(\blockmatx{\tilde{B}}^{\,\text{v}}\left(\textbf{w}\right)\right)\blockvec{\textbf{q}}\,d\vec{\xi} \\
=& -\int_{\mathrm{E}}\left(\vec{\nabla}_{\xi}\textbf{w}\right)^{T}\left(\blockmatx{\tilde{B}}^{\,\text{v}}\left(\textbf{w}\right)\right)\,\left(\blockvec{\textbf{q}}\right)\,d\vec{\xi} \\
& +\int_{\partial \mathrm{E}}\textbf{w}^{T}\left\{ \left(\blockmatx{\tilde{B}}^{\,\text{v}}\left(\textbf{w}\right)\right)\,\blockvec{\textbf{q}}\cdot\hat{n}\right\} \,dS,
\end{split} 
\end{align}
where $\hat{n}$ is the normal of the domain $\mathrm{E}$. Then we multiply the equation \eqref{eq:ViscousPart} by 
\begin{equation}
\frac{1}{J}\left(\blockmatx{\tilde{B}}^{\,\text{v}}\left(\textbf{w}\right)\right)^{T}\vec{\nabla}_{\xi}\textbf{w},
\end{equation}
integrate over $\mathrm{E}$ and use \eqref{ViscousBlockFlux3}. This results in    
\begin{align}\label{ContinuousEntropy2}
\begin{split} 
\int_{\mathrm{E}}\left(\vec{\nabla}_{\xi}\textbf{w}\right)^{T}\left(\blockmatx{\tilde{B}}^{\,\text{v}}\left(\textbf{w}\right)\right)\,\left(\blockvec{\textbf{q}}\right)\,d\vec{\xi}
=& \int_{\mathrm{E}}\left(\blockvec{\textbf{q}}\right)^{T}\left(\blockmatx{\tilde{B}}^{\,\text{v}}\left(\textbf{w}\right)\right)^{T}\left(\vec{\nabla}_{\xi}\textbf{w}\right)\,d\vec{\xi} \\
=& \int_{\mathrm{E}}\left(\blockvec{\textbf{q}}\right)^{T}\blockmatx{B}^{\,\text{v}}\left(\textbf{w}\right)\left(\blockmatx{M}\vec{\nabla}_{\xi}\textbf{w}\right)\,d\vec{\xi} \\
=& \int_{\mathrm{E}}\frac{1}{J}\left(\blockmatx{M}\vec{\nabla}_{\xi}\textbf{w}\right)^{T}\blockmatx{B}^{\,\text{v}}\left(\textbf{w}\right)\left(\blockmatx{M}\vec{\nabla}_{\xi}\textbf{w}\right)\,d\vec{\xi}
\geq 0, 
\end{split} 
\end{align}
where we used that the block matrix $\blockmatx{B}^{\,\text{v}}$ is symmetric and positive semidefinite. Finally, we combine the equations \eqref{ContinuousEntropy1}, \eqref{ContinuousEntropy2} and obtain 
\begin{equation}\label{ContinuousEntropy3}
\pderivative{}{\tau}\int_{\mathrm{E}}Js\,d\vec{\xi}\leq-\int_{\partial \mathrm{E}}\left(\vec{\tilde{f}}^{s}-\vec{\tilde{\nu}}s\right)^{T}\hat{n}\,dS+\int_{\partial \mathrm{E}}\textbf{w}^{T}\left(\blockmatx{\tilde{B}}^{\,\text{v}}\left(\textbf{w}\right)\blockvec{\textbf{q}}\cdot\hat{n}\right)\,dS.
\end{equation}

\section{Discontinuous Galerkin spectral element method (DGSEM)}\label{Sec:DGSEMNSE}
\subsection{Building blocks for the spectral element approximation}\label{sec:Spectral}
The spectral element approximation is based on a nodal approach with Lagrange basis functions $\left\{ \ell_{j}\left(\cdot\right)\right\} _{j=0}^{N}$ constructed from Legendre Gauss Lobatto (LGL) points $\left\{ \xi_{i}\right\} _{i=0}^{N}$. We note that $\xi_{0}=-1$ and $\xi_{N}=1$. The Lagrange basis functions satisfy the cardinal property  
\begin{equation}\label{CardinalProperty}
\ell_{i}\left(\xi_{j}\right)=\delta_{ji}.
\end{equation}  
On the reference element $\mathrm{E}=\left[-1,1\right]^{3}$ the solution and fluxes of the system \eqref{eq:NSELawBlockRefE} are approximated by tensor product Lagrange polynomials of degree $N$, e.g., 
\begin{subequations}\label{UnkonwnNSE}
\begin{align}
J\left(\xi^{1},\xi^{2},\xi^{3},t\right)\approx\mathcal{J}\left(\xi^{1},\xi^{2},\xi^{3},t\right):=&\sum_{i,j,k=0}^{N}\mathcal{J}_{ijk}\left(t\right)\ell_{i}\left(\xi^{1}\right)\ell_{j}\left(\xi^{2}\right)\ell_{k}\left(\xi^{3}\right). \label{SpatialPolynomialMetricQuantity} \\
\textbf{u}\left(\xi^{1},\xi^{2},\xi^{3},t\right)\approx\textbf{U}\left(\xi^{1},\xi^{2},\xi^{3},t\right):=&
\sum_{i,j,k=0}^{N}\textbf{U}_{ijk}\left(t\right)\ell_{i}\left(\xi^{1}\right)\ell_{j}\left(\xi^{2}\right)\ell_{k}\left(\xi^{3}\right), \label{SpatialPolynomial} \\
\blockvec{\textbf{q}}\left(\xi^{1},\xi^{2},\xi^{3},t\right)\approx\blockvec{\textbf{Q}}\left(\xi^{1},\xi^{2},\xi^{3},t\right):=&\sum_{i,j,k=0}^{N}\blockvec{\textbf{Q}}{}_{ijk}\left(t\right)\ell_{i}\left(\xi^{1}\right)\ell_{j}\left(\xi^{2}\right)\ell_{k}\left(\xi^{3}\right). \label{SpatialPolynomialQ} 
\end{align}
\end{subequations}
From now on, polynomial approximations are highlighted by capital letters, e.g., $\textbf{U}$ is an approximation for the state vector $\textbf{u}$,    $\blockvec{\textbf{Q}}$ is an approximation for the solution of equation \eqref{eq:ViscousPart} and $\textbf{F}_{\iota}$, $\iota=1,2,3$, are approximations for the fluxes $\textbf{f}_{\iota}$, $\iota=1,2,3$. The approximation for the determinant $J$ of the Jacobian matrix $\vec{\nabla}_{\xi}\vec{\chi}$ is highlighted by $\mathcal{J}$. Furthermore, the interpolation operator for a function $\textbf{g}$ is given by 
\begin{equation}\label{InterpolationOperator}  
\interpolation{N}{\left(\textbf{g}\right)}\left(\xi^{1},\xi^{2},\xi^{3}\right)=\sum_{i,j,k=0}^{N}\textbf{g}_{ijk}\ell_{i}\left(\xi^{1}\right)\ell_{j}\left(\xi^{2}\right)\ell_{k}\left(\xi^{3}\right),
\end{equation}
where $\textbf{g}_{ijk}:=\textbf{g}\left(\xi_{i}^{1},\xi_{j}^{2},\xi_{k}^{3}\right)$ and $\left\{ \xi_{i}^{1}\right\} _{i=0}^{N}$, $\left\{\xi_{i}^{2}\right\} _{i=0}^{N}$, $\left\{\xi_{i}^{3}\right\} _{i=0}^{N}$ are sets of LGL points. Derivatives are approximated by exact differentiation of the polynomial interpolants. In general we have $ \left(\interpolation{N}{\left(g\right)}\right)'\neq \interpolation{N-1}{\left(g'\right)}$ (cf. e.g. \cite{CHQZ:2006,Kopriva2009}), as differentiation and interpolation only commute if there are no interpolation errors. 

\subsubsection{Discrete integrals}
Integrals are approximated by a tensor product extension of a $2N-1$ accurate LGL quadrature formula. Hence, interpolation and quadrature nodes are collocated. In one spatial dimension the LGL quadrature formula is given by   
\begin{equation}\label{GLQ}
\int\limits_{-1}^{1}g\left(\xi\right)\,d\xi\approx\sum_{i=0}^{N}\omega_{i}g\left(\xi_{i}\right)=\sum_{i=0}^{N}\omega_{i}g_{i},  
\end{equation}
where $\omega_{i}$, $i=0,\dots,N$, are the quadrature weights and $\xi_{i}$, $i=0,\dots,N$, are the LGL quadrature points. The formula \eqref{GLQ} motivates the definition of the inner product notation 
\begin{equation}\label{Innerproduct}
\left\langle \textbf{f},\textbf{g}\right\rangle _{N}:=\sum_{i=0}^{N}\sum_{j=0}^{N}\sum_{k=0}^{N}\omega_{i}\omega_{j}\omega_{k}\textbf{f}_{ijk}^{T}\textbf{g}_{ijk}=\sum_{i,j,k=0}^{N}\omega_{ijk}\textbf{f}_{ijk}^{T}\textbf{g}_{ijk}
\end{equation}
for two functions $\textbf{f}$ and $\textbf{g}$. We note that the inner product \eqref{Innerproduct} satisfies        
\begin{equation}\label{MagicProperty}
\left\langle \interpolation{N}{\left(\textbf{g}\right)},\boldsymbol{\varphi}\right\rangle _{N}=\left\langle \textbf{g},\boldsymbol{\varphi}\right\rangle _{N}, \qquad \forall \boldsymbol{\varphi}\in\mathbb{P}^{N}\left(\mathrm{E},\R^{5}\right), 
\end{equation}
where 
\begin{equation}
\mathbb{P}^{N}\left(\mathrm{E},\R^{5}\right):=\left\{ \boldsymbol{\varphi}\mid\boldsymbol{\varphi}=\left[\varphi_{1},\varphi_{2},\varphi_{3},\varphi_{4},\varphi_{5}\right]^{T},\quad\varphi_{i}\in\mathbb{P}^{N}\left(\mathrm{E}\right),\quad i=1,2,3,4,5\right\}, 
\end{equation} 
and $\mathbb{P}^{N}\left(\mathrm{E}\right)$ is the space of tensor product polynomials with three dimensional domain. Furthermore, for a block vector $\blockvec{\textbf{F}}$ and test functions $\boldsymbol{\varphi}\in\mathbb{P}^{N}\left(\mathrm{E},\R^{5}\right)$, we define the discrete surface integral  
\begin{align}\label{DiscreteSpatialSurface1}
\begin{split}
\int\limits _{\partial\mathrm{E},N}\boldsymbol{\varphi}^{T}\left\{ \blockvec{\textbf{F}}\cdot\hat{n}\right\} \,dS=\sum_{\iota=1}^{3}\int\limits _{\partial\mathrm{E},N}\boldsymbol{\varphi}^{T}\mathbf{F}_{\iota}\hat{n}^{\iota}\,dS,
\end{split}
\end{align}
where $\hat{n}=\left[\hat{n}^{1},\hat{n}^{2},\hat{n}^{3}\right]^{T}$ is the unit outward normal at the faces of the reference element $\mathrm{E}$ and  
\begin{align}\label{DiscreteSpatialSurface2}
\begin{split}
\int\limits _{\partial\mathrm{E},N}\boldsymbol{\varphi}^{T}\mathbf{F}_{1}\hat{n}^{1}\,dS:=&\sum_{j,k=0}^{N}\omega_{j}\omega_{k}\left(\boldsymbol{\varphi}_{Njk}^{T}\left(\mathbf{F}_{1}\right)_{Njk}-\boldsymbol{\varphi}_{0jk}^{T}\left(\mathbf{F}_{1}\right)_{0jk}\right), 
\\ 
\int\limits _{\partial\mathrm{E},N}\boldsymbol{\varphi}^{T}\mathbf{F}_{2}\hat{n}^{2}\,dS:=&\sum_{i,k=0}^{N}\omega_{i}\omega_{k}\left(\boldsymbol{\varphi}_{iNk}^{T}\left(\mathbf{F}_{2}\right)_{iNk}-\boldsymbol{\varphi}_{i0k}^{T}\left(\mathbf{F}_{2}\right)_{i0k}\right), \\
\int\limits _{\partial\mathrm{E},N}\boldsymbol{\varphi}^{T}\mathbf{F}_{3}\hat{n}^{3}\,dS:=&\sum_{i,j=0}^{N}\omega_{i}\omega_{j}\left(\boldsymbol{\varphi}_{ijN}^{T}\left(\mathbf{F}_{3}\right)_{ijN}-\boldsymbol{\varphi}_{ij0}^{T}\left(\mathbf{F}_{3}\right)_{ij0}\right).
\end{split}
\end{align}
\subsubsection{Discrete metric identities} 
The contravariant coordinate vectors need to be discretized in such a way that the metric identities \eqref{MetricIdentities} are satisfied on the discrete level, too. Kopriva \cite{Kopriva2006} introduced the conservative curl form to approximate the metric terms. In this approach the coefficients of the volume weighted contravariant coordinate vectors $J\vec{a}^{\,\iota}$, $\iota =1,2,3$, are computed by
\begin{equation}\label{DiscreteContravariantVectors}
J\vec{a}_{\beta}^{\,\iota}\approx\mathbb{J}\vec{a}_{\beta}^{\,\iota}:=-\hat{x}_{\iota}\cdot\nabla_{\xi}\times\left(\interpolation{N}{\left(\chi_{\gamma}\nabla_{\delta}\chi_{m}\right)}\right),\quad\iota=1,2,3,\quad\beta=1,2,3,\quad\left(\beta,\gamma,\delta\right)\ \text{cyclic}.
\end{equation}
Here $\vec{\chi}=\left[\chi_{1},\chi_{2},\chi_{3}\right]^{T}$ represents the mapping from the physical element to the reference element and $\hat{x}_{\iota}$ is the unit vector in the $\iota$-th Cartesian coordinate direction. The representation \eqref{DiscreteContravariantVectors} ensures that   
\begin{equation}\label{DiscreteMetricIdentities1}
\sum_{\iota=1}^{3}\pderivative{\,\interpolation{N}{\left(\mathbb{J}\vec{a}^{\,\iota}\right)}}{\xi^{\iota}}=0.
\end{equation}

\subsubsection{SBP operator}   
The spectral element approximation with LGL points for interpolation and quadrature  gives a SBP operator  $\matx{Q}=\matx{M}\,\matx{D}$ with the mass matrix $\matx{M}$ and the derivative matrix $\matx{D}$. The mass matrix and the derivative matrix have the entries 
\begin{equation}\label{SBPcoefficient}
\mathcal{M}_{ij}=\omega_{i}\delta_{ij}, \qquad 
\mathcal{D}_{ij}=\ell_{j}'\left(\xi_{i}\right), \qquad 
i,j=0,\dots,N.
\end{equation}
A SBP operator satisfies the property     
\begin{equation}\label{SBP}
\matx{Q}+\matx{Q}^T=\matx{B},  
\end{equation}  
where $\matx{B}=\text{diag}\left(-1,0,\dots,0,1\right)$, \cite{Fernandez2014,Gassner2013,kreiss1}.   

\subsection{Split form ALE DGSEM for the Euler equations}\label{sec:DGSEMEuler}
The construction of the split form ALE DGSEM for the Euler equations is analogous to \cite{Schnuecke2018}. First, we replace the solution $\textbf{u}$ by \eqref{SpatialPolynomial}, the Jacobian $J$ by \eqref{SpatialPolynomialMetricQuantity} and approximate the fluxes by the interpolation operator \eqref{InterpolationOperator}. Next, we multiply the equation \eqref{GCL} with test functions $\varphi\in\mathbb{P}^{N}\left(\mathrm{E}\right)$ and the transformed Euler equations with 
$\boldsymbol{\varphi}\in\mathbb{P}^{N}\left(\mathrm{E},\R^{5}\right)$, integrate the resulting equations and use integration-by-parts to separate boundary and volume contributions. The volume integrals in the variational form are approximated with the LGL quadrature. Then, we insert numerical surface fluxes $\vec{\tilde{\nu}}^{*}$ and $\blockvec{\tilde{\textbf{G}}}{}^{*}$ at the spatial element interfaces. Afterwards, we use the SBP property \eqref{SBP} for the volume contribution to get the standard ALE DGSEM in strong form:
\begin{subequations}\label{StandardDGSEM}
\begin{align}
\left\langle \pderivative{\mathcal{J}}{\tau},\varphi\right\rangle_{N} 
=& \phantom{-} 
\left\langle \vec{\nabla}_{\xi}\cdot\interpolation{N}{\left(\vec{\tilde{\nu}}\right)},\varphi\right\rangle _{N}+\int\limits _{\partial \mathrm{E},N}\varphi\left(\tilde{\nu}_{\hat{n}}^{*}-\tilde{\nu}_{\hat{n}}\right)\,dS, \qquad \forall \varphi\in\mathbb{P}^{N}\left(\mathrm{E}\right), \label{StandardDGSEM:D-GLC} \\   
\left\langle \pderivative{\left(\mathcal{J}\textbf{U}\right)}{\tau},\boldsymbol{\varphi}\right\rangle_{N} 
=& -\left\langle \vec{\nabla}_{\xi}\cdot\interpolation{N}{\left(\blockvec{\tilde{\textbf{g}}}\right)},\boldsymbol{\varphi}\right\rangle _{N}
-\int\limits _{\partial \mathrm{E},N}\boldsymbol{\varphi}^{T}\left(\tilde{\textbf{G}}_{\hat{n}}^{*}-\tilde{\textbf{G}}_{\hat{n}}\right)\,dS,
\qquad \forall 
\boldsymbol{\varphi}\in\mathbb{P}^{N}\left(\mathrm{E},\R^{5}\right). 
\end{align}
\end{subequations}
It is known that the interpolation of $\vec{\tilde{\nu}}$ and the nonlinear flux $\blockvec{\tilde{\textbf{g}}}$ (see \eqref{G-flux}) causes aliasing errors in the standard strong form \cite{Gassner2017, Gassner2016}. Thus, we follow \cite{Gassner2016} and introduce a special derivative operator and change the divergence in equation \eqref{StandardDGSEM:D-GLC} to
\begin{align}\label{SpatialDerivativeProjectionOperator1}
\begin{split}
\Dprojection{N}\cdot\vec{\tilde{\nu}}_{ijk}^{\#}:=\sum_{m=0}^{N}
& \phantom{+} \,2\,\mathcal{D}_{im}\left(\avg{\vec{\nu}}_{\left(i,m\right)jk}\cdot\avg{\mathbb{J}\vec{a}^{1}}_{\left(i,m\right)jk}\right) \\
& + \,2\,\mathcal{D}_{jm}\left(\avg{\vec{\nu}}_{i\left(j,m\right)k}\cdot\avg{\mathbb{J}\vec{a}^{2}}_{i\left(j,m\right)k}\right) \\
& + \,2\,\mathcal{D}_{km}\left(\avg{\vec{\nu}}_{ij\left(m,k\right)}\cdot\avg{\mathbb{J}\vec{a}^{3}}_{ij\left(m,k\right)}\right),
\end{split}
\end{align}
with the volume averages 
\begin{align}\label{VolumeAverages}
\avg{\ast}_{(i,m)jk}:=\frac{1}{2}\left[\left(\ast\right)_{ ijk}+\left(\ast\right)_{ mjk}\right].
\end{align}
The derivative projection operator for the Euler fluxes is computed as
\begin{align}\label{SpatialDerivativeProjectionOperator2}
\begin{split}
\Dprojection{N}\cdot\blockvec{\tilde{\textbf{G}}}{}_{ijk}^{\#}:=& \sum_{m=0}^{N}\quad\,2\,\mathcal{D}_{im} \left(\blockvec{\textbf{G}}{}^{\#}\left(\vec{\nu}_{ijk},\vec{\nu}_{mjk},\textbf{U}_{ijk},\textbf{U}_{mjk}\right)\cdot\avg{\mathbb{J}\vec{a}^{1}}_{\left(i,m\right)jk}\right)\\
& \quad \ \ \quad+2\,\mathcal{D}_{jm} 
\left(\blockvec{\textbf{G}}{}^{\#}\left(\vec{\nu}_{ijk},\vec{\nu}_{imk},\textbf{U}_{ijk},\textbf{U}_{imk}\right)\cdot\avg{\mathbb{J}\vec{a}^{2}}_{i\left(j,m\right)k}\right) \\
& \quad \ \ \quad+2\,\mathcal{D}_{km}\left(\blockvec{\textbf{G}}{}^{\#}\left(\vec{\nu}_{ijk},\vec{\nu}_{ijm},\textbf{U}_{ijk},\textbf{U}_{ijm}\right)\cdot\avg{\mathbb{J}\vec{a}^{3}}_{ij\left(k,m\right)}\right).
\end{split}
\end{align}
The functions $\tilde{\textbf{G}}{}_{\iota}^{\#}$, $\iota=1,2,3$, in \eqref{SpatialDerivativeProjectionOperator2} are Cartesian numerical volume fluxes. In Appendix \ref{NumericalFlux} possible choices for two point volume flux functions are given. These fluxes discretely recover corresponding split formulations, see \cite{Gassner2016}. It is important that the volume fluxes are consistent and symmetric such that for all $i,j,k,m=0,\dots,N$ holds 
\begin{equation*}
\textbf{G}{}_{\iota}^{\#}\left(\vec{\nu}_{ijk},\vec{\nu}_{mjk},\textbf{U}_{ijk},\textbf{U}_{ijk}\right)=\mathbf{F}_{\iota}\left(\textbf{U}_{ijk}\right)-\avg{\vec{\nu}_{\iota}}_{\left(i,m\right)jk}\textbf{U},\qquad\iota=1,2,3,
\end{equation*}
\begin{equation}\label{FluxSymmetric}
\textbf{G}{}_{\iota}^{\#}\left(\vec{\nu}_{ijk},\vec{\nu}_{mjk},\textbf{U}_{ijk},\textbf{U}_{mjk}\right)=\textbf{G}{}_{\iota}^{\#}\left(\vec{\nu}_{mjk},\vec{\nu}_{ijk},\textbf{U}_{mjk},\textbf{U}_{ijk}\right),\qquad\iota=1,2,3,
\end{equation}  
where $\mathbf{F}_{\iota}\left(\textbf{U}_{ijk}\right)$  is the advective Euler flux \eqref{AdvectiveFluxes} evaluated at the LGL points. Then, for each element $e_{\alpha}(t)$ the semi-discrete split form ALE DGSEM for the Euler equations is represented on the reference element $\mathrm{E}$ by: 
\begin{subequations}\label{MovingMeshDGSEMEuler}
\begin{align}
\left\langle 
\pderivative{\mathcal{J}}{\tau},\varphi\right\rangle_{N} 
=& \phantom{-} \left\langle \Dprojection{N}\cdot\vec{\tilde{\nu}}^{\#},\varphi\right\rangle _{N}+\int\limits _{\partial \mathrm{E},N}\varphi\left(\tilde{\nu}_{\hat{n}}^{*}-\tilde{\nu}_{\hat{n}}\right)\,dS, \qquad \forall \varphi\in\mathbb{P}^{N}\left(\mathrm{E}\right), \label{DGSEM:D-GLC} \\ 
\left\langle \pderivative{\left(\mathcal{J}\textbf{U}\right)}{\tau},\boldsymbol{\varphi}\right\rangle_{N} 
=& -\left\langle \Dprojection{N}\cdot\blockvec{\tilde{\textbf{G}}}{}^{\#},\boldsymbol{\varphi}\right\rangle _{N}
-\int\limits _{\partial \mathrm{E},N}\boldsymbol{\varphi}^{T}\left(\tilde{\textbf{G}}_{\hat{n}}^{*}-\tilde{\textbf{G}}_{\hat{n}}\right)\,dS,
\qquad \forall 
\boldsymbol{\varphi}\in\mathbb{P}^{N}\left(\mathrm{E},\R^{5}\right). \label{DGSEM:CL} 
\end{align}
\end{subequations}
It is important to mention that the split form 
ALE DGSEM \eqref{MovingMeshDGSEMEuler} is discretely conservative. 
\noindent
The unit outward facing normal vector and surface element on the element side are constructed from the element metrics by
\begin{equation}\label{ComputationNormal}
\vec{n}:=\frac{1}{\hat{s}}\sum_{\iota=1}^{3}\left(\mathbb{J}\vec{a}^{\iota}\right)\hat{n}^{\iota},\qquad\hat{s}:=\left|\sum_{\iota=1}^{3}\left(\mathbb{J}\vec{a}^{\iota}\right)\hat{n}^{\iota}\right|.
\end{equation}
Thus, the quantity $\tilde{\nu}_{\hat{n}}$ in \eqref{DGSEM:D-GLC} and the flux $\tilde{\textbf{G}}_{\hat{n}}$ in \eqref{DGSEM:CL} are defined by 
\begin{align}
\tilde{\nu}_{\hat{n}}=\left(\hat{s}\vec{n}\right)\cdot\vec{\nu}=\sum_{\iota=1}^{3}\hat{n}^{\iota}\left(\mathbb{J}a_{1}^{\iota}\nu_{1}+\mathbb{J}a_{2}^{\iota}\nu_{2}+\mathbb{J}a_{3}^{\iota}\nu_{3}\right), \label{GridvelocityNormal} \\
\tilde{\textbf{G}}_{\hat{n}}=\left(\hat{s}\vec{n}\right)\cdot\blockvec{\textbf{G}}=\sum_{\iota=1}^{3}\hat{n}^{\iota}\left(\mathbb{J}a_{1}^{\iota}\textbf{G}_{1}+\mathbb{J}a_{2}^{\iota}\textbf{G}_{2}+\mathbb{J}a_{3}^{\iota}\textbf{G}_{3}\right).\label{BlockvectorNormal}
\end{align}
To define the numerical surface fluxes in \eqref{DGSEM:D-GLC} and \eqref{DGSEM:CL}, we introduce notation for states at the LGL nodes along an interface between two spatial elements to be a primary ``$-$'' and complement the notation with a secondary ``$+$'' to denote the value at the LGL nodes on the opposite side. Then the orientated jump and
the averages at the interfaces are defined by  
\begin{equation}\label{SurfaceJumpMean}
\jump{\ast}:=\left(\ast\right)^{+}-\left(\ast\right)^{-},
\quad \text{and} \quad 
\avg{\ast}:=\frac{1}{2}\left[\left(\ast\right)^{+}+\left(\ast\right)^{-}\right].  
\end{equation}
When applied to vectors, the average and jump operators are evaluated separately for each vector component. Then the normal vector $\vec{n}$ is defined uniquely to point from the ``$-$'' to the ``$+$'' side. This notation allows to compute the  contravariant surface numerical fluxes in \eqref{DGSEM:D-GLC} as
\begin{equation}
\tilde{\nu}_{\hat{n}}^{*}=\hat{s}\left(n_{1}\avg{v_{1}}+n_{2}\avg{v_{2}}+n_{3}\avg{v_{3}}\right).
\end{equation}
The contravariant surface numerical fluxes in \eqref{DGSEM:CL} are given by
\begin{equation}\label{ContravariantSurfaceFlux}
\tilde{\textbf{G}}_{\hat{n}}^{*}=\hat{s}\left(n_{1}\textbf{G}_{1}^{*}+n_{2}\textbf{G}_{2}^{*}+n_{3}\textbf{G}_{3}^{*}\right).
\end{equation}
In order to construct an EC or KEP ALE DGSEM for the Euler equations, the Cartesian fluxes $\textbf{G}_{\iota}^{*}$, $\iota=1,2,3$, should be consistent with $\textbf{F}_{\iota}-\nu_{\iota}\textbf{U}$, symmetric and EC or KEP in the sense that Tadmor's \cite{Tadmor2003} discrete entropy or Jameson's \cite{Jameson2008} discrete criteria are satisfied. In Appendix \ref{NumericalFlux} common choices for these two point flux functions are given. We note that these variants are virtually dissipation free and may become unstable for solutions with very steep gradients. Hence, we might add dissipation to the surface fluxes $\textbf{G}_{\iota}^{*}$. Then the numerical surface fluxes can be computed by  
\begin{equation}\label{CartesianSurfaceFlux}
\textbf{G}_{\iota}^{*}:=\textbf{G}_{\iota}^{\star}-\frac{1}{2}\,\matx{H}_{\iota}\,\jump{\textbf{W}},\qquad \iota=1,2,3,
\end{equation}
where the quantities $\textbf{W}$ are the interpolated entropy variables  \eqref{EntropyVariables} evaluated in the LGL points. The Cartesian flux \eqref{CartesianSurfaceFlux} is constructed with a consistent and symmetric two point flux $\textbf{G}_{\iota}^{\star}$ and a matrix dissipation operator $\matx{H}_{\iota}$. The dissipation operator $\matx{H}_{\iota}$ is a symmetric positive definite matrix of the form 
\begin{equation}\label{EntropyBasedDissipation}
\matx{H}_{\iota}=\hat{\matx{R}}_{\iota}\,\left|\Lambda_{\iota}\right|\,\hat{\matx{R}}_{\iota}^{T},\qquad\hat{\matx{R}}_{\iota}=\matx{R}_{\iota}\matx{T}_{\iota},\qquad \iota=1,2,3,
\end{equation}
where the matrices $\matx{R}_{\iota}$, $\matx{T}_{\iota}$, depend on the averaged values of the states $\textbf{U}^{-}$, $\textbf{U}^{+}$ and they are consistent with the right eigenvector matrix, which corresponds to the flux Jacobian matrices from the advective fluxes \eqref{AdvectiveFluxes}, and a diagonal scaling matrix. The matrix $\left|\Lambda_{\iota}\right|$ depends on the eigenvalues of the flux Jacobian matrices from the advective fluxes \eqref{AdvectiveFluxes}. In order to construct an ES or KED ALE DGSEM for the Euler equations, the matrix $\matx{H}_{\iota}$ needs to be a symmetric positive definite matrix. In \cite[Appendix C.3]{Schnuecke2018} a suitable matrix dissipation operator is given. The construction of such a matrix dissipation operator is based on the fact that there are block diagonal scaling matrices such that the Hessian matrix of the entropy \eqref{EntropyNSE} can be represented by scaled right eigenvector matrices (cf. Merriam \cite{Merriam1989}). 
 
\subsection{Analysis of the discrete kinetic and internal energy balance}\label{Sec:KESEuler}
On the continuous level, the Euler equations describe the kinetic energy balance \eqref{ContinuousKESEuler} and internal energy balance \eqref{ContinuousInnerEuler} for smooth solutions. The ALE DGSEM \eqref{MovingMeshDGSEMEuler} should have discrete analogues of these  balance laws. The next Theorem provides an identity for the evolution of the discrete kinetic energy. A proof of this identity is given in Appendix \ref{ProofTheorem3_1}.
\begin{thm-hand}[3.1.]
Suppose the flux functions $\textbf{G}{}_{\iota}^{\#}$, $\iota=1,2,3$, in the derivative projection operator \eqref{SpatialDerivativeProjectionOperator2} and the surface flux functions $\textbf{G}{}_{\iota}^{*}$, $\iota=1,2,3$, \eqref{CartesianSurfaceFlux} satisfy Jameson's \cite{Jameson2008} conditions 
\begin{align}\label{Jameson}
\textbf{G}{}_{\iota}^{\upsilon+1,\#}=\textbf{G}{}_{\iota}^{1,\#}\avg{u_{\upsilon}}+\avg{p}\delta_{\iota\upsilon},  \qquad 
\mathbf{G}{}_{\iota}^{\upsilon+1,*}=\mathbf{G}{}_{\iota}^{1,*}\avg{u_{\upsilon}}+\avg{p}^{\star}\delta_{\iota\upsilon},\qquad\upsilon=1,2,3,
\end{align}
where $\avg{p}^{\star}$ can be any consistent numerical trace approximation of the pressure. Then the split form ALE DGSEM \eqref{MovingMeshDGSEMEuler} satisfies for each element $e_{\alpha}\left(t\right)$, $\alpha=1,\dots,K$, the identity  
\begin{align}\label{LocalDisKEP}
\begin{split}
\pderivative{}{\tau}\left\langle \interpolation{N}{\left(\Bbbk\right)},\mathcal{J}\right\rangle _{N}=& 
-\frac{1}{2}\sum_{\iota=1}^{3}\left\langle \left(\pderivative{\interpolation{N}{\left(p\right)}}{\xi^{\iota}}\right)\mathbb{J}\vec{a}^{\iota}+p\left(\pderivative{\interpolation{N}{\left(\mathbb{J}\vec{a}^{\iota}\right)}}{\xi^{\iota}}\right)+\left(\pderivative{\interpolation{N}{\left(p\mathbb{J}\vec{a}^{\iota}\right)}}{\xi^{\iota}}\right),\vec{u}\right\rangle _{N}\\
&-\sum_{\iota=1}^{3}\ \int\limits _{\partial\mathrm{E},N}\hat{s}n_{\iota}\left[\frac{1}{2}\overline{\left|\vec{u}\right|}^{2}\mathbf{G}_{\iota}^{1,*}+\left(\avg{p}^{\star}-p^{-}\right)u_{\iota}^{-}\right]\,dS,
\end{split}
\end{align}
where the kinetic energy $\Bbbk$ is given by \eqref{KineticEnergy}, the Cartesian surface fluxes $\mathbf{G}_{\iota}^{1,*}$ are consistent with 
\begin{equation}
\mathbf{G}_{\iota}^{1}= 
\rho\left(u_{\iota}-\nu_{\iota}\right),\qquad \iota=1,2,3,
\end{equation}
and 
\begin{equation}\label{NormAVG}
\overline{\left|\vec{u}\right|}^{2}:=\sum_{\iota=1}^{3}2\avg{u_{\iota}}^{2}-\avg{u_{\iota}^{2}}.
\end{equation}  
The quantities $u_{\iota}^{-}$, $\iota=1,2,3$, and $p^{-}$ are interior contributions form the element $e_{\alpha}\left(t\right)$, $\alpha=1,\dots,K$, and the ``$-$'' has to be understood as in \eqref{SurfaceJumpMean}.
\end{thm-hand}
\begin{rem-hand}[3.2.]
Since the contravariant coordinate vectors $J\vec{a}^{\iota}\approx\mathbb{J}\vec{a}^{\iota}$, $\iota=1,2,3$, are discretized by \eqref{DiscreteContravariantVectors}, the discrete metric identities \eqref{DiscreteMetricIdentities1} are satisfied and we obtain the identity 
\begin{align}\label{DiscreteMetric_ZeroIntegral}
\sum_{\iota=1}^{3}\left\langle p\left(\pderivative{\interpolation{N}{\left(\mathbb{J}\vec{a}^{\iota}\right)}}{\xi^{\iota}}\right),\vec{u}\right\rangle _{N}=0.
\end{align}
Furthermore, the contravariant coordinate vectors are constant on a Cartesian mesh. Hence, the result in Theorem 3.1 is consistent with other results about kinetic energy preserving DG methods \cite{Gassner2016,Ranocha2018} when a Cartesian mesh is investigated.      
\end{rem-hand}
It is important to mention that in Theorem 3.1 Jameson's conditions for the volume flux functions $\textbf{G}{}_{\iota}^{\#}$, $\iota=1,2,3$, differ from the conditions for the surface flux functions $\textbf{G}{}_{\iota}^{*}$, $\iota=1,2,3$, since the proof in Appendix \ref{ProofTheorem3_1} requires that the pressure in the discrete momentum equations ($\mathbf{G}{}_{\iota}^{\upsilon+1,\#}$, $\iota, \upsilon=1,2,3$) of the volume flux functions is given by the average operator $\avg{p}$. It is not clear, if this restriction for the pressure in the discrete momentum equations is sufficient or necessary to obtain a discrete equation for the evolution of the discrete kinetic energy $\interpolation{N}{\left(\Bbbk\right)}$. However, Gassner et al. \cite[Fig. 4 and Fig. 5]{Gassner2016} investigated the evolution of the kinetic energy for the inviscid Taylor-Green vortex (TGV) test case \cite{shu2005numerical} on a static mesh. The flux from Chandrashekar \cite{Chandrashekar2013} (see Appendix \ref{CH}) was used as volume flux in the DGSEM scheme and a decrease of the total kinetic energy was observed. On the other hand, the total kinetic energy is preserved when the pressure term $\avg{\rho}/\avg{\frac{\rho}{p}}$ is replaced by $\avg{p}$ in the Chandrashekar flux and the corresponding flux (see the flux \eqref{CH_AVG_Flux}) is used as volume flux in the DGSEM scheme. We observe the same behavior on a moving mesh in our numerical experiments (see Figure~\ref{fig:TGV_EC_KEP_Euler} in Section \ref{EC_KEP_Test}).  

\subsubsection{Discrete internal energy balance}
The evolution of the discrete kinetic energy needs to be consistent with the evolution of the discrete total energy which is given by   
\begin{align}\label{LocalDisTotalEnergy}
\pderivative{}{\tau}\left\langle \interpolation{N}{\left(\rho e+\frac{1}{2}\rho\left|\vec{u}\right|^{2}\right)},\mathcal{J}\right\rangle _{N}=-\int\limits _{\partial\mathrm{E},N}\tilde{\mathbf{G}}_{\hat{n}}^{5,*}\,dS=-\sum_{\sigma=1}^{3}\ \int\limits _{\partial\mathrm{E},N}\hat{s}\left(n_{\iota}\mathbf{G}_{\iota}^{5,*}\right)\,dS,
\end{align}
where the Cartesian surface fluxes $\mathbf{G}_{\iota}^{5,*}$ are consistent with    
\begin{equation}
\mathbf{G}_{\iota}^{5}=\rho\left(u_{\iota}-\nu_{\iota}\right)\left[e+\frac{1}{2}\left|\vec{u}\right|^{2}\right]+pu_{\iota},\qquad\iota=1,2,3.
\end{equation}
The equation \eqref{LocalDisTotalEnergy} results from the split form ALE DGSEM \eqref{MovingMeshDGSEMEuler} when the equation \eqref{DGSEM:CL} is tested with a constant state. In particular, as in the continuous case (see Section \ref{KineticEnergyEuler}), the discrete total and kinetic energy give the evolution of the discrete internal energy. The equations \eqref{LocalDisKEP} and \eqref{LocalDisTotalEnergy} provide        
\begin{align}\label{LocalDisInternalEnergy}
\begin{split}
\pderivative{}{\tau}\left\langle \interpolation{N}{\left(\rho e\right)},\mathcal{J}\right\rangle _{N} 
=& \phantom{+} \pderivative{}{\tau}\left\langle \interpolation{N}{\left(\rho e+\frac{1}{2}\rho\left|\vec{u}\right|^{2}\right)},\mathcal{J}\right\rangle _{N}-
\pderivative{}{\tau}\left\langle \interpolation{N}{\left(\Bbbk\right)},J\right\rangle _{N} \\ 
=& \phantom{-}
\frac{1}{2}\sum_{\iota=1}^{3}\left\langle \left(\pderivative{\interpolation{N}{\left(p\right)}}{\xi^{\iota}}\right)\mathbb{J}\vec{a}^{\iota}+p\left(\pderivative{\interpolation{N}{\left(\mathbb{J}\vec{a}^{\iota}\right)}}{\xi^{\iota}}\right)+\left(\pderivative{\interpolation{N}{\left(p\mathbb{J}\vec{a}^{\iota}\right)}}{\xi^{\iota}}\right),\vec{u}\right\rangle _{N}\\
&-\sum_{\iota=1}^{3}\ \int\limits _{\partial\mathrm{E},N}\hat{s}n_{\iota}\left[\mathbf{G}_{\iota}^{5,*}-\frac{1}{2}\overline{\left|\vec{u}\right|}^{2}\mathbf{G}_{\iota}^{1,*}-\left(\avg{p}^{*}-p^{-}\right)u_{\iota}^{-}\right]\,dS
\\
=&-\sum_{\iota=1}^{3}\left(\left\langle p,\pderivative{\interpolation{N}{\left(\mathbb{J}\vec{a}^{\iota}\cdot\vec{u}\right)}}{\xi^{\iota}}\right\rangle _{N}-\left\langle S\left(p,\mathbb{J}\vec{a}^{\iota}\right),\vec{u}\right\rangle _{N}\right) \\
&+\sum_{\iota=1}^{3}\ \int\limits _{\partial\mathrm{E},N}\hat{s}n_{\iota}\left[\mathbf{G}_{\iota}^{5,*}-\frac{1}{2}\overline{\left|\vec{u}\right|}^{2}\mathbf{G}_{\iota}^{1,*}-\avg{p}^{*}u_{\iota}^{-}\right]\,dS,
\end{split}
\end{align}
with   
\begin{equation}\label{ConservedZero}
S\left(p,\mathbb{J}\vec{a}^{\iota}\right)=\frac{1}{2}\left[p\left(\pderivative{\interpolation{N}{\left(\mathbb{J}\vec{a}^{\iota}\right)}}{\xi^{\iota}}\right)+\left(\pderivative{\interpolation{N}{\left(p\mathbb{J}\vec{a}^{\iota}\right)}}{\xi^{\iota}}\right)-\left(\pderivative{\interpolation{N}{\left(p\right)}}{\xi^{\iota}}\right)\mathbb{J}\vec{a}^{\iota}\right], \qquad \iota=1,2,3, 
\end{equation}
where the last equality results from the SBP property \eqref{SBP}. 
Therefore, a split form ALE DGSEM \eqref{MovingMeshDGSEMEuler}, which preserves the kinetic and internal energy such that the total energy is correctly balanced, can be constructed, if it is ensured that:    
\begin{itemize}
\item[\textbf{(R1)}] The volume flux functions $\textbf{G}{}_{\iota}^{\#}$, $\iota=1,2,3$, in the derivative operator \eqref{SpatialDerivativeProjectionOperator2} and the surface flux functions $\textbf{G}{}_{\iota}^{\*}$, $\iota=1,2,3$, satisfy the generalized Jameson's \cite{Jameson2008} conditions \eqref{Jameson} in Theorem 3.1.
\item[\textbf{(R2)}] The discrete pressure trace approximation $\avg{p}^{\star}$, mass fluxes $\mathbf{G}_{\iota}^{1,*}$ and energy fluxes  $\mathbf{G}_{\iota}^{1,*}$ need to be chosen such that the equation \eqref{LocalDisInternalEnergy} provides a discrete analogue of the equation \eqref{ContinuousInnerEuler}.  
\end{itemize}
In Appendix \ref{NumericalFlux} certain numerical flux functions are listed. By plugging these fluxes in the equation \eqref{LocalDisInternalEnergy}, we obtain the following equations for the discrete internal energy evolution:\footnote{In order to evaluate \eqref{LocalDisInternalEnergy} for the different flux functions, we use the identities:  
\begin{align*}
2\avg{p}\avg{u_{\beta}}-\avg{pu_{\beta}}-\avg{p}u_{\beta}^{-}=\frac{1}{2}p^{-}\jump{u_{\iota}},\qquad\text{and}\qquad\avg{p}\avg{u_{\beta}}-\avg{p}u_{\beta}^{-}=\frac{1}{2}\avg{p}\jump{u_{\beta}},\qquad\beta=1,2,3,
\end{align*}
where $u_{\beta}^{-}$ and $p^{-}$ are interior contributions form the element $e_{\alpha}\left(t\right)$, $\alpha=1,\dots,K$, and the ``$-$'' has to be understood as in \eqref{SurfaceJumpMean}.}   
\begin{description}

\item[Discrete internal energy, Pirozzoli (PI) flux \cite{Pirozzoli2010}:] 
The three dimensional PI-flux is given in Appendix \ref{PI}. This flux provides  the following equation for the discrete internal energy  
\begin{align}\label{LocalDisInternalEnergyPI}
\begin{split}
\pderivative{}{\tau}\left\langle \interpolation{N}{\left(\rho e\right)},\mathcal{J}\right\rangle _{N} 
=&  -\sum_{\iota=1}^{3}\left(\left\langle p,\pderivative{\interpolation{N}{\left(\mathbb{J}\vec{a}^{\iota}\cdot\vec{u}\right)}}{\xi^{\iota}}\right\rangle _{N}-\left\langle S\left(p,\mathbb{J}\vec{a}^{\iota}\right),\vec{u}\right\rangle _{N}\right)\\
& -\sum_{\iota=1}^{3}\sum_{\beta=1}^{3}\ \int\limits _{\partial\mathrm{E},N}\hat{s}n_{\iota}\avg{\rho}\avg{u_{\iota}-\nu_{\iota}}\left[\avg{u_{\beta}^{2}}-\avg{u_{\beta}}^{2}\right]\,dS \\
& -\sum_{\iota=1}^{3}\ \int\limits _{\partial\mathrm{E},N}\hat{s}n_{\iota}\left[\avg{\rho}\avg{\frac{p}{\rho}}-\avg{p}\right]\avg{u_{\iota}}\,dS \\  
& -\sum_{\iota=1}^{3}\ \int\limits _{\partial\mathrm{E},N}\hat{s}n_{\iota}\left[\avg{\rho}\avg{u_{\iota}-\nu_{\iota}}\avg{e}+\frac{1}{2}\avg{p}\jump{u}\right]\,dS. 
\end{split}
\end{align}

We note that in the continuous case for a sufficiently smooth density $\rho$, velocity $\vec{u}=\left[u_{1},u_{2},u_{3}\right]^{T}$ and pressure $p$ the terms 
\begin{align}
& \avg{u_{\beta}^{2}}-\avg{u_{\beta}}^{2}, \qquad \beta =1,2,3, \label{VelocityConservedZero} \\
& \avg{\rho}\avg{\frac{p}{\rho}}-\avg{p} \label{PressureConservedZero}
\end{align}
cancels out to zero. These terms are spectrally small, but they will have a noticeable influence in regions of rapid variation of the solution. 
\item[Discrete internal energy Kennedy and Gruber (KG) flux \cite{Kennedy2008}:] 
The three dimensional KG flux is given in Appendix \ref{KG}. This flux provides the following equation for the discrete internal energy  
\begin{align}\label{LocalDisInternalEnergyKG}
\begin{split}
\pderivative{}{\tau}\left\langle \interpolation{N}{\left(\rho e\right)},\mathcal{J}\right\rangle _{N} 
=&  -\sum_{\iota=1}^{3}\left(\left\langle p,\pderivative{\interpolation{N}{\left(\mathbb{J}\vec{a}^{\iota}\cdot\vec{u}\right)}}{\xi^{\iota}}\right\rangle _{N}-\left\langle S\left(p,\mathbb{J}\vec{a}^{\iota}\right),\vec{u}\right\rangle _{N}\right)\\
& -\sum_{\iota=1}^{3}\sum_{\beta=1}^{3}\ \int\limits _{\partial\mathrm{E},N}\hat{s}n_{\iota}\avg{\rho}\avg{u_{\iota}-\nu_{\iota}}\left[\avg{u_{\beta}^{2}}-\avg{u_{\beta}}^{2}\right]\,dS \\
&-\sum_{\iota=1}^{3}\ \int\limits _{\partial\mathrm{E},N}\hat{s}n_{\iota}\left[\avg{\rho}\avg{u_{\iota}-\nu_{\iota}}\avg{e}+\frac{1}{2}\avg{p}\jump{u}\right]\,dS. 
\end{split}
\end{align}

\item[Discrete internal energy, Kuya, Totani and Kawai (KTK) flux \cite{kuya2018kinetic}:]
The three dimensional KTK-flux is given in Appendix \ref{Kuya}. This flux provides the following equation for the discrete internal energy  
\begin{align}\label{LocalDisInternalEnergyKTK}
\begin{split}
\pderivative{}{\tau}\left\langle \interpolation{N}{\left(\rho e\right)},\mathcal{J}\right\rangle _{N} 
=& -\sum_{\iota=1}^{3}\left(\left\langle p,\pderivative{\interpolation{N}{\left(\mathbb{J}\vec{a}^{\iota}\cdot\vec{u}\right)}}{\xi^{\iota}}\right\rangle _{N}-\left\langle S\left(p,\mathbb{J}\vec{a}^{\iota}\right),\vec{u}\right\rangle _{N}\right) \\
&-\sum_{\iota=1}^{3}\ \int\limits _{\partial\mathrm{E},N}\hat{s}n_{\iota}\left[\avg{\rho}\avg{u_{\iota}-\nu_{\iota}}\avg{e}+\frac{1}{2}p^{-}\jump{u}\right]\,dS,
\end{split}
\end{align}
where $p^{-}$ is the interior pressure contribution form the element $e_{\alpha}\left(t\right)$, $\alpha=1,\dots,K$. The exterior pressure contribution $p^{+}$ does not appear in the discrete integral along the element interfaces in the discrete internal energy \eqref{LocalDisInternalEnergyKTK} like in the equation \eqref{LocalDisInternalEnergyKG} for the KG-flux. In certain situations, this has an impact on the preservation of the discrete energy ratio. Our numerical experiments in Figure~\ref{fig:TGV_EC_KEP_Euler} in Section \ref{EC_KEP_Test} support this statement. However, by replacing the term $2\avg{p}\avg{u_{\iota}}-\avg{pu_{\iota}}$ with $\avg{p}\avg{u_{\iota}}$, $\iota=1,2,3$, in the energy contribution 
($\mathbf{G}{}_{\iota}^{5,\text{KTK}}$) of the KTK-flux \eqref{KTK_Flux}, we obtain the modified Kuya, Totani and Kawai (M\_KTK) flux 
\begin{align}\label{HF_Flux}
\begin{split} 
\mathbf{G}{}_{\iota}^{1,\text{M\_KTK}\phantom{+1,}}=&\phantom{+} \avg{\rho}\avg{u_{\iota}-\nu_{\iota}}, \\
\mathbf{G}{}_{\iota}^{\upsilon+1,\text{M\_KTK}}=& \phantom{+}
\avg{\rho}\avg{u_{\iota}-\nu_{\iota}}\avg{u_{\upsilon}}+\avg{p}\delta_{\iota\upsilon},\qquad\upsilon=1,2,3, \\
\mathbf{G}{}_{\iota}^{5,\text{M\_KTK}\phantom{+1,}}=&\phantom{+}  
\avg{\rho}\avg{u_{\iota}-\nu_{\iota}}\left[\avg{e}+\frac{1}{2}\overline{\left|\vec{u}\right|}^{2}\right]+\avg{p}\avg{u_{\iota}}.
\end{split}
\end{align}  
For this flux the discrete internal energy becomes    
\begin{align}\label{LocalDisInternalEnergyHF}
\begin{split}
\pderivative{}{\tau}\left\langle \interpolation{N}{\left(\rho e\right)},\mathcal{J}\right\rangle _{N} 
=& -\sum_{\iota=1}^{3}\left(\left\langle p,\pderivative{\interpolation{N}{\left(\mathbb{J}\vec{a}^{\iota}\cdot\vec{u}\right)}}{\xi^{\iota}}\right\rangle _{N}-\left\langle S\left(p,\mathbb{J}\vec{a}^{\iota}\right),\vec{u}\right\rangle _{N}\right) \\
& -\sum_{\iota=1}^{3}\ \int\limits _{\partial\mathrm{E},N}\hat{s}n_{\iota}\left[\avg{\rho}\avg{u_{\iota}-\nu_{\iota}}\avg{e}+\frac{1}{2}\avg{p}\jump{u}\right]\,dS.
\end{split}
\end{align}
Here the quantity \eqref{VelocityConservedZero} and the pressure contribution are the same as for the KG-flux in the discrete integral along the element interfaces. The numerical experiments in Section \ref{EC_KEP_Test} show that this slight modification in the KTK-flux \eqref{KTK_Flux} provides a meaningful preservation of the discrete energy ratio.  

\item[Discrete internal energy, Ranocha (RA) flux \cite{Ranocha2018}:]
In order to work with the RA-flux, we introduce the logarithmic average 
\begin{equation}\label{LogMean}
\avg{a}^{\text{log}}:=\begin{cases}
\frac{\jump{a}}{\jump{\log\left(a\right)}}, & \text{if }a^{-}\neq a^{+},\\
\avg{a}, & \text{if }a^{-}=a^{+}
\end{cases}
\end{equation} 
for a state $a$ with positive right limit $a^{+}>0$ and positive left limit $a^{-}>0$. 
The three dimensional RA-flux is given in Appendix \ref{RA}. We note that this flux is EC, since it satisfies Tadmor's \cite{Tadmor2003}  discrete entropy conditions such that 
\begin{equation}\label{Tadmor1}
\jump{\mathbf{w}}^{T}\mathbf{G}{}_{\iota}^{\text{RA}}=\jump{\rho u_{\iota}}-\avg{\nu_{\iota}}\jump{\rho},\qquad\iota=1,2,3.
\end{equation} 
Furthermore, it gives the following equation for the discrete internal energy  
\begin{align}\label{LocalDisInternalEnergyRA}
\begin{split}
\pderivative{}{\tau}\left\langle \interpolation{N}{\left(\rho e\right)},\mathcal{J}\right\rangle _{N} 
=& -\sum_{\iota=1}^{3}\left(\left\langle p,\pderivative{\interpolation{N}{\left(\mathbb{J}\vec{a}^{\iota}\cdot\vec{u}\right)}}{\xi^{\iota}}\right\rangle _{N}-\left\langle S\left(p,\mathbb{J}\vec{a}^{\iota}\right),\vec{u}\right\rangle _{N}\right)\\
&-\sum_{\iota=1}^{3}\ \int\limits _{\partial\mathrm{E},N}\hat{s}n_{\iota}\left[\avg{\rho}^{\text{log}}\avg{u_{\iota}-\nu_{\iota}}\frac{1}{\avg{\frac{1}{e}}^{\text{log}}}-\frac{1}{2}p^{-}\jump{u}\right]\,dS,
\end{split}
\end{align}
where $p^{-}$ is the interior pressure contribution form the element $e_{\alpha}\left(t\right)$, $\alpha=1,\dots,K$. Here the same issue as in the discrete internal energy equation \eqref{LocalDisInternalEnergyKTK} for the KTK-flux appears. The exterior pressure contribution $p^{+}$ does not appear in the discrete integral along the element interfaces. This seems to negatively impact the preservation of the discrete energy.

\item[Discrete internal energy, Chandrashekar (CH) flux \cite{Chandrashekar2013}:] The three dimensional CH-flux is given in Appendix \ref{CH}. The CH-flux is an EC one, since it satisfies the equation \eqref{Tadmor1}, too. However, the restrictions for the volume flux in Theorem 3.1 are not satisfied for this flux, since the pressure is computed by $\avg{\rho}/\avg{\frac{\rho}{p}}$  in the momentum contribution 
($\mathbf{G}{}_{\iota}^{\upsilon+1,\text{CH}}$, $\iota, \upsilon=1,2,3$). Thus, it cannot be used in the equation \eqref{LocalDisInternalEnergy} to compute the discrete internal energy. If we replace the term $\avg{\rho}/\avg{\frac{\rho}{p}}$  by $\avg{p}$ in the momentum contribution, we obtain the flux   
\begin{align}\label{CH_AVG_Flux}
\begin{split} 
\textbf{G}{}_{\iota}^{1,\text{M\_CH}\phantom{+1,}}=& \phantom{+}
\avg{\rho}^{\text{log}}\avg{u_{\iota}-\nu_{\iota}}, \\
\textbf{G}{}_{\iota}^{\upsilon+1,\text{M\_CH}}=& \phantom{+}
\avg{\rho}^{\text{log}}\avg{u_{\iota}-\nu_{\iota}}\avg{u_{\upsilon}}+\avg{p}\delta_{\iota\upsilon},\qquad\upsilon=1,2,3, \\
\textbf{G}{}_{\iota}^{5,\text{M\_CH}\phantom{+1,}}=& \phantom{+}
\avg{\rho}^{\text{log}}\left[\frac{1}{\avg{\frac{1}{e}}^{\text{log}}}+\frac{1}{2}\overline{\left|\vec{u}\right|}^{2}\right]\avg{u_{\iota}-\nu_{\iota}}+\frac{\avg{\rho}}{\avg{\frac{\rho}{p}}}\avg{u_{\iota}}.
\end{split}
\end{align}
The flux \eqref{CH_AVG_Flux} does not satisfy the entropy condition \eqref{Tadmor1}, but the restrictions of Theorem 3.1. For this flux the discrete internal energy is given by 
\begin{align}\label{LocalDisInternalEnergyMCH}
\begin{split}
\pderivative{}{\tau}\left\langle \interpolation{N}{\left(\rho e\right)},\mathcal{J}\right\rangle _{N} 
=& -\sum_{\iota=1}^{3}\left(\left\langle p,\pderivative{\interpolation{N}{\left(\mathbb{J}\vec{a}^{\iota}\cdot\vec{u}\right)}}{\xi^{\iota}}\right\rangle _{N}-\left\langle S\left(p,\mathbb{J}\vec{a}^{\iota}\right),\vec{u}\right\rangle _{N}\right)\\
&-\sum_{\iota=1}^{3}\ \int\limits _{\partial\mathrm{E},N}\hat{s}n_{\iota}\left[\frac{\avg{\rho}}{\avg{\frac{\rho}{p}}}-\avg{p}\right]\avg{u_{\iota}}\,dS\\
&-\sum_{\iota=1}^{3}\ \int\limits _{\partial\mathrm{E},N}\hat{s}n_{\iota}\left[\avg{\rho}^{\text{log}}\avg{u_{\iota}-\nu_{\iota}}\frac{1}{\avg{\frac{1}{e}}^{\text{log}}}-\frac{1}{2}p^{-}\jump{u}\right]\,dS.
\end{split}
\end{align}
We note that in the continuous case for a sufficiently smooth density $\rho$,  and pressure $p$ the terms 
\begin{align}
\frac{\avg{\rho}}{\avg{\frac{\rho}{p}}}-\avg{p} \label{PressureConservedZero1}
\end{align}
cancel out to zero. This term is spectrally small, but it will be noticeable in regions of rapid variation of the solution. In Section \ref{EC_KEP_Test}, we investigate the CH-Flux and the modified flux \eqref{CH_AVG_Flux} for the inviscid TGV test case on a moving curved mesh. The results are given in Figure~\ref{fig:TGV_EC_KEP_Euler}, it can be seen that the the total kinetic energy is preserved for the modified flux \eqref{CH_AVG_Flux} and decreases for the CH-Flux. For the total entropy the contrary behavior is observed. The same behavior of the total kinetic energy and entropy was observed in numerical experiments on a static mesh in \cite{Gassner2016}.    
\end{description}


\subsection{Split form ALE DGSEM for the NSE}\label{sec:DGSEMNSE}
In this section, we  extend the split form ALE DGSEM for the Euler equations \eqref{MovingMeshDGSEMEuler} to solve the NSE \eqref{eq:NSELawBlock} on moving curved elements with the approximation of the viscous terms from \cite{Gassner2017}.  

The unknowns of the transformed NSE \eqref{eq:NSELawBlockRefE} $J$, $\textbf{u}$ and $\blockvec{\textbf{q}}$ on the time-independent reference element $\mathrm{E}$ are approximated by $\mathcal{J}$, $\textbf{U}$, $\blockvec{\textbf{Q}}$. These quantiles are given by \eqref{UnkonwnNSE}. The fluxes and the viscous block matrix are approximated by the interpolation operator \eqref{InterpolationOperator}. The approximation of the equation \eqref{GCL} is the same as for the Euler equations (see Section \ref{sec:DGSEMEuler}) and leads to the equation \eqref{DGSEM:D-GLC}. The equation \eqref{eq:NSE} is also treated as for the Euler equations and leads to the equation 
\begin{align}\label{DGSEM:NSE1} 
\begin{split}
\left\langle \pderivative{\left(\mathcal{J}\textbf{U}\right)}{\tau},\boldsymbol{\varphi}\right\rangle _{N}
=&-\left\langle \Dprojection{N}\cdot\blockvec{\tilde{\textbf{G}}}{}^{\#},\boldsymbol{\varphi}\right\rangle _{N}-\int\limits _{\partial\mathrm{E},N}\boldsymbol{\varphi}^{T}\left(\tilde{\textbf{G}}_{\hat{n}}^{*}-\tilde{\textbf{G}}_{\hat{n}}\right)\,dS \\ 
&
+\left\langle \vec{\nabla}_{\xi}\cdot\interpolation{N}{\left(\blockmatx{\tilde{B}}^{\,\text{v}}\left(\mathbf{W}\right)\blockvec{\textbf{Q}}\right)},\mathbf{\boldsymbol{\varphi}}\right\rangle _{N}
+
\int\limits _{\partial E,N}\boldsymbol{\varphi}^{T}\left\{ \left[\avg{\blockmatx{\tilde{B}}^{\,\text{v}}\left(\textbf{W}\right)\blockvec{\textbf{Q}}}-\blockmatx{\tilde{B}}^{\,\text{v}}\left(\textbf{W}\right)\blockvec{\textbf{Q}}\right]\cdot\hat{n}\right\} \,dS, 
\end{split}
\end{align}
for all test functions $\boldsymbol{\varphi}\in\mathbb{P}^{N}\left(\mathrm{E},\R^{5}\right)$. The derivative operator in \eqref{DGSEM:NSE1} is again computed by \eqref{SpatialDerivativeProjectionOperator2}, and the flux $\tilde{\textbf{G}}_{\hat{n}}$ is computed by \eqref{BlockvectorNormal}. Furthermore, for $i,j,k=0,\dots,N$  
\begin{equation}\label{DiscreteEntropyVariables}
\mathbf{W}_{ijk}:=\left[\interpolation{N}{\left(\mathbf{w}\right)}\right]_{ijk}=\left[\frac{\gamma-\varsigma_{ijk}}{\gamma-1}-\frac{\rho_{ijk}}{2p_{ijk}}\left|\vec{u}\right|^{2},\frac{\rho_{ijk}\left(\vec{u}\right)_{ijk}}{p_{ijk}},-\frac{\rho_{ijk}}{p_{ijk}}\right]^{T}.
\end{equation}
Next, we multiply the equation \eqref{eq:ViscousPart} with test functions 
\begin{equation}\label{BlockTestFunction}
\blockvec{\boldsymbol{\varphi}}=\begin{bmatrix}\boldsymbol{\varphi}_{1}^{T},\boldsymbol{\varphi}_{2}^{T},\boldsymbol{\varphi}_{3}^{T}\end{bmatrix}^{T}
\qquad
\boldsymbol{\varphi}{}_{i}\in\mathbb{P}^{N}\left(\mathrm{E},\R^{5}\right), \qquad i=1,2,3,
\end{equation}
integrate the resulting equations and use integration-by-parts to separate boundary and volume contributions. The volume integrals in the variational form are approximated with the LGL quadrature. Then, we insert numerical surface fluxes to approximate 
$\blockmatx{\tilde{B}}^{\,\text{v}}\left(\textbf{W}\right)\blockvec{\textbf{Q}}$ at the spatial element interfaces. Finally, the SBP property \eqref{SBP} is used for the volume contribution to get the strong form: 
\begin{align}\label{DGSEM:NSE2}
\begin{split}
\left\langle \mathcal{J}\blockvec{\mathbf{Q}},\blockvec{\boldsymbol{\varphi}}\right\rangle _{N}= & \phantom{-} \left\langle \vec{\nabla}_{\xi}\interpolation{N}{\left(\mathbf{W}\right)},\blockmatx{M}^{T}\blockvec{\boldsymbol{\varphi}}\right\rangle _{N}+\int\limits _{\partial E,N}\left(\avg{\mathbf{W}}-\mathbf{W}\right)^{T}\left\{ \blockmatx{M}^{T}\blockvec{\boldsymbol{\varphi}}\cdot\hat{n}\right\} \,dS,  
\end{split}
\end{align}
for all test functions \eqref{BlockTestFunction}. Thus, the split form ALE DGSEM for the NSE is given by the equations \eqref{DGSEM:D-GLC}, \eqref{DGSEM:NSE1} and \eqref{DGSEM:NSE2}. 


\subsection{Discrete entropy stability}\label{Sec:EntropyStableDGSEM}
The NSE equations are equipped with the entropy/entropy flux pair \eqref{EntropyNSE} and satisfy the entropy inequality \eqref{ContinuousEntropy3}. \\
\noindent   
In the following, we assume that the volume flux functions $\tilde{\textbf{G}}{}_{\iota}^{\#}$, $\iota=1,2,3$, in the derivative projection operator \eqref{SpatialDerivativeProjectionOperator2} satisfy in the LGL points for $i,j,k,m=0,\dots,N$, the discrete entropy conditions: 
\begin{align}\label{EntropyConditionDGECFlux}
\begin{split}
\jump{\textbf{W}}_{(i,m)jk}^{T}\textbf{G}_{1}^{\#}\left(\vec{\nu}_{ijk},\vec{\nu}_{mjk},\textbf{U}_{ijk},\textbf{U}_{mjk}\right)=\jump{\rho\vec{u}}_{(i,m)jk}-\avg{\vec{\nu}}_{\left(i,m\right)jk}\jump{\rho}_{(i,m)jk},\\ 
\jump{\textbf{W}}_{i(j,m)k}^{T}\textbf{G}_{2}^{\#}\left(\vec{\nu}_{ijk},\vec{\nu}_{imk},\textbf{U}_{ijk},\textbf{U}_{imk}\right)=\jump{\rho\vec{u}}_{i(j,m)k}-\avg{\vec{\nu}}_{i\left(j,m\right)k}\jump{\rho}_{i(j,m)k}, \\
\jump{\textbf{W}}_{ij(k,m)}^{T}\textbf{G}_{3}^{\#}\left(\vec{\nu}_{ijk},\vec{\nu}_{ijm},\textbf{U}_{ijk},\textbf{U}_{ijm}\right)=\jump{\rho\vec{u}}_{ij(k,m)}-\avg{\vec{\nu}}_{ij\left(k,m\right)}\jump{\rho}_{ij(k,m)}.
\end{split}
\end{align}
Moreover, we assume that the Cartesian surface flux functions $\tilde{\textbf{G}}{}_{\iota}^{*}$, $\iota=1,2,3$, satisfy the discrete entropy conditions: 
\begin{equation}\label{Tadmor0}
\jump{\mathbf{w}}^{T}\mathbf{G}{}_{\iota}^{*}\leq\jump{\rho u_{\iota}}-\avg{\nu_{\iota}}\jump{\rho},\qquad\iota=1,2,3.
\end{equation} 
We note that the flux functions in the Appendices \ref{RA}, \ref{CH}, satisfies the conditions \eqref{EntropyConditionDGECFlux}. A flux with the property \eqref{Tadmor0} can be constructed from the fluxes in the Appendices \ref{RA}, \ref{CH}, and a suitable dissipation operator $\matx{H}_{\iota}$ \eqref{EntropyBasedDissipation}.   

In \cite[Section 3.6]{Schnuecke2018}, it has been proven that the restrictions \eqref{EntropyConditionDGECFlux} and \eqref{Tadmor0} lead to an ES ALE DGSEM \eqref{MovingMeshDGSEMEuler}. The following two identities with the entropy/entropy flux pair \eqref{EntropyNSE} are a direct consequence from the proof in \cite{Schnuecke2018}: 
\begin{align}\label{DGSEM:EntropyPreservation}
\begin{split}
\left\langle \pderivative{\left(\mathcal{J}\textbf{U}\right)}{\tau},\mathbf{W}\right\rangle _{N} 
=& \phantom{+} \pderivative{}{\tau}\left\langle \interpolation{N}{\left(s\right)},J\right\rangle _{N}+\left\langle \Dprojection{N}\cdot\vec{\tilde{\nu}},\rho\right\rangle _{N} +\int\limits _{\partial E,N}\left(\tilde{\nu}_{\hat{n}}^{*}-\tilde{\nu}_{\hat{n}}\right)\rho\,dS,
\end{split}
\end{align}
\begin{equation}\label{DGSEM_VolumeCont}
\left\langle \Dprojection{N}\cdot\blockvec{\tilde{\textbf{G}}}{}^{\#},\mathbf{W}\right\rangle _{N} =
\int\limits _{\partial E,N}\left(\tilde{f}_{\hat{n}}^{s}-\tilde{\nu}_{\hat{n}}s\right)\,dS-\left\langle \Dprojection{N}\cdot\vec{\tilde{\nu}},\rho\right\rangle _{N}, 
\end{equation}
where $\textbf{W}$ is given by \eqref{DiscreteEntropyVariables} and 
\begin{equation}
\tilde{f}_{\hat{n}}^{s}=\left(\hat{s}\vec{n}\right)\cdot\vec{f}^{s},\qquad\vec{f}^{s}=\left[f_{1}^{s},f_{2}^{s},f_{3}^{s}\right]^{T}.
\end{equation} 
On the other hand in \cite[Section 4]{Gassner2017}, it has been proven that the original Bassi and Rebay \cite{Bassi1997} scheme (BR1) for the discretization of second order viscous terms on static meshes is stable in the context of the DGSEM approximation with LGL points. By the same analysis as in \cite[Section 4.2.2]{Gassner2017}, the following identity can be proven   
\begin{equation}\label{DGSEM_SurfCont}
\left\langle \vec{\nabla}_{\xi}\cdot\interpolation{N}{\left(\blockmatx{\tilde{B}}^{\,\text{v}}\left(\mathbf{W}\right)\blockvec{\mathbf{Q}}\right)},\mathbf{W}^{T}\right\rangle _{N}\leq\int\limits _{\partial E,N}\avg{\mathbf{W}}^{T}\left\{ \left(\blockmatx{\tilde{B}}^{\,\text{v}}\left(\mathbf{W}\right)\blockvec{\mathbf{Q}}\right)\cdot\hat{n}\right\} \,dS.
\end{equation}
Therefore, we obtain the identity 
\begin{align}\label{DGSEM:EntropyInequality}
\begin{split}
\pderivative{}{\tau}\left\langle \interpolation{N}{\left(s\right)},J\right\rangle _{N}
\leq& -\int\limits _{\partial E,N}\left(\rho\left[\tilde{u}_{\hat{n}}-\tilde{\nu}_{\hat{n}}^{*}\right]-\mathbf{W}^{T}\tilde{\mathbf{G}}_{\hat{n}}^{*}\right)\,dS\\ 
&+\int\limits _{\partial E,N}\mathbf{W}^{T}\left\{ \left[\avg{\blockmatx{\tilde{B}}^{\,\text{v}}\left(\textbf{W}\right)\blockvec{\textbf{Q}}}-\blockmatx{\tilde{B}}^{\,\text{v}}\left(\textbf{W}\right)\blockvec{\textbf{Q}}\right]\cdot\hat{n}\right\} \,dS \\
& +\int\limits _{\partial E,N}\avg{\mathbf{W}}^{T}\left\{ \left(\blockmatx{\tilde{B}}^{\,\text{v}}\left(\mathbf{W}\right)\blockvec{\mathbf{Q}}\right)\cdot\hat{n}\right\} \,dS
\end{split}
\end{align}
by using $\textbf{W}$ (see \eqref{DiscreteEntropyVariables}) as test function in equation \eqref{DGSEM:NSE2} and the identities \eqref{DGSEM:EntropyPreservation}, \eqref{DGSEM_VolumeCont}  and \eqref{DGSEM_SurfCont}. A summation of this equation provides a discrete analogue of the equation \eqref{ContinuousEntropy3} in the context of the spectral element approximation when the NSE is investigated with periodic boundary conditions.  
\begin{rem-hand}[3.3.]
In a similar way a KED DGSEM can be constructed for the NSE. Then, the flux functions $\textbf{G}{}_{\iota}^{\#}$, $\iota=1,2,3$, in the derivative projection operator \eqref{SpatialDerivativeProjectionOperator2} and the surface flux functions $\textbf{G}{}_{\iota}^{*}$, $\iota=1,2,3$, should satisfy the conditions \eqref{CartesianSurfaceFlux}.  
\end{rem-hand}

\section{Numerical results}\label{sec:numResults}

The proposed split form ALE DGSEM is implemented in the open source high order DG 
solver FLEXI\footnote{http://www.flexi-project.org}~\cite{krais2019flexi}.  It 
provides the necessary framework for the implementation of different split forms for 
high order unstructured meshes, was successfully applied to under-resolved 
simulations in fluid dynamics before~\cite{beck2014high,flad2016simulation} and 
shows excellent scaling properties, making it a suitable choice for large-scale
simulations, as presented in the next chapter.

\subsection{Experimental convergence rates}

In this section, the convergence behavior under mesh refinement of the split form ALE DGSEM is assessed for different flux functions
using the method of manufactured solutions. The results are used to verify that the proposed methodology retains its high order
accuracy on moving grids. For the following simulations, we assume a solution of the form
\begin{align}\label{MFS:Euler3D}
\begin{split}
\rho\left(\vec{x},t\right)
=& \quad 2+0.1\sin\left(\pi(x_{1}+x_{2}+x_{3}-0.6t\right)), \\ 
\rho u_{\iota}\left(\vec{x},t\right)
=& \quad 2+0.1\sin\left(\pi(x_{1}+x_{2}+x_{3}-0.6t)\right),\qquad\iota=1,2,3, 
\\
E\left(\vec{x},t\right)=& \quad\left[2+0.1\sin\left(\pi(x_{1}+x_{2}+x_{3}-0.6t)\right)\right]^{2},
\end{split}
\end{align}
with total energy $E=e+\frac{1}{2}\rho\left|\vec{u}\right|^{2}$ on the domain $[-1,1]^3$ and compute the residual when \eqref{MFS:Euler3D} is inserted into the Euler equations. The resulting
terms are then used as sources for the simulations, and are discretized as a solution independent part of the computation. 

\begin{figure}
    \centering
\begin{tikzpicture}
\def\xsize{8cm}

\begin{scope}
\begin{axis}[
legend cell align={left},
legend style={fill opacity=0.8, draw opacity=1, text opacity=1, at={(0.03,0.97)}, anchor=north west, draw=white!80.0!black,nodes={scale=0.8, transform shape}},
log basis x={10},
log basis y={10},
tick align=outside,
tick pos=left,
x grid style={white!69.01960784313725!black},
xmin=0.0544094102060077, xmax=1.14869835499704,
xmode=log,
xtick style={color=black},
y grid style={white!69.01960784313725!black},
ymin=1.58379684144504e-8, ymax=0.829770434474778,
ymode=log,
ytick style={color=black},
title={$N=3$},
xlabel={$h$},
ylabel={$L^2(\rho)$},
width=0.47\linewidth
]
\addplot
table {%
1                  0.32472134156961E-00001
0.5                0.10372984686567E-00001
0.25               0.59876439114993E-00003
0.125              0.24792105688542E-00004
0.0833333333333333 0.13818857294637E-00005
};
\addlegendentry{M\textunderscore KTK}
\addplot
table {%
1                  0.49178693068282E-00001
0.5                0.10753581246525E-00001
0.25               0.71020383307959E-00003
0.125              0.24950520574466E-00004
0.0833333333333333 0.13284204686632E-00005
};
\addlegendentry{CH}
\addplot
table {%
1                  0.32614776075428E-00001
0.5                0.10386132468305E-00001
0.25               0.59871696519161E-00003
0.125              0.24792104504906E-00004
0.0833333333333333 0.13818852164842E-00005
};
\addlegendentry{KG}
\addplot
table {%
1                  0.32551381651132E-00001
0.5                0.10338460937626E-00001
0.25               0.59934247385576E-00003
0.125              0.24696772068622E-00004
0.0833333333333333 0.13761242786106E-00005
};
\addlegendentry{KTK}
\addplot
table {%
1                  0.33030668904638E-00001
0.5                0.10450579167904E-00001
0.25               0.60520608824025E-00003
0.125              0.25325548080013E-00004
0.0833333333333333 0.14178742176779E-00005
};
\addlegendentry{PI}
\addplot
table {%
1                  0.50410585613431E-00001
0.5                0.10812152413767E-00001
0.25               0.73820729060528E-00003
0.125              0.32935331561974E-00004
0.0833333333333333 0.15934720443680E-00005
};
\addlegendentry{RA}
\logLogSlopeTriangle{0.62}{0.25}{0.45}{4}{black};
\end{axis}
\end{scope}

\begin{scope}[xshift=\xsize]
\begin{axis}[
legend cell align={left},
legend style={fill opacity=0.8, draw opacity=1, text opacity=1, at={(0.03,0.97)}, anchor=north west, draw=white!80.0!black,nodes={scale=0.8, transform shape}},
log basis x={10},
log basis y={10},
tick align=outside,
tick pos=left,
x grid style={white!69.01960784313725!black},
xmin=0.0544094102060077, xmax=1.14869835499704,
xmode=log,
xtick style={color=black},
y grid style={white!69.01960784313725!black},
ymin=1.58379684144504e-8, ymax=0.829770434474778,
ymode=log,
ytick style={color=black},
title={$N=4$},
xlabel={$h$},
ylabel={$L^2(\rho)$},
width=0.47\linewidth
]
\addplot
table {%
1                  0.20000576143502E-00001
0.5                0.15271308255255E-00002
0.25               0.38781457425397E-00004
0.125              0.10783531009122E-00005
0.0833333333333333 0.41263034994883E-00007
};
\addlegendentry{M\textunderscore KTK}
\addplot
table {%
1                  0.17800989516228E-00001
0.5                0.16545124153948E-00002
0.25               0.40008739154108E-00004
0.125              0.10318676045724E-00005
0.0833333333333333 0.35895265834119E-00007
};
\addlegendentry{CH}
\addplot
table {%
1                  0.19998010260730E-00001
0.5                0.15271665274415E-00002
0.25               0.38781581949539E-00004
0.125              0.10783527524152E-00005
0.0833333333333333 0.41263034274154E-00007
};
\addlegendentry{KG}
\addplot
table {%
1                  0.19973218716477E-00001
0.5                0.15229452748103E-00002
0.25               0.38647024449703E-00004
0.125              0.10744181952269E-00005
0.0833333333333333 0.41073808982333E-00007
};
\addlegendentry{KTK}
\addplot
table {%
1                  0.19715751906352E-00001
0.5                0.15293854249747E-00002
0.25               0.38681803611928E-00004
0.125              0.11597839770542E-00005
0.0833333333333333 0.44002810305127E-00007
};
\addlegendentry{PI}
\addplot
table {%
1                  0.17649412117726E-00001
0.5                0.17286729935944E-00002
0.25               0.49586657286233E-00004
0.125              0.18933920475610E-00005
0.0833333333333333 0.54445518200573E-00007
};
\addlegendentry{RA}
\logLogSlopeTriangle{0.62}{0.25}{0.30}{5}{black};
\end{axis}
\end{scope}

\end{tikzpicture}
    \caption{Experimental convergence rates for different flux functions. Manufactured solution test case with standing wave mesh
movement, for both $N=3$ and $N=4$. Shown are the $L^2$ error norms of the density $\rho$ over the mesh size $h$.}
    \label{fig:hConv}
\end{figure}
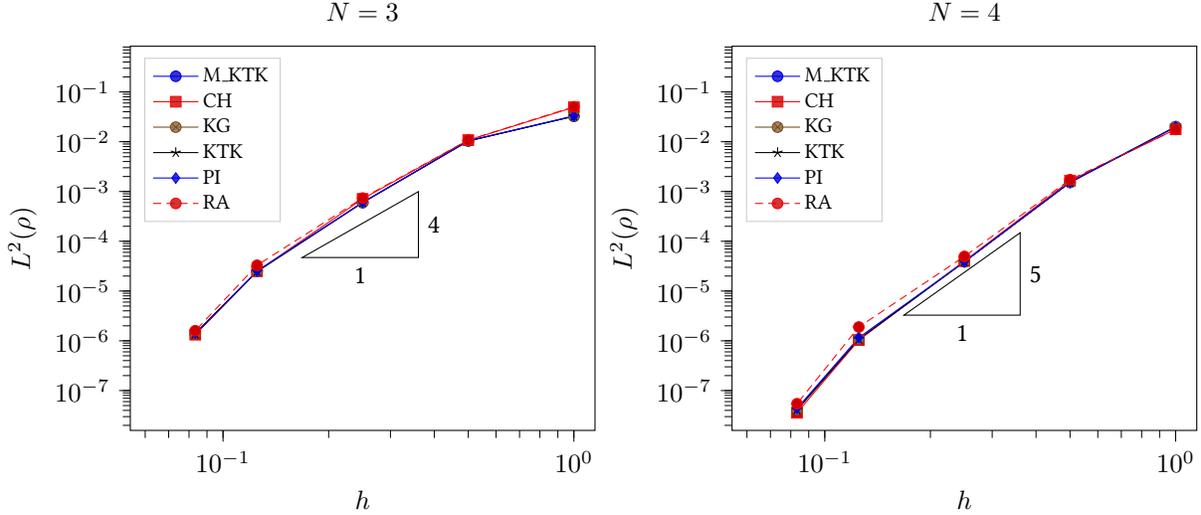

The simulations are performed on initially Cartesian grids with an increasing number of elements 
\begin{equation}
K=2^{3},4^{3},8^{3},16^{3},32^{3}.
\end{equation}
Since we are interested in the behavior for high order meshes, we represent the boundary curves of the elements with polynomials of
degree $2$, denoted as $N_{geo}=2$, for all test cases. All meshes are undergoing a forced periodic motion, corresponding to a
standing wave. The position $\vec{x}$ of a grid point at time $t$ can be described by
the equation
\begin{equation}\label{GridpointDistribution3D}
\vec{x}\left(t\right)=\vec{x}\left(0\right)+0.05 L \sin\left(2\pi t\right)\sin\left(\frac{2\pi}{L}x_{1}\left(0\right)\right)\sin\left(\frac{2\pi}{L}x_{2}\left(0\right)\right)\sin\left(\frac{2\pi}{L}x_{3}\left(0\right)\right),
\end{equation}
where $\vec{x}(t=0)$ is the position of the grid point in the non-deformed configuration of the mesh, and $L$ is the length of the
domain. Note that \eqref{GridpointDistribution3D} provides a time-dependent domain which is divided in $K$ non-overlapping elements in each time point $t$. Then each element is mapped on the reference element. The mappings are element local polynomials that generate a watertight mesh. Thus, the mapping is piecewise polynomial as is common for DG approximations. Furthermore, we note that the mesh velocities are computed by exact differentiation of the above equation \eqref{GridpointDistribution3D}. The simulation is advanced until final time $T=5$, and we employ the five stage fourth order low-storage explicit RK method (RK4(3)5[2R+]) from Kennedy, Carpenter and Lewis \cite{kennedy2000low} for time-integration. Since the stability region of the explicit time integration scheme is
restricted by the Courant-Friedrichs-Lewy (CFL) condition, the allowable time step $\Delta t$ is computed as in \cite{Cockburn2001} 
\begin{equation}\label{CFL1}
\frac{\Delta t}{\underset{1\leq\alpha\leq K}{\min}\left|h_{\alpha}\left(t^{n}\right)\right|}\leq\frac{C_{\text{CFL}}}{\left(2N+1\right)\lambda_{\text{max}}},
\end{equation}
where $h_{\alpha}(t)$ is the size of element $e_\alpha(t)$, $\lambda_{max}$ the fastest signal velocity of the Euler equations and
$C_{\text{CFL}}$ is set to $C_{\text{CFL}}=0.1$. The surface fluxes are used with Roe-type dissipation terms \eqref{EntropyBasedDissipation}, which can be found in \cite[Appendix C.3]{Schnuecke2018}. Figure~\ref{fig:hConv} shows the results obtained for the manufactured solution test cases. The quantity used for comparison is the
$L^2$ norm of the error in the density $\rho$, when compared with the manufactured solution \eqref{MFS:Euler3D}. To compute the
integrals required for the $L^2$ norm, we employ Gauss quadrature on a supersampled version of the solution field, with $M=13$
integration points per direction. We show results for both $N=3$ and $N=4$ and a variety of numerical flux functions. As for smooth solutions the difference in the split forms is spectrally small, it is not surprising that all the different formulations lead to very similar results for these well resolved convergence test cases. All of them retain the design order of accuracy, confirming the high order approximation property of the split form ALE DGSEM.

\subsection{Numerical validation of the entropy and kinetic energy analysis}\label{EC_KEP_Test}

In this section, the KEP and EC property of the different split forms (choice of 
volume flux functions) is assessed. As a test case, we choose the Taylor-Green 
vortex (TGV) \cite{shu2005numerical}, which is an important canonical problem for 
laminar-turbulent transition and turbulent flows. The flow is
defined on the domain $\Omega=\left[0,2\pi\right]^{3}$ with periodic boundary conditions on all sides, and its initial conditions
are given by 
\begin{align}\label{TGVInitial}
\begin{split}
\rho\left(\vec{x},0\right)=&\quad 1, \\ 
\vec{u}\left(\vec{x},0\right)=& \quad\left[\text{sin}\left(x_{1}\right)\text{cos}\left(x_{2}\right)\text{cos}\left(x_{3}\right),-\text{cos}\left(x_{1}\right)\text{sin}\left(x_{2}\right)\text{cos}\left(x_{3}\right),0\right]^{T}, \\
p\left(\vec{x},0\right)=&\quad p_{0}+\frac{1}{16}\left(\text{cos}\left(2x_{1}\right)+\text{cos}\left(2x_{2}\right)\right)\left(\text{cos}\left(2x_{3}\right)+2\right), 
\end{split}
\end{align}
where for all following simulations, $p_{0}$ is chosen such that the Mach number becomes $\textrm{Ma}=\frac{1}{\sqrt{\gamma p_{0}}}=0.1$, and the influence of
compressibility is relatively small. The flow is evolved until final time $T=13$, which is past the critical point for stability reached at
$t \approx 9$, where the maximum of turbulent dissipation occurs. The mesh is again initially Cartesian with $K=16^3$ elements,
$N_{geo}=2$ and forced to perform the periodic motion given by \eqref{GridpointDistribution3D}. The constant $C_{\text{CFL}}$ in \eqref{CFL1} is set to $C_{\text{CFL}}=0.5$, and again the RK4(3)5[2R+] explicit time stepping scheme is employed.

We chose the TGV as a test case since (i) it has simple periodic boundary conditions, removing the influence of e.g. walls,
(ii) is challenging for the stability of the scheme, especially if the inviscid case is investigated, and (iii) it is analytically
isentropic if viscous effects are neglected, allowing to confirm the EC property of the scheme.
\begin{figure}
    \centering
    \input{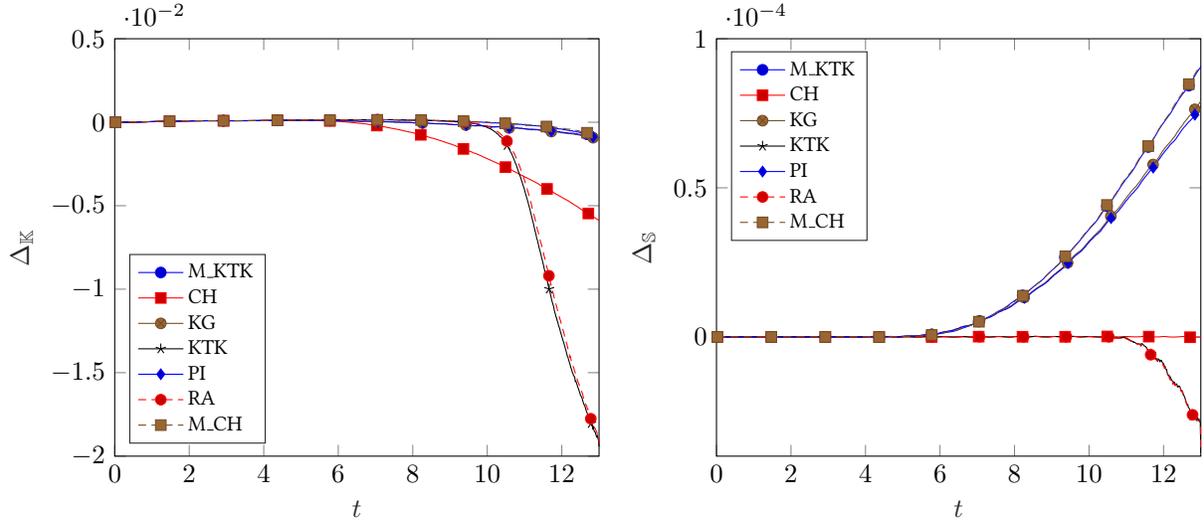}
    \caption{Results of the inviscid TGV computations with $N=3$ on a moving mesh. Shown is the temporal evolution of both the error
	in kinetic energy $\Delta_{\mathbb{K}}$ and entropy $\Delta_{\mathbb{S}}$ for all considered flux functions.}
    \label{fig:TGV_EC_KEP_Euler}
\end{figure}
In a first series of numerical experiments, the behavior of the integral error for both the global entropy
\begin{equation}\label{ErrorTotalEntropy3D}
\Delta_{\mathbb{S}}(T)=\mathbb{S}\left(T\right)-\mathbb{S}\left(0\right),\qquad\mathbb{S}\left(\tau\right):=\sum_{k=1}^{K}\left\langle \interpolation{N}{\left(s\left(\tau\right)\right)},\mathcal{J}\left(\tau\right)\right\rangle _{N}\approx\int_{\mathrm{E}}s\left(\tau\right)J\left(\tau\right)\,d\vec{\xi},\qquad\forall\tau\in\left[0,T\right],
\end{equation} 
with $s=-\frac{1}{\gamma-1}\rho\log\left(p\rho^{-\gamma}\right)$, and the global kinetic energy
\begin{equation}\label{ErrorTotalEKin3D}
\Delta_{\mathbb{K}}(T)=\mathbb{K}\left(T\right)-\mathbb{K}\left(0\right),\qquad\mathbb{K}\left(\tau\right):=\sum_{k=1}^{K}\left\langle \interpolation{N}{\left(\Bbbk\left(\tau\right)\right)},\mathcal{J}\left(\tau\right)\right\rangle _{N}\approx\int_{\mathrm{E}}\Bbbk\left(\tau\right)J\left(\tau\right)\,d\vec{\xi},\qquad\forall\tau\in\left[0,T\right],
\end{equation} 
with $\Bbbk=\frac{1}{2}\rho\left|\vec{u}\right|^{2}$, is investigated for the inviscid TGV (corresponding to Reynolds number $\textrm{Re} \rightarrow \infty$). It is important to note that the numerical integration is again collocated at the corresponding LGL nodes, as only for this discrete integration the conservation properties hold. Since we are at first concerned with the inviscid TGV, no molecular dissipation is present. As
also no additional dissipation terms are added to the surface fluxes, all dissipative effects must be attributed to the numerical
scheme.  This allows us to judge the conservation and preservation properties of the method. The temporal evolution of both mentioned quantities for
all considered flux formulations is plotted in Figure~\ref{fig:TGV_EC_KEP_Euler}.

\begin{figure}
    \centering
    \input{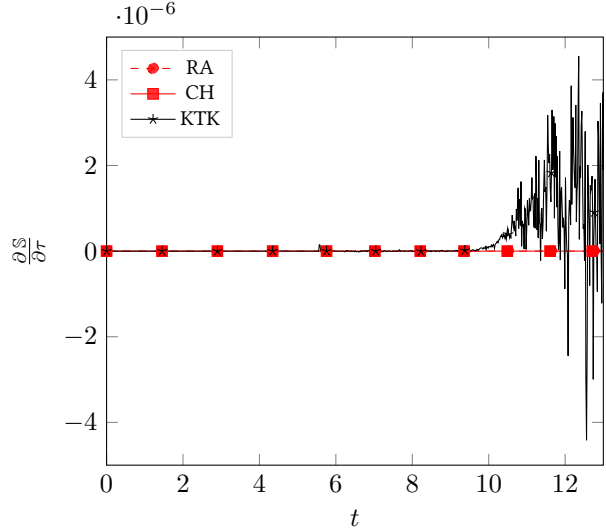}
	\caption{Temporal derivative of the semi-discrete integral entropy \eqref{eqn:DisSDEvolution1} for the inviscid TGV test case, using $N=3$ and the
	standing wave mesh motion.}
    \label{fig:dSdt}
\end{figure}

First we consider the evolution of the kinetic energy. It should be re-iterated that the kinetic energy is not a conserved quantity for compressible flows, but as the Mach number is very small, only minor deviations are expected. The fluxes from Pirozzoli and Kennedy \& Gruber nearly keep the kinetic energy constant over time. Only very late in the simulation, one can observe a small decrease. The flux functions of Ranocha and Kuya, Totani \& Kaway keep the kinetic energy constant slightly past time $t\thickapprox 9$, which is the critical point for stability in the TGV test, but past this point the fluxes show a marked decrease in kinetic energy. The modified flux \eqref{HF_Flux} only differs in the discrete momentum equations ($\mathbf{G}{}_{\iota}^{\upsilon+1,\#}$, $\iota, \upsilon=1,2,3$) from the flux of Kuya et al. This flux shows the same behavior as the fluxes from Pirozzoli and Kennedy \& Gruber. Here the influence of the pressure contribution in the discrete internal energy equations \eqref{LocalDisInternalEnergyPI}, \eqref{LocalDisInternalEnergyKG}, \eqref{LocalDisInternalEnergyKTK} and \eqref{LocalDisInternalEnergyHF} on the preservation of the discrete energy ratio can be observed numerically. The flux originating from Chandrashekar is designed to be KEP and EC in the sense that the discrete criterion from Jameson \cite{Jameson2008} and Tadmor \cite{Tadmor2003} are satisfied. We observe a clear decay in kinetic energy, when this flux is used in the derivative operator \eqref{SpatialDerivativeProjectionOperator2}. This observation is consistent with the experience on static grids made by Gassner et al.~\cite{Gassner2016}, who identified the specific split form of the the pressure term in the Euler equations as the cause for that behavior. Modifying the flux according to \eqref{CH_AVG_Flux} leads to an improved kinetic energy persevering behaviour, but at the cost of loosing conservation in entropy. Our numerical experiments thus support the notion that the discretization of the pressure is the most significant contribution to the balance between internal and kinetic energy. On the other hand, the flux of Chandrashekar shows the expected conservation property in the integral entropy.  Furthermore, the fluxes of Kuya et al. and Ranocha manage to conserve the entropy until very late in the simulation, where both then start to decay. The results here seem to indicate that in practice, those fluxes show an ES and not truly conservative behavior for the
considered test case, at least in the time range that contains the largest under-resolution. It is remarkable that the fluxes of
Kuya et al. and Ranocha show a similar behavior for the entropy, since the flux of Kuya et al. is not designed to be EC in the sense
that the discrete criterion from Tadmor \cite{Tadmor2003} is satisfied. As this behaviour demands further investigation, we now
directly consider the semi-discrete evolution of the entropy, to exclude any influence of the chosen time-integration scheme. To
this end, we mimic the continuous entropy analysis on the semi-discrete level to gain an expression for the semi-discrete evolution of the entropy. Repeating from the equations \eqref{EC:Condition1} and \eqref{EC:Condition2}, the continuous entropy evolution on a moving mesh is
\begin{align}
\begin{split}
\pderivative{\left(Js\right)}{\tau} 
=& \phantom{-} \textbf{w}^{T}\pderivative{\left(J\textbf{u}\right)}{\tau}-\left(\pderivative{J}{\tau}\right)\left(\textbf{w}^{T}\textbf{u}-s\right) \\
=& -\textbf{w}^{T}\vec{\nabla}_{\xi}\cdot\blockvec{\tilde{\textbf{g}}}-\left(\vec{\nabla}_{\xi}\cdot\vec{\tilde{\nu}}\right)\rho \\
=&-\vec{\nabla}_{\xi}\cdot\left(\vec{\tilde{f}}^{s}-\vec{\tilde{\nu}}s\right),
\end{split}
\end{align}
where we used the equations \eqref{GCL} and \eqref{eq:NSE} without the viscous part and the identity $\textbf{w}^{T}\textbf{u}-s=\rho$ in the penultimate step. If the system is investigated with appropriate boundary conditions like periodic boundary conditions, an integration over $\mathrm{E}$ gives    
\begin{equation}\label{eqn:SDEvolution}
\pderivative{}{\tau}\int_{\mathrm{E}}Js\,d\vec{\xi}=-\int_{\mathrm{E}}\vec{\nabla}_{\xi}\cdot\left(\vec{\tilde{f}}^{s}-\vec{\tilde{\nu}}s\right)\,d\vec{\xi}=0.
\end{equation}
Repeating the steps that lead to the expression for the continuous equations in the semi-discrete case, one gains for each element $e_{\alpha}(t)$, $\alpha =1,\dots,K$ 
\begin{align}\label{eqn:DisSDEvolution}
\begin{split}
\pderivative{}{\tau}\left\langle \interpolation{N}{\left(s\right)},\mathcal{J}\right\rangle _{N}
=& \phantom{-} 
\left\langle \pderivative{\left(\mathcal{J}\textbf{U}\right)}{\tau},\mathbf{W}\right\rangle _{N}-\left\langle \pderivative{\mathcal{J}}{\tau},\mathbf{W}^{T}\mathbf{U}-\interpolation{N}{\left(s\right)}\right\rangle _{N} \\
=&-\left\langle \Dprojection{N}\cdot\blockvec{\tilde{\textbf{G}}}{}^{\#},\mathbf{W}\right\rangle _{N}-\int\limits _{\partial\mathrm{E},N}\mathbf{W}^{T}\left(\tilde{\textbf{G}}_{\hat{n}}^{*}-\tilde{\textbf{G}}_{\hat{n}}\right)\,dS \\
&-\left\langle \Dprojection{N}\cdot\vec{\tilde{\nu}}^{\#},\mathbf{W}^{T}\mathbf{U}-\interpolation{N}{\left(s\right)}\right\rangle _{N}-\int\limits _{\partial\mathrm{E},N}\left(\mathbf{W}^{T}\mathbf{U}-\interpolation{N}{\left(s\right)}\right)\left(\tilde{\nu}_{\hat{n}}^{*}-\tilde{\nu}_{\hat{n}}\right)\,dS,
\end{split}
\end{align}
where the equations \eqref{DGSEM:CL} and \eqref{DGSEM:D-GLC} have been used in the last step. Is the split form ALE DGSEM \eqref{MovingMeshDGSEMEuler} investigated with periodic boundary conditions and the volume and surface fluxes are chosen   to be EC, a summation of equation \eqref{eqn:DisSDEvolution} over all elements gives up to machine precision the following discrete analogue of equation \eqref{eqn:SDEvolution}   
\begin{equation}\label{eqn:DisSDEvolution1}
\pderivative{\,\mathbb{S}}{\tau}=\pderivative{}{\tau}\sum_{k=1}^{K}\left\langle \interpolation{N}{\left(s\right)},\mathcal{J}\right\rangle _{N}=0, 
\end{equation}
where $\mathbb{S}$ is given as in \eqref{ErrorTotalEntropy3D}. In Figure~\ref{fig:dSdt}, we show this quantity for the three
flux functions of Chandrashekar, Ranocha and Kuya et al. It can be seen that the fluxes CH and RA show the expected entropy conserving behavior, meaning that the integral change in entropy is zero down to the accuracy expected for finite-precision calculations ($<10^{-14}$). On the other hand, the KTK flux clearly shows that it was not constructed as an entropy conserving flux. It must thus be concluded that the observed decay in entropy for the RA flux is a consequence of the fully discrete system, a
behavior that definitely merits further research in the future. The observed behavior of the EC fluxes directly confirm the claim
that the split form ALE DGSEM is conservative in entropy in the semi-discrete sense for two-point fluxes that have been designed 
following the discrete entropy criteria \eqref{EntropyConditionDGECFlux} for the volume fluxes.

The KEP fluxes of Pirozzoli and Kennedy \& Gruber on the other hand introduce an actual increase in entropy, which could hint at a potential instability in certain situations. The same observation can be made for the modification of the flux of Kuya et al. \eqref{HF_Flux} as well as the modification of the flux of Chandrashekar \eqref{CH_AVG_Flux}. The observed behavior of the various flux functions are all in excellent agreement with the experiments conducted by Gassner et
al.~\cite{Gassner2016} with a similar setup, but on static grids. It can thus be concluded that the split form ALE DGSEM proposed in
this manuscript does indeed extend the conservation properties of the static method to curvilinear moving grids.

\begin{figure}
    \centering
    \input{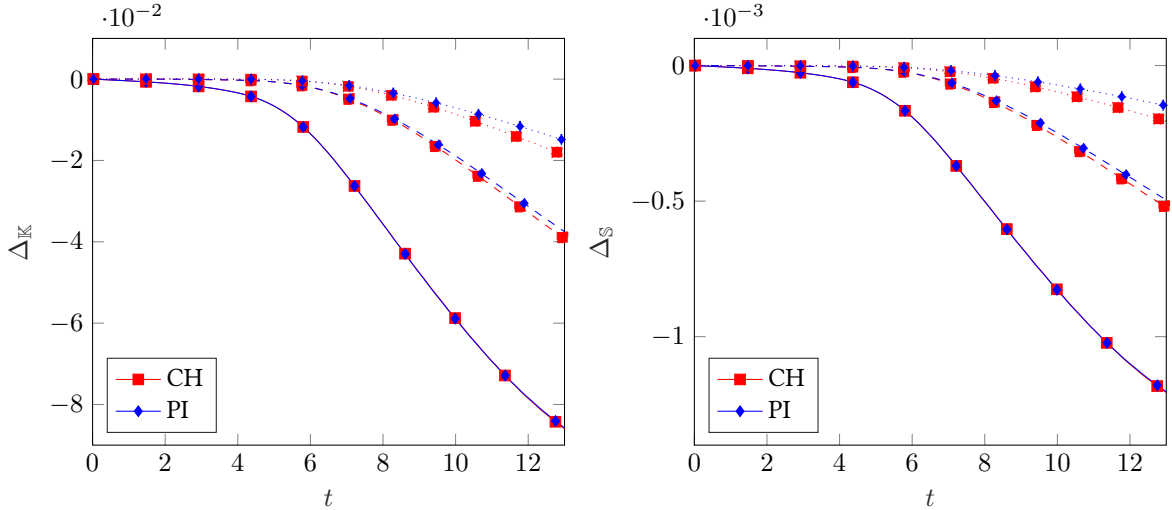}
	\caption{Results of the TGV computations at $\textrm{Re}=1600$ (solid lines), $\textrm{Re}=16,000$ (dashed lines) and $\textrm{Re}=50,000$ (dotted
	lines) with $N=3$ on a moving mesh. Shown is the temporal evolution of both the error in kinetic energy
	$\Delta_{\mathbb{K}}$ and entropy $\Delta_{\mathbb{S}}$ for the flux of Pirozzoli (see Appendix \ref{PI}) and the flux of Chandrashekar (see Appendix \ref{CH}).}
    \label{fig:TGV_EC_KEP_NS}
\end{figure}
So far, only the Euler equations have been considered. Since we are also concerned with the behavior including viscous
contributions, the TGV simulations are repeated for the full NSE. The constant dynamic viscosity $\mu$ (see Section \ref{NSE}) of the fluid is set such that three different Reynolds number of 
\begin{equation}
\textrm{Re}=\frac{\rho \left|\vec{u}_{0}\right|L}{\mu}=\begin{cases}
1600\\
16,000\\
50,000
\end{cases}
\end{equation}
are achieved, where $\left|\vec{u}_{0}\right|$ is the magnitude of the initial velocity in \eqref{TGVInitial} and $L$ a characteristic length chosen as $1$. We note that additional time step restrictions are introduced by the viscous fluxes, but those are not dominant in the considered case, and the
time step is still defined by the CFL restriction \eqref{CFL1}. The molecular dissipation should have a dissipative effect on both the kinetic energy
and the entropy of the system. In Figure~\ref{fig:TGV_EC_KEP_NS}, the temporal evolution of both quantities is shown for the EC flux
of Chandrashekar and the KEP flux following Pirozzoli. It becomes immediately clear that the molecular dissipation dominates any
numerical effect previously observed in the inviscid simulations. For the smallest Reynolds number, associated with the largest
contribution of molecular dissipation, the decay for both quantities is an order of magnitude larger than the difference observed
between flux functions for the Euler case, such that the results for both fluxes are virtually indistinguishable. If the Reynolds
number is increased, some differences start to show. For both the entropy and the kinetic energy, the dissipation observed with the
flux of Pirozzoli is slightly smaller compared with the flux of Chandrashekar. As the molecular dissipation is a function of the
resolved gradients in the flow field, differences in the detail of the flow can explain the observed differences. As expected, the
viscous terms have a dissipative effect in all considered cases.

\begin{figure}
    \centering
    \input{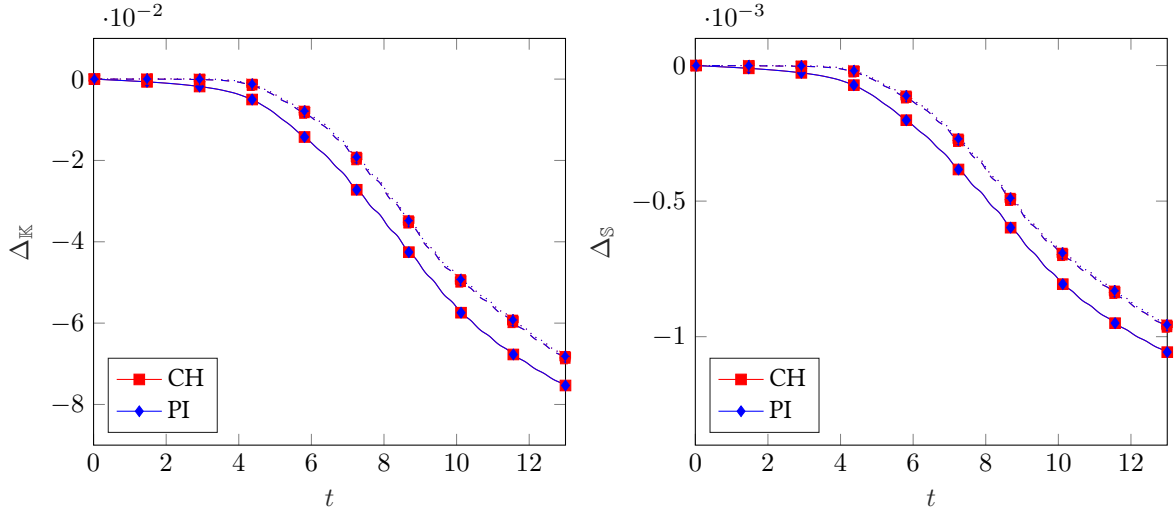}
	\caption{Results of the TGV computations at $\textrm{Re}=1600$ (solid lines), $\textrm{Re}=16,000$ (dashed lines) and $\textrm{Re}=50,000$ (dotted
	lines) with $N=3$ on a moving mesh. Shown is the temporal evolution of both the error in kinetic energy
	$\Delta_{\mathbb{K}}$ and entropy $\Delta_{\mathbb{S}}$ for the flux of Pirozzoli (see Appendix \ref{PI}) and the flux
	of Chandrashekar (see Appendix \ref{CH}), including surface dissipation.}
    \label{fig:TGV_EC_KEP_NS_SurfaceDiss}
\end{figure}
Additional stabilization might still be necessary, especially if the solution features steep gradients. To this end, typically the
surface fluxes are augmented by a matrix dissipation term, see eqn.~\eqref{CartesianSurfaceFlux}. We repeat the simulations of the
viscsous TGV with surface dissipation, to gain insight into the dissipative behaviour of those stabilization terms. In
Figure~\ref{fig:TGV_EC_KEP_NS_SurfaceDiss}, the results for the fluxes of Pirozzoli and Chandrashekar are shown. While the
behaviour for the lowest Reynolds number and thus best resolution of the flow field only slightly changes, the surface dissipation
dominates over the molecular viscous effects for the two considered higher Reynolds numbers. The effect of the surface dissipation
is so dominant that the results are virtually indistinguishable, as the molecular dissipation is already nearly negligible compared
with the advective contribution for the higher Reynolds numbers considered. This clearly shows the stabilizing effect of the
augmented surface fluxes. In certain cases, the dissipation added by the surface fluxes can be used to mimic the effect of a subgrid
scale model for underresolved simulations, a practice which will be applied in the context of large eddy simulation in the next
chapter.

\section{Simulation of transitional flow past a plunging airfoil}\label{Sim_P_Air}

In this section, a complex application of our novel split form ALE DGSEM is presented. As our test case, we choose the low Reynolds number flow around an airfoil that undergoes an unsteady
plunging motion. The different fluid dynamics processes are important for both the understanding of flapping flight and the flows around micro
unmanned vehicles. Such flows have thus attracted considerable attention from engineers and scientists in the past, but efficient
simulation of those cases still remains a challenge. This is due to the complex flow-field that emerges for specific
configurations, which are characterized by large parts of laminar flow, dynamic-stall processes, laminar separation bubbles and
breakdown to turbulence. It is thus necessary to employ unsteady and three-dimensional methods to capture all relevant interactions.
Here, an implicit large eddy simulations (iLES) of the flow around a SD7003 airfoil that is forced to perform a sinusoidal plunging motion is considered. Both experimental~\cite{mcgowan2008computation} and numerical~\cite{visbal2009high} data exist for comparison, such that this test
case can be used to evaluate the performance of the iLES split form ALE DGSEM approach for moving meshes in a complex scenario.

\subsection{Numerical setup}

\begin{figure}
\begin{center}
	\includegraphics[width=0.5\textwidth]{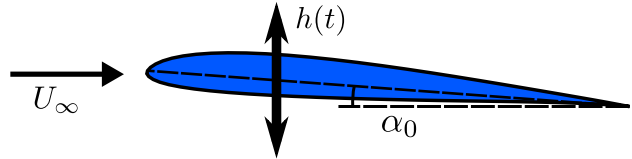}
	\caption{\label{fig:schematicSD7003} Schematic view of the forced plunging motion.}
\end{center}
\end{figure}
The flow around a SD7003 airfoil of chord length $c$ at $\textrm{Re}=40,000$ and $\textrm{Ma}=0.1$ is considered. As can be seen from
sketch~\ref{fig:schematicSD7003}, the airfoil is forced to perform a plunging motion perpendicular to the incoming free stream
velocity $U_{\infty}$. The time-dependent displacement $h(t)$ can be expressed as
\begin{equation}
	h(t)=h_0 \text{sin} \left(  2kF(t)t  \right),
\end{equation}
with the amplitude of the plunging motion $h_0$, the non dimensional frequency $k=\pi f c /U_{\infty}$ and the ramping function
$F(t)$. This function is used to start the simulation from the fully developed flow around a static airfoil and smoothly transition
to the full amplitude of the oscillation, and following Visbal~\cite{visbal2009high} is chosen as
\begin{equation}
	F(t)=1-e^{-at}, \quad a=9.2.
\end{equation}
From several available flow configurations, a non dimensional frequency of $k=3.93$, amplitude of $h_0=0.05c$ and a static angle of
attack $\alpha_0=4 \degree$ were selected, since those lead to a complex, truly three-dimensional flow field. With this settings,
the angle of attack changes due to the plunging motion by up to $21.5 \degree$, which leads to a flow separation at the leading edge
and subsequent breakdown of the leading edge vortex. The details of the flow will be described later.

\begin{figure}
\begin{center}
	\includegraphics[width=0.7\textwidth]{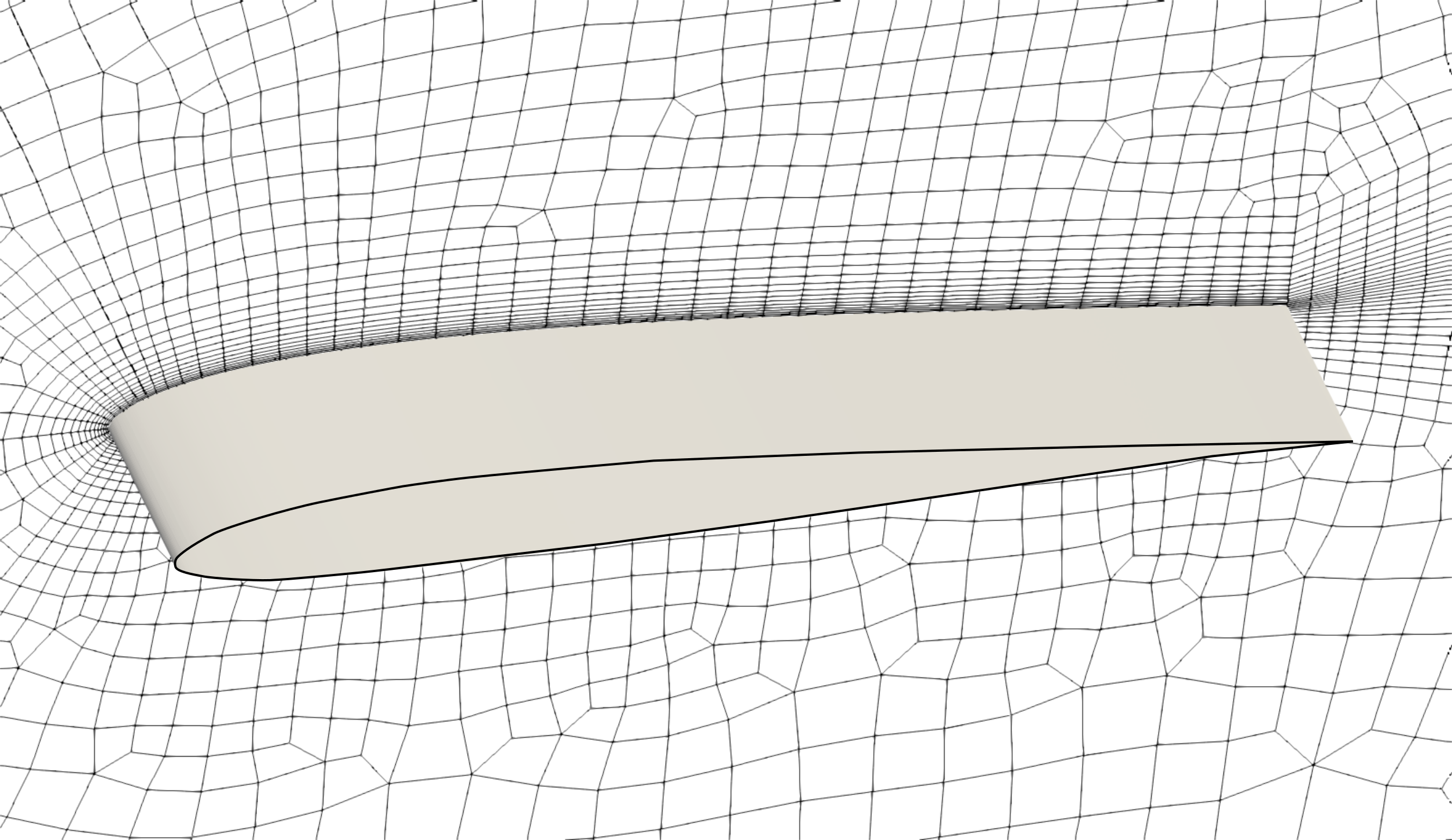}
	\caption{\label{fig:sd7003mesh} Close-up view of the SD7003 airfoil and a slice through the mesh.}
\end{center}
\end{figure}
The spanwise extend of the airfoil is taken as $0.2c$, and periodic boundary conditions are employed in the spanwise direction.
According to~\cite{visbal2009high}, the chosen spanwise extend is enough to neglect an influence of the periodic boundary conditions
on the flow field. On the outer boundary (located $50-100c$ away from the airfoil), the free stream values according to the chosen flow conditions are prescribed, while the airfoil itself is modelled as an adiabatic wall. The construction of entropy stable boundary conditions is an
active field of research, and specific formulations are provided e.g. by Svärd and Özcan~\cite{parsani2015entropy}, Parsani et
al.~\cite{svard2014entropy} or Dalcin et al.~\cite{dalcin2019conservative}. 
We observed no stability issues with the wall boundary conditions for this particular example.

A layer of structured cells is employed around the airfoil (see Figure~\ref{fig:sd7003mesh}), to guarantee
optimal grid quality in the vicinity of the boundary layer. The first cell has a height of about $0.0015c$, and the length of the
cells vary from $0.004c$ at the leading edge to $0.03c$ towards the trailing edge. The structured layer extends $0.1c$ in the
wall-normal direction and consists of $15$ cells, and a grid stretching is used to rapidly increase the height of the cells away
from the wall. On the circumference of the airfoil, $88$ grid cells are used. Outside of the structured layer, the mesh becomes
unstructured and the grid spacing increases rapidly, except for a part reaching around $5c$ into the wake. The rapid increase in
grid spacing helps to efficiently dampen any disturbances before they reach the outer boundaries. All in all, $5849$ cells are
used in a two-dimensional slice of the grid.  In the spanwise direction, the grid is extruded in a structured manner using $10$
cells, leading to a grid spacing of $0.02c$. Thus, in total, the grid consists of $58,490$ cells. The curved geometry is represented
with polynomials of degree $4$ ($N_{geo}=4$).

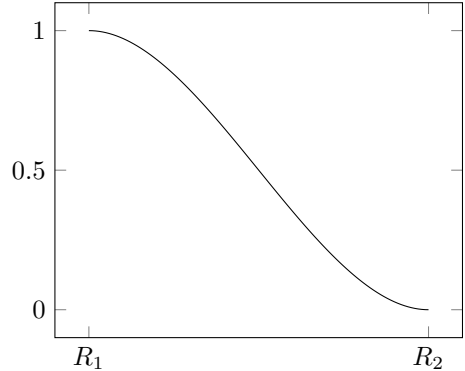
\begin{figure}
\begin{center}

	\begin{tikzpicture}
	\begin{axis}[ xtick={0,1.},
	              xticklabels={$R_1$,$R_2$},
	              width=7cm]
	\addplot [domain=0:1, samples=100]{2.0*x^3-3.0*x^2+1};
	\end{axis}
	\end{tikzpicture}

	\caption{\label{fig:BlendFunc} Simple polynomial blending function for the mesh movement around a rigidly moving body.}
\end{center}
\end{figure}

The movement of the mesh is prescribed as an analytical function. In principle, the whole mesh could be moved rigidly with the
plunging airfoil, but in most scenarios the outer boundaries of the considered domain are best kept stationary. Thus, a blending
approach was implemented, where the mesh close to the moving geometry should move as a rigid body, to keep the desired mesh quality
in the critical areas. We define a radial zone ranging from radius $R_1$ to $R_2$ around the center of the airfoil, where
the mesh closer than $R_1$ should move rigidly with the body and further away than $R_2$ should remain stationary. Then, the
movement of the mesh depending on the radial distance $r$ of the currently considered point to the center of the airfoil can be described by 
\begin{equation}
	\vec{h}^*(r,t) = \begin{cases}
		\vec{h}(t)            & r \le R_1 \\
		\vec{h}(t) \cdot p(r) & R_1 < r < R_2 \\
		0                     & r \ge R_2
	\end{cases},
\end{equation}
where $p(r)$ is a polynomial of third order. This polynomial is designed to fulfill the two conditions $p(R_1)=1$ and $p(R_2)=0$ to
ensure a continuous transition between the inner and the outer zone. Additionally, $p'(R_1)=p'(R_2)=0$ is also required to provide
for a smooth transition between the zones. A plot of this simple blending function can be found in Figure~\ref{fig:BlendFunc}. In theory,
higher order polynomials could be used to set higher derivatives at the interval boundaries equal to zero and create an even
smoother transition, but we observed no clear beneficial effect of that.

A simulation using $N=7$ is performed, leading to around $30$ million degrees of freedom (DOF) per solution variable. We use the kinetic
energy preserving split formulation of Pirozzoli (see Section \ref{PI}). As has been shown by Flad and Gassner~\cite{flad2017use},
this approach leads to an efficient scheme for iLES on a static mesh, and we transfer this approach to problems with moving
domains. For the numerical surface fluxes, we add Roe-type dissipation terms \eqref{EntropyBasedDissipation}, which can be found in \cite[Appendix C.3]{Schnuecke2018}, to the KEP fluxes, leading to a scheme that is KED. To compare the effectiveness and accuracy of the scheme, we also perform a simulation where the EC flux function after Chandrashekar is employed, also with the same Roe type-dissipation, leading to a scheme that is ES. For both simulations, the constant $C_{\text{CFL}}$ in \eqref{CFL1} is set to $C_{\text{CFL}}=0.8$, and the fourth order accurate RK scheme after Kennedy, Carpenter and Louis~\cite{kennedy2000low} (RK4(3)5[2R+]) is used for time integration.

For comparison, the reference simulation by Visbal~\cite{visbal2009high} employed a sixth-order finite difference scheme on a mesh
with $26$ million degrees of freedom. Thus, as the accuracy of the high order schemes is comparable and the number of DOF as well, we expect a rather similar resolution.
\subsection{Results}

The simulation was advanced for $26$ periods of the plunging motion, starting from a fully developed flow around the static
airfoil. After $17$ plunging periods, the flow was considered fully periodic. The following phase-averaged data are created by first averaging in the homogeneous spanwise direction and subsequently over the $9$ periods used to gather
statistics. To denote the different phases, the same notation as in~\cite{visbal2009high} is used: A phase of $\phi=0$ corresponds
to the maximal upward displacement during the plunging motion, $\phi=0.25$ to no displacement and maximal downward velocity,
$\phi=0.5$ to maximal downward displacement and $\phi=0.75$ to no displacement and maximal upward velocity.

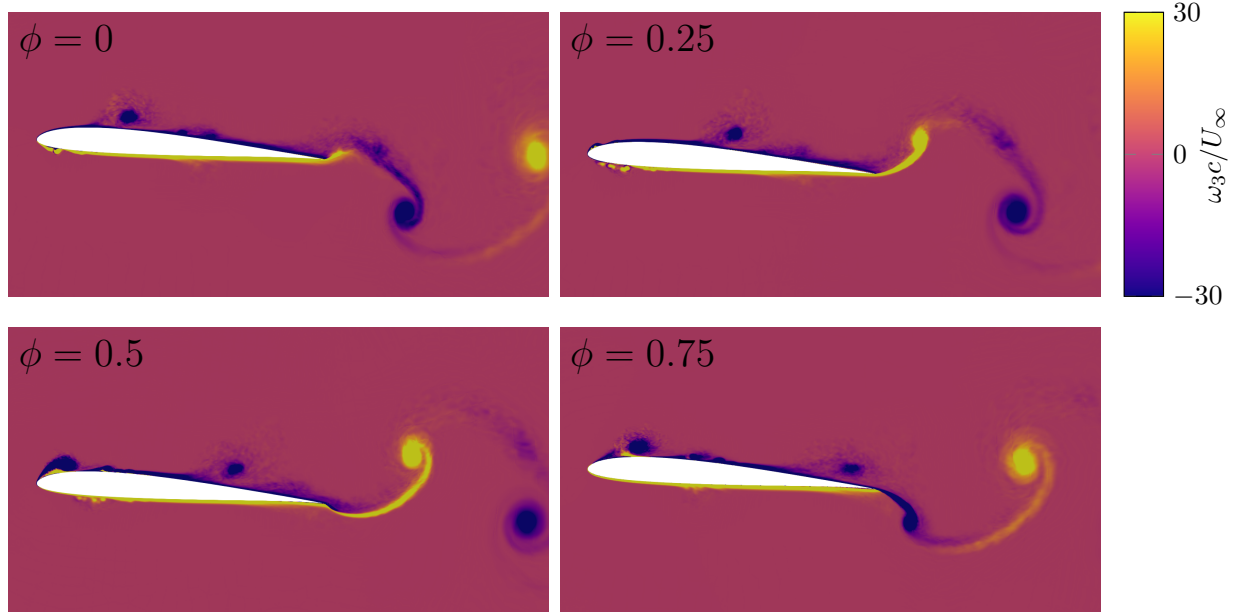
\begin{figure}
\begin{center}
	\def\xsize{7cm}
\def\ysize{4cm}

\begin{tikzpicture}

\begin{scope}
\begin{axis}[
    hide axis,
    axis equal image,
    enlargelimits=false,
    axis on top=true,
    ymin=-0.5,
    ymax=0.5,
    xmin=-0.1,
    xmax=1.8
    ]
    \addplot graphics[xmin=-0.138,ymin=-0.592,xmax=1.970,ymax=0.627] {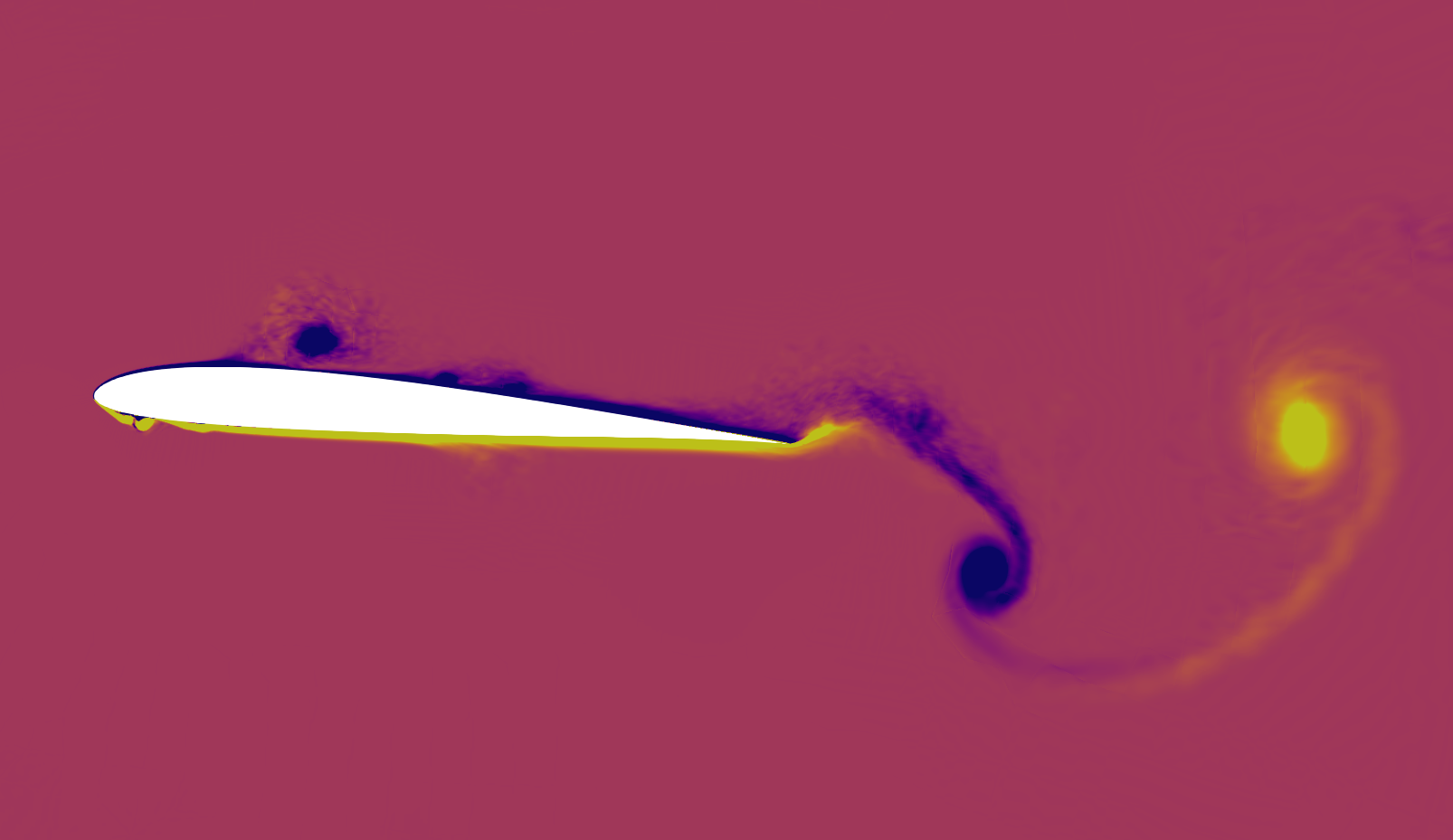};
    \node at (axis cs:-0.1,0.4 )[anchor=west]  {\Large $\phi=0$};
\end{axis}
\end{scope}

\begin{scope}[xshift=\xsize]
\begin{axis}[
    hide axis,
    axis equal image,
    enlargelimits=false,
    axis on top=true,
    ymin=-0.5,
    ymax=0.5,
    xmin=-0.1,
    xmax=1.8,
    colorbar,
    colormap name=plasma,
    point meta min=-30,
    point meta max=30,
    colorbar style={
	    ylabel=$\omega_3 c/U_{\infty}$,
	    y label style={yshift=1.5em},
	ytick={-30,0,30}
    }]
    \addplot graphics[xmin=-0.138,ymin=-0.592,xmax=1.970,ymax=0.627] {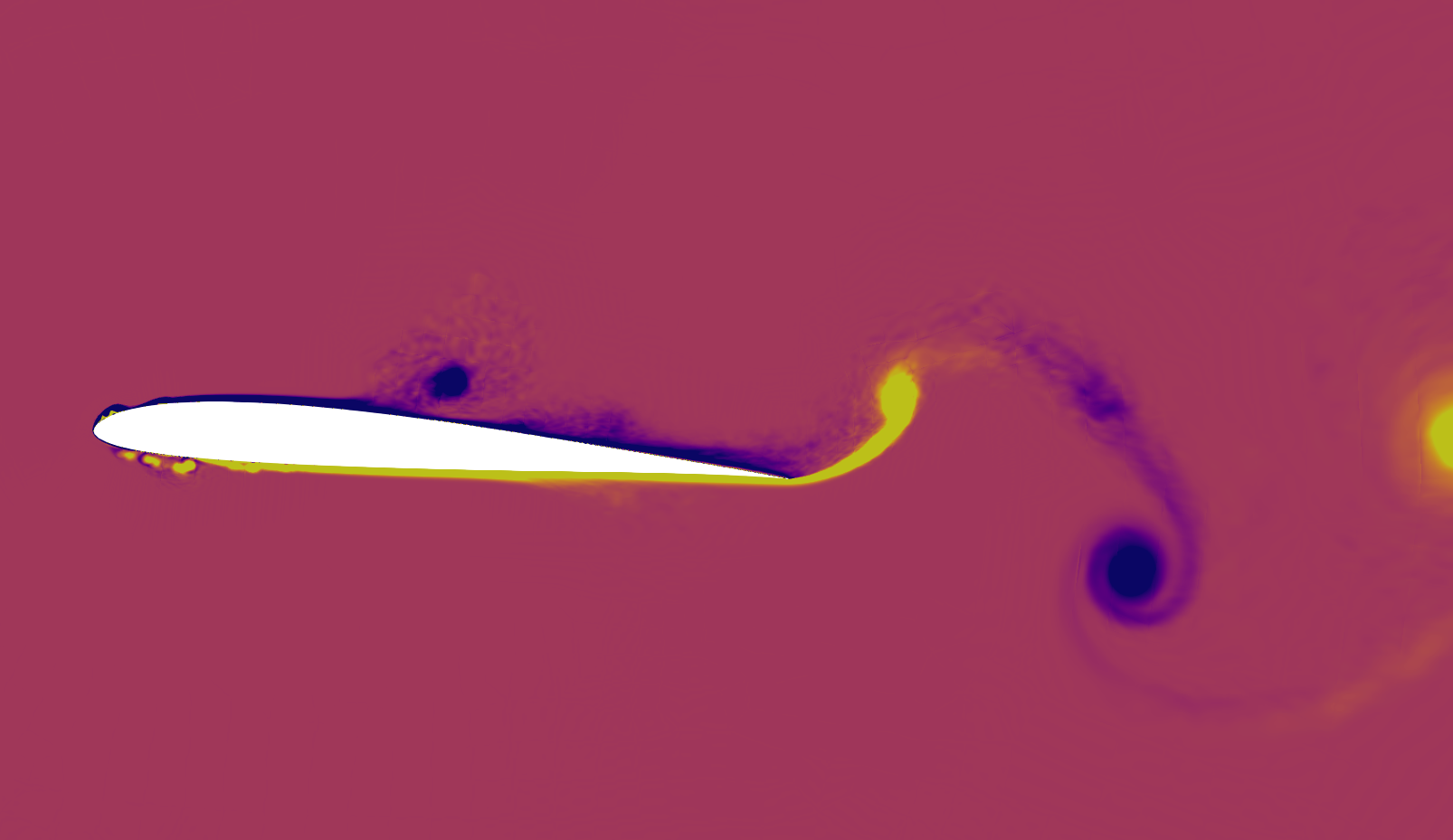};
    \node at (axis cs:-0.1,0.4 )[anchor=west]  {\Large $\phi=0.25$};
\end{axis}
\end{scope}

\begin{scope}[yshift=-\ysize]
\begin{axis}[
    hide axis,
    axis equal image,
    enlargelimits=false,
    axis on top=true,
    ymin=-0.5,
    ymax=0.5,
    xmin=-0.1,
    xmax=1.8
    ]
    \addplot graphics[xmin=-0.138,ymin=-0.592,xmax=1.970,ymax=0.627] {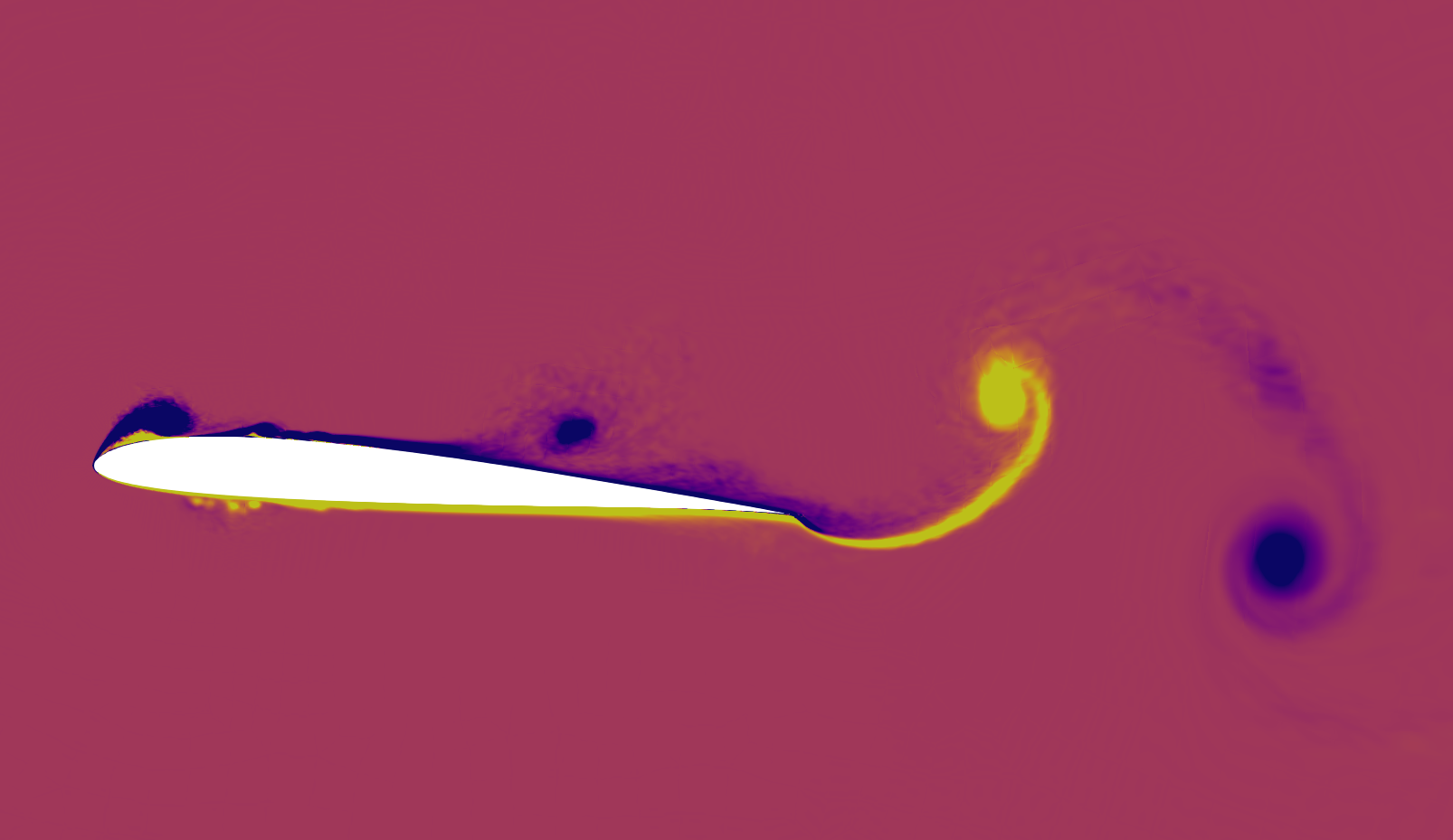};
    \node at (axis cs:-0.1,0.4 )[anchor=west]  {\Large $\phi=0.5$};
\end{axis}
\end{scope}

\begin{scope}[xshift=\xsize,yshift=-\ysize]
\begin{axis}[
    hide axis,
    axis equal image,
    enlargelimits=false,
    axis on top=true,
    ymin=-0.5,
    ymax=0.5,
    xmin=-0.1,
    xmax=1.8
    ]
    \addplot graphics[xmin=-0.138,ymin=-0.592,xmax=1.970,ymax=0.627] {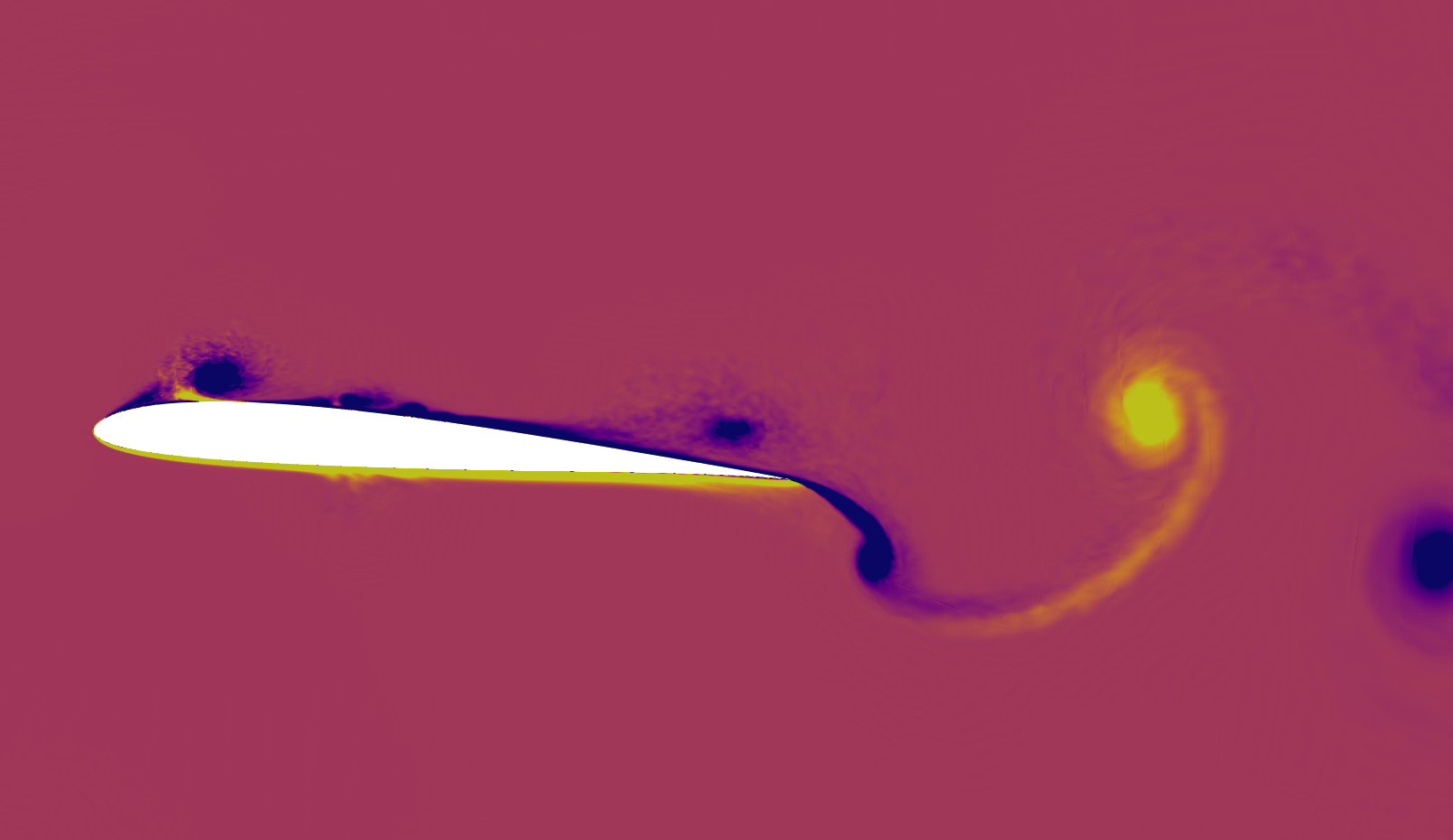};
    \node at (axis cs:-0.1,0.4 )[anchor=west]  {\Large $\phi=0.75$};
\end{axis}
\end{scope}

\end{tikzpicture}%

	\caption{\label{fig:vortZcomp}Phase-averaged spanwise vorticity during different phases of the plunging motion.}
\end{center}
\end{figure}
To gain an overview of the occurring physical processes, Figure~\ref{fig:vortZcomp} shows the phase-averaged spanwise vorticity during
different phases of the plunging motion. If not mentioned otherwise, the results are obtained from the simulation using the KEP
flux of Pirozzoli, as the qualitative comparison of the results from both simulations reveal no major differences. The flow field is
characterized by large areas of laminar flow, and concentrated regions of transitional or turbulent flow. Focusing on the upper
side of the airfoil, the generation of such a region can be traced to the leading edge during the downward part of the motion. At
$\phi=0.5$, the flow at the leading edge can be seen in the process of separation, generating a vortical structure that is
subsequently convected downstream above the airfoil. Two of these regions can be seen at the same time on the airfoil, due to the
relatively slow convective velocity compared to the plunging frequency. The flow separates because the downward plunging motion
induces an increase in angle of attack, surpassing the maximal allowable angle for attached flow, a phenomenon known as dynamic
stall. During the upward part of the motion, the flow around the leading edge re-attaches and purely laminar flow is again obtained.
A similar process can be observed on the lower side during opposite phases of the motion, albeit the generated vortical structures
are much smaller due to the fact that the static angle of attack is positive, and thus the minimal angle of attack is not as
critical as the maximum angle.

\begin{figure}
\begin{center}
	\def\xsize{7cm}
\def\ysize{5.5cm}

\begin{tikzpicture}
\end{tikzpicture}

\begin{tikzpicture}

\begin{scope}
\begin{axis}[
    hide axis,
    axis equal image,
    enlargelimits=false,
    axis on top=true
    ]
    \addplot graphics[xmin=0,ymin=0,xmax=0.225,ymax=0.17044] {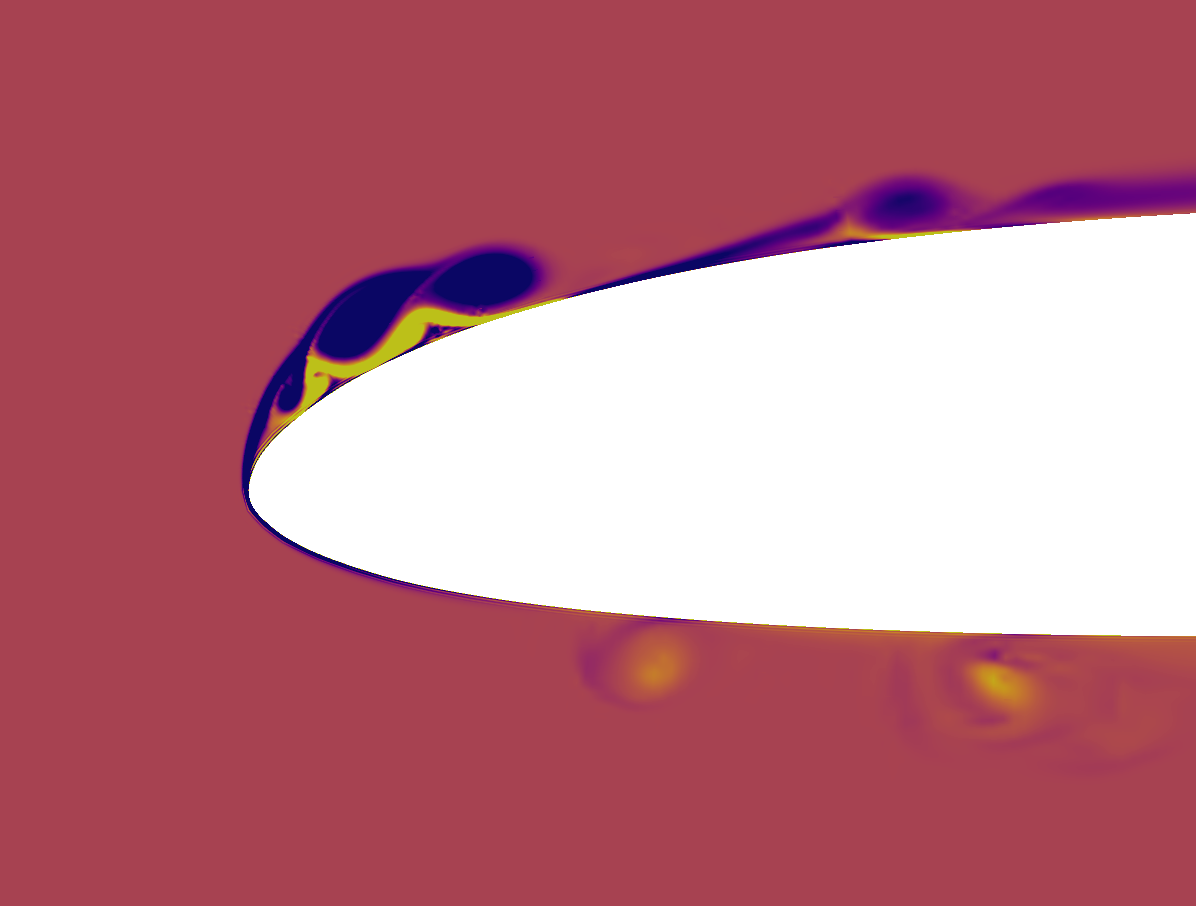};
    \node at (axis cs:0,0.15 )[anchor=west]  {\Large $\phi=0.3$};
\end{axis}
\end{scope}

\begin{scope}[xshift=\xsize]
\begin{axis}[
    hide axis,
    axis equal image,
    enlargelimits=false,
    axis on top=true,
    colorbar,
    colormap name=plasma,
    point meta min=-550,
    point meta max=450,
    colorbar style={
	    ylabel=$\omega_3 c/U_{\infty}$,
	    y label style={yshift=1.5em},
	ytick={-550,0,450}
    }]
    \addplot graphics[xmin=0,ymin=0,xmax=0.225,ymax=0.17044] {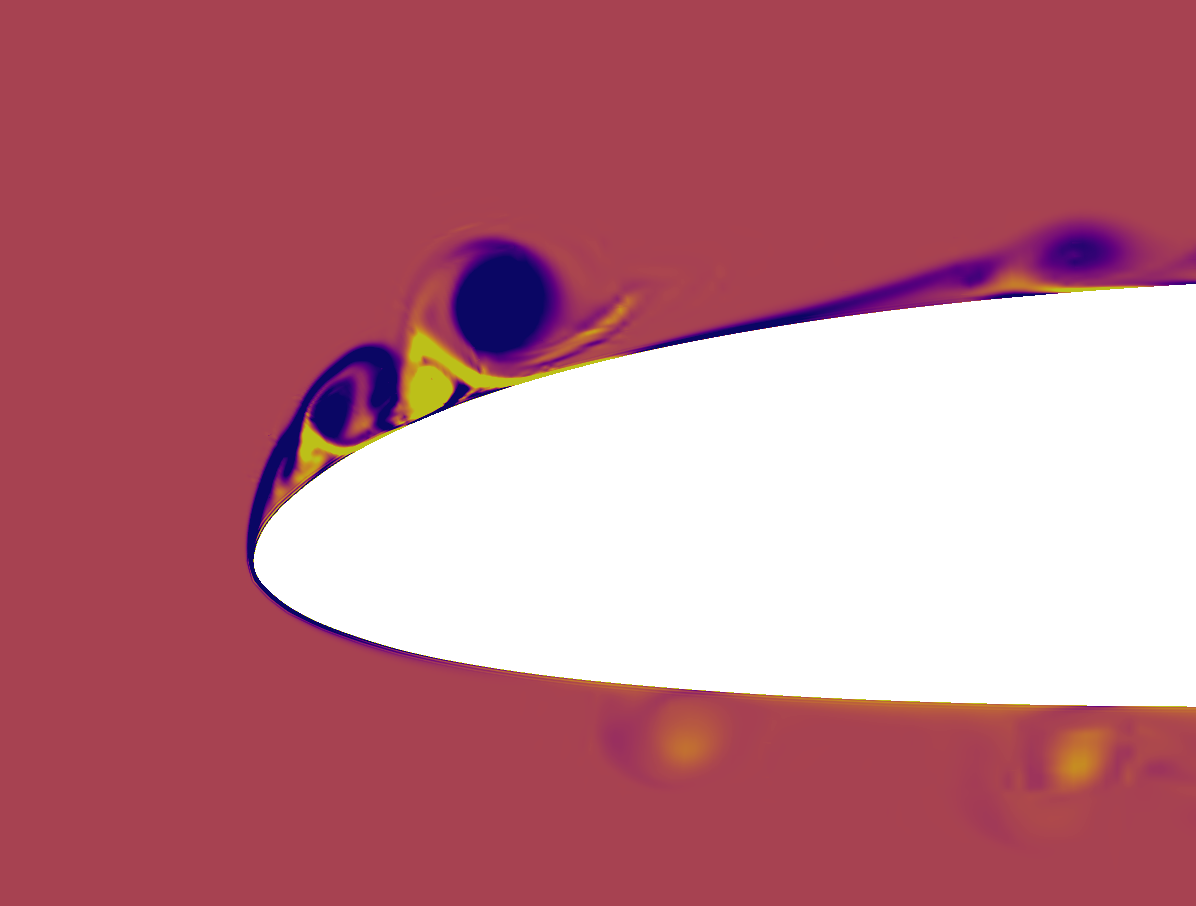};
    \node at (axis cs:0,0.15 )[anchor=west]  {\Large $\phi=0.35$};
\end{axis}
\end{scope}

\begin{scope}[yshift=-\ysize]
\begin{axis}[
    hide axis,
    axis equal image,
    enlargelimits=false,
    axis on top=true
    ]
    \addplot graphics[xmin=0,ymin=0,xmax=0.225,ymax=0.17044] {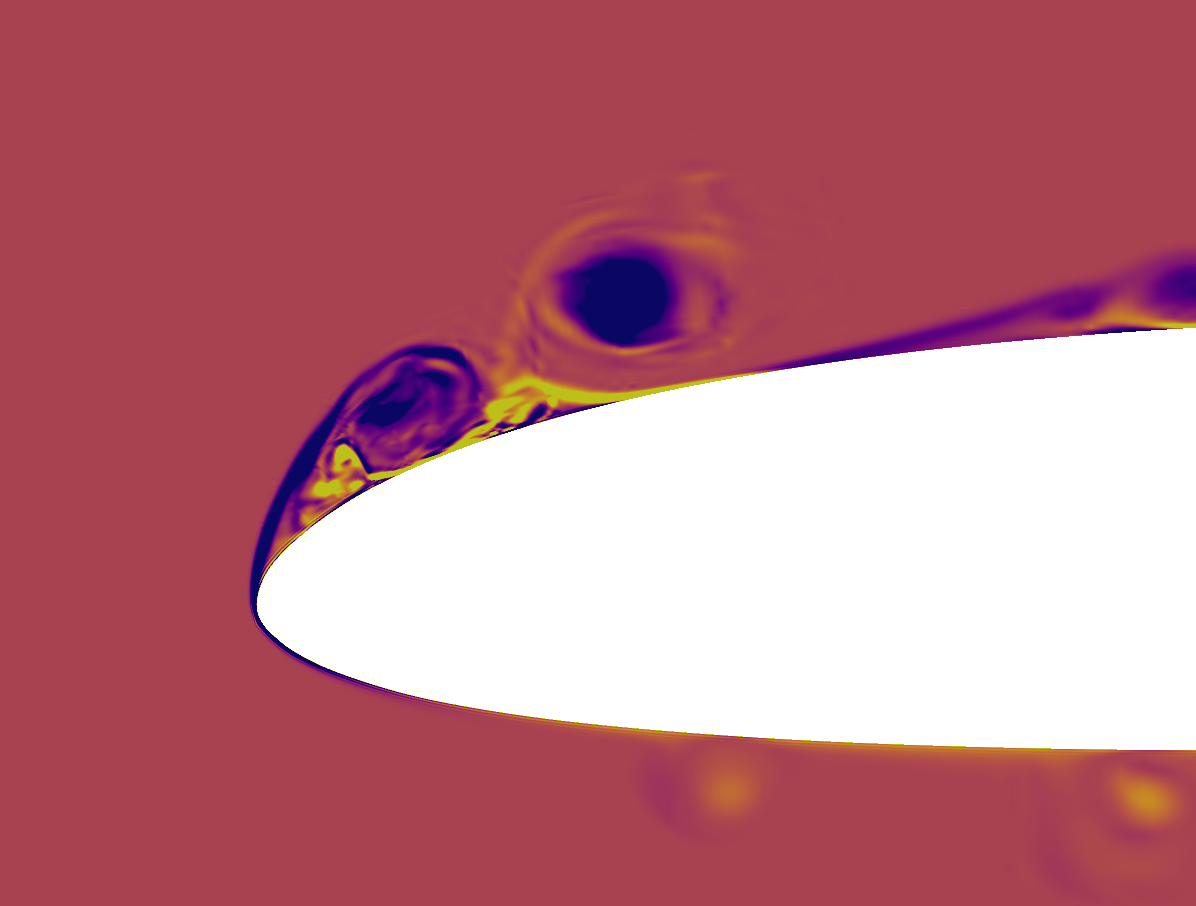};
    \node at (axis cs:0,0.15 )[anchor=west]  {\Large $\phi=0.3875$};
\end{axis}
\end{scope}

\begin{scope}[xshift=\xsize,yshift=-\ysize]
\begin{axis}[
    hide axis,
    axis equal image,
    enlargelimits=false,
    axis on top=true
    ]
    \addplot graphics[xmin=0,ymin=0,xmax=0.225,ymax=0.17044] {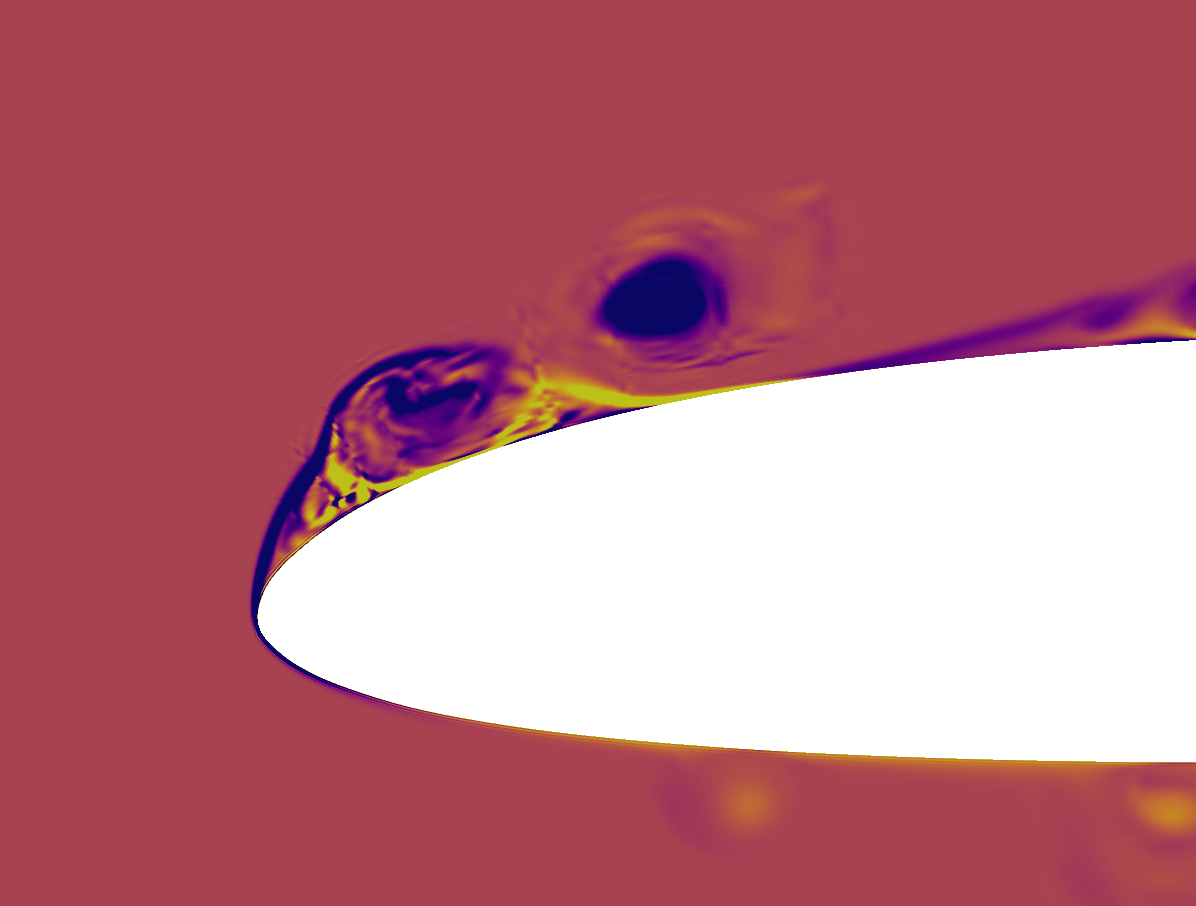};
    \node at (axis cs:0,0.15 )[anchor=west]  {\Large $\phi=0.4$};
\end{axis}
\end{scope}

\end{tikzpicture}%
	\caption{\label{fig:vortZcompClose}Instantaneous spanwise vorticity during the breakdown of the leading-edge vortex.}
\end{center}
\end{figure}
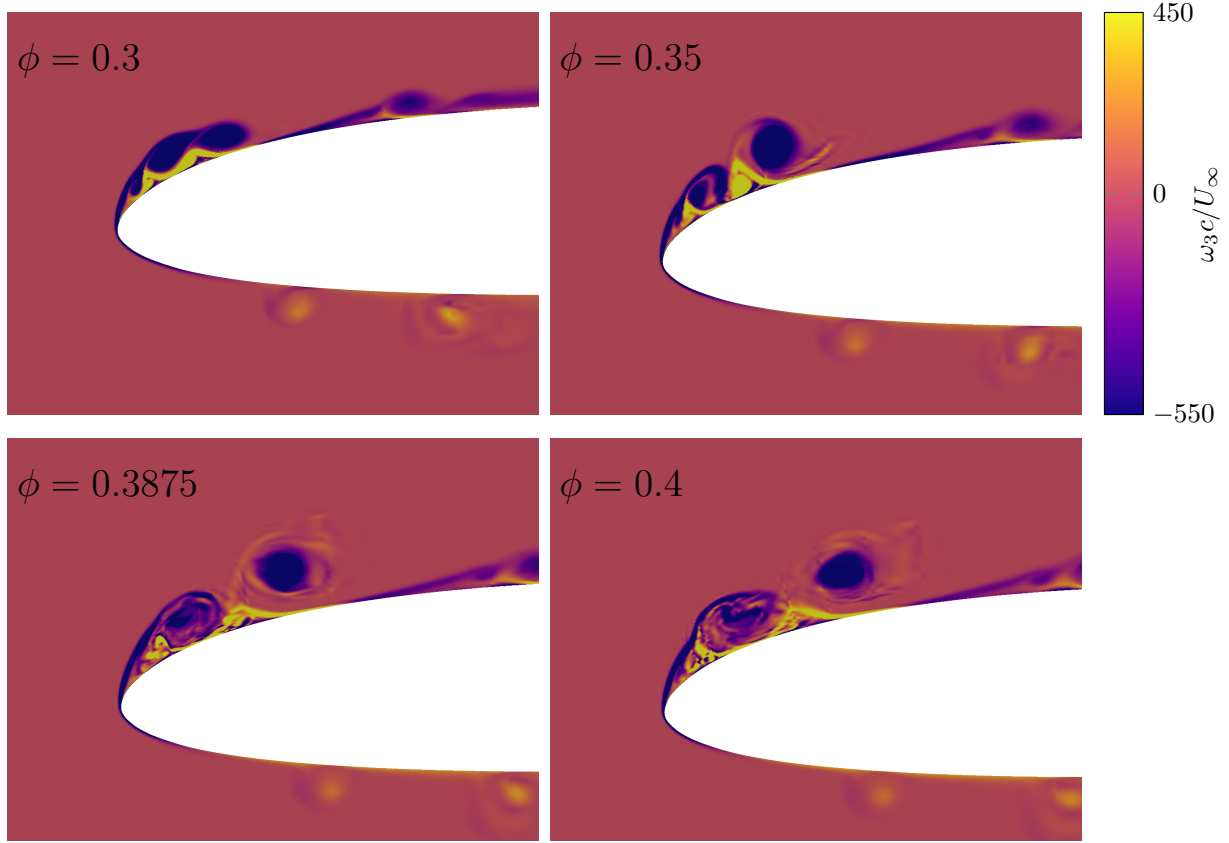
\begin{figure}
\begin{center}
	\includegraphics[width=0.7\textwidth]{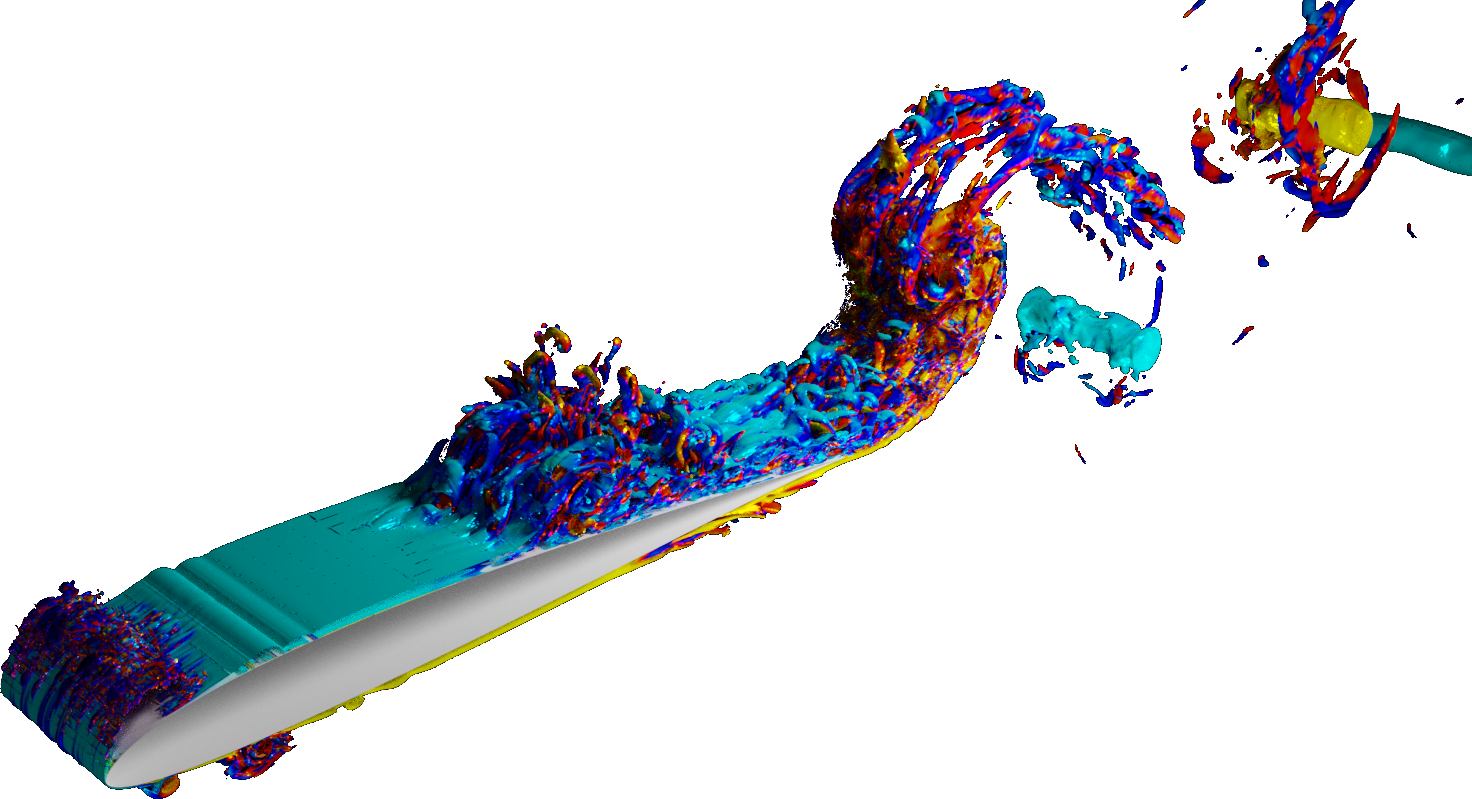}
	\caption{\label{fig:sd7003Render} Render of the three-dimensional flow-field at $\phi=0.3875$. Isocontour of vorticity
	magnitude $|\omega|=40$, colored by spanwise vorticity.}
\end{center}
\end{figure}

Since the process of the breakdown of the dynamic stall vortex happens rather abruptly, a more detailed look can be found
in Figure~\ref{fig:vortZcompClose}. Here, the instantaneous spanwise vorticity is visualized at a two-dimensional slice through the flow
field, such that the dynamics of the process can be inspected. At the beginning of the formation of the vortex ($\phi=0.3$), the
flow is still laminar and two-dimensional. Five distinct vortex-cores can be identified at that stage, three of them rotating
clockwise (negative vorticity) and two embedded, counter-clockwise rotating ones. In the following, these vortices become unstable
and start to break down. This is accompanied by the formation of three-dimensional and transitional structures, which becomes
apparent in the three-dimensional render of the flow field found in Figure~\ref{fig:sd7003Render}. The breakdown of the dynamic stall
vortex has just started, and the flow has become three-dimensional. An area of laminar flow follows further downstream, until the vortex that has been
shed in the previous cycle again constitutes a region of transitional flow. The reported structures are in good agreement with the
results from Visbal~\cite{visbal2009high}, while slight deviations can be observed in the details of the breakdown process. Since
this is the result of an instability, it is expected to be very sensitive to small differences in e.g. discretization. Also, the
exact development of the breakdown depends on the respective two-dimensional slice of the flow one inspects, since the process is
three-dimensional.

In the rendered image~\ref{fig:sd7003Render}, showing isocontours of the magnitude of the vorticity $\vec{\omega}=\nabla \times
\vec{u}$, one can also observe the vortices which are convected downstream in the wake of the airfoil. They alternate between an
upper and a lower vortex, with opposing signs in spanwise vorticity: the upper vortex is rotating counter-clockwise, while the lower
one is rotating clockwise. The vortices dominate the development of the velocity profile in the wake of the airfoil, and those
profiles at $x_1/c=1.5$ (where $x_1$ is the streamwise coordinate) are compared to reference data in Figure~\ref{fig:velCompSD7003} for
the streamwise velocity component $u_1$. This analysis is also used to quantitatively compare the results from the ES (EC flux of Chandrashekar with Roe type-dissipation) and the KED (KEP flux of Pirozzoli with Roe type-dissipation) flux simulation, which are nearly identical. In general, the agreement to both experimental and numerical data is very good. For
$\phi=0$ and $\phi=0.5$, a clear jet-like structure of the profiles can be observed. The fluid has been accelerated in the
streamwise direction, which indicates that the airfoil has experienced a force in the opposing direction, thus creating thrust. This
can be quantified by the mean drag coefficient
\begin{equation}
	C_d=\frac{F_d}{S \cdot \frac{1}{2} \rho_{\infty} U_{\infty}^2},
\end{equation}
with the force acting on the airfoil in streamwise direction $F_d$, the free-stream density $\rho_{\infty}$ and the surface of the airfoil
$S$. For the present simulation, $C_d=-0.088$, the negative value confirming that the motion generates net thrust. The value is also
in good agreement to the simulation by Visbal~\cite{visbal2009high}, where $C_d=-0.082$ was obtained. The ability of a sinusoidal
plunging airfoil to generate thrust is well know and often denoted as the \textit{Knoller-Betz} effect, see
e.g.~\cite{jones1998experimental}.
\begin{figure}
\begin{center}
	\includestandalone[width=\linewidth]{figures/velCompSD7003}
	\caption{\label{fig:velCompSD7003}Comparison of phase-averaged profiles of the streamwise velocity component $u_1$ in the near wake of the airfoil
	($x_1/c=1.5$, where $x_1$ is the streamwise coordinate) over the stream-normal coordinate $x_2$.
	Experimental data from McGowan \etal~\cite{mcgowan2008computation}
	(\protect\tikz[x=2cm,y=2cm,baseline=-0.5ex]\protect\draw[thick,dotted,fill=black] (0,0) -- (0.3,0);),
	simulation from Visbal~\cite{visbal2009high} (\protect\tikz[x=2cm,y=2cm,baseline=-0.5ex]{\protect\draw
	[thick,dashed,color=red] (0,0) -- (0.3,0); 
	})
	and current results using KED fluxes (\protect\tikz[x=2cm,y=2cm,baseline=-0.5ex]\protect\draw [thick,solid,color=blue] (0,0)
	-- (0.3,0);) or ES fluxes (\protect\tikz[x=2cm,y=2cm,baseline=-0.5ex]\protect\draw [thick,solid,color=green] (0,0) -- (0.3,0);).}
\end{center}
\end{figure}

The temporal evolution of the drag coefficient and the lift coefficient
\begin{equation}
	C_l=\frac{F_l}{S \cdot \frac{1}{2} \rho_{\infty} U_{\infty}^2},
\end{equation}
with the force acting on the airfoil perpendicular to the incoming velocity $F_l$, can be seen in Figure~\ref{fig:velForcesSD7003}. Both
quantities show the expected periodic behavior and are in good agreement with the reference simulation. 
\begin{figure}
\begin{center}
	\includestandalone[width=\linewidth]{figures/sd7003compForces}
	\caption{\label{fig:velForcesSD7003}Comparison of the temporal evolution of drag- and lift-coefficients.
	Simulation from Visbal~\cite{visbal2009high} (\protect\tikz[x=2cm,y=2cm,baseline=-0.5ex]{\protect\draw
	[thick,dashed,color=red] (0,0) -- (0.3,0); \protect\draw [solid,color=red,fill=red] (0.15,0) +(-0.03,-0.03)rectangle +(0.03,0.03);} )
	and current results with KED fluxes (\protect\tikz[x=2cm,y=2cm,baseline=-0.5ex]\protect\draw [thick,solid,color=blue] (0,0) --
	(0.3,0);).}
\end{center}
\end{figure}

\begin{table}[]
	\centering
\begin{tabular}{l|cccc}
\toprule[1.5pt]
                & standard & consistent integration & ES fluxes & KED fluxes \\
 \midrule
	PID $[10^{-6} s]$ & $2.15$      & $8.03$                 & $3.27$      & $2.84$ \\
\bottomrule[1.5pt]
\end{tabular}
	\caption{Achieved performance index (PID) for the SD7003 test case with four different moving mesh DG
	methods.\label{tbl:PIDSD7003}}
\end{table}
While the results for both chosen flux functions are remarkably similar, a difference in computational efficiency exists, as the
entropy conservative flux functions require the computation of logarithmic means \eqref{LogMean}, which are computationally expensive compared with
arithmetic means. It is thus worthwhile to compare the  performance of the different schemes with each other. To this end, we
measured the performance index (PID), defined as
\begin{equation}
   \text{PID} = \frac{\text{wall clock time} \cdot \#\text{cores}}{\#\text{DOF} \cdot \#\text{time steps}}.
\end{equation}
It measures the time it takes to advance a single DOF by one time step with the RK scheme. Table~\ref{tbl:PIDSD7003}
compares the two schemes with each other, and adds some reference values. The standard ALE DGSEM \eqref{StandardDGSEM} is the cheapest of the
considered methods, as it does not require the computation of two-point flux functions. It must be noted however, that a simulation
using the standard approach becomes unstable after only a few time steps, since it is plagued by aliasing in the considered,
under-resolved case. A common method to alleviate the problem is to use consistent integration or overintegration~\cite{gassner2013accuracy, kirby2003aliasing, kopriva2018stability,mengaldo2015dealiasing}, where the
numerical integration is evaluated using a quadrature rule of a higher polynomial degree instead of collocating integration and
interpolation. Choosing this approach with an integration rule using $\frac{3}{2}N$, as it is recommended for nearly incompressible
flows, leads to a performance decrease by nearly a factor of $4$ compared to the collocation method. In contrast, the ALE DGSEM using
KEP fluxes is only $\approx 30 \%$ more expensive than the standard DGSEM, while allowing a stable and accurate simulation. As
already mentioned, the ES fluxes are slightly more expensive and lead to a $\approx 50 \%$ increased compute time compared with
standard DG.

All in all, the results from the iLES split form ALE DGSEM approach proposed herein compare well with both the available experimental data and the
reference simulation. This is the case for both a qualitative comparison of the general flow structure as well as for the
quantitative analysis of the available data. The method was able to accurately predict the complex flow field emerging from the
chosen setup, including the rapid breakdown of laminar vortices into fine-scale structures and the resulting transitional process.
The efficiency of the method is shown as the computational time is only slightly increased from the unstable, but cheap
standard DGSEM approach. 

\section{Conclusions}\label{sec:conc}
High order accurate DG methods might be affected by aliasing errors due to the
non-linearity of the flux functions when solving under-resolved turbulent 
vortex dominated flows. One possibility to avoid aliasing issues in the 
discretization is the construction of KED or ES high order split form DG methods 
\cite{Gassner2016,Gassner2017,flad2017use, 
winters2018comparative}. 

In this work, KEP and EC high order ALE DGSEM for the Euler and KED and ES 
high order ALE DGSEM for the NSE were analyzed. Here the key element in the 
approximation is the flux form volume integral of Fisher and Carpenter 
\cite{Fisher2013} or split form DG framework of Gassner et al.~\cite{Gassner2016, Gassner2017}. 
As in \cite{Schnuecke2018} the modification \eqref{EntropyConditionDGECFlux} of the 
discrete entropy criterion from Tadmor \cite{Tadmor2003} was used to construct EC 
volume flux functions in the moving mesh context. In order to construct a provably 
KEP ALE DGSEM for the Euler equations or a KED ALE DGSEM for the NSE, it is required 
that the volume flux in the derivative operator 
\eqref{SpatialDerivativeProjectionOperator2} satisfies the following more 
restrictive version of the discrete KEP condition from Jameson 
\cite{Jameson2008}  
\begin{align}\label{JamesonConclusion}
\textbf{G}{}_{\iota}^{\upsilon+1,\#}=\textbf{G}{}_{\iota}^{1,\#}
\avg{u_{\upsilon}}+\avg{p}\delta_{\iota\upsilon},\qquad \iota,\upsilon=1,2,3,
\end{align}
where the Cartesian surface fluxes $\mathbf{G}_{\iota}^{1,*}$ are consistent 
with the momentum fluxes $\rho\left(u_{\iota}-\nu_{\iota}\right)$, $
\iota=1,2,3$, in the transformed NSE \eqref{eq:NSELawBlockRefE}. Note that the 
given proof for Theorem 3.1 (see Appendix \ref{ProofTheorem3_1}) needs the 
restriction \eqref{JamesonConclusion}, but it is not clear if 
\eqref{JamesonConclusion} is a necessary condition. On the other hand, in 
Section \ref{EC_KEP_Test} numerical experiments with the TGV are presented, 
Figure 4.2 show that the flux of Chandrashekar in Appendix \ref{CH} does 
not provide a KEP ALE DGSEM. This flux satisfies  
\begin{align}\label{JamesonConclusionCH}
\textbf{G}{}_{\iota}^{\upsilon+1,\text{CH}}=\textbf{G}{}_{\iota}^{1,\text{CH}}
\avg{u_{\upsilon}}+\avg{\rho}/\avg{\frac{\rho}{p}}\delta_{\iota\upsilon},\qquad
\iota,\upsilon=1,2,3,
\end{align}
which is on a static mesh an approved KEP condition to construct first order 
KEP FV and FD methods \cite{Chandrashekar2013}. This behavior was also observed by 
Gassner et al.~\cite{Gassner2016} on static grids. Likewise, the numerical 
experiments in Section \ref{EC_KEP_Test} Figure 4.2 show that the fluxes of Kuya, 
Totani, Kawai in Appendix \ref{Kuya} and Ranocha in Appendix \ref{RA} do not 
preserve the kinetic energy and entropy until final time $T=13$, although these 
flux functions satisfy the restriction \eqref{JamesonConclusion}. This behavior can 
be explained by the discrete internal energy equations  
\eqref{LocalDisInternalEnergyKTK}, \eqref{LocalDisInternalEnergyRA}, since 
these equations include a different discrete pressure contribution than the 
internal energy equations  \eqref{LocalDisInternalEnergyPI}, 
\eqref{LocalDisInternalEnergyKG} for the fluxes of Pirozzoli in Appendix 
\ref{PI} and Kennedy \& Gruber in Appendix \ref{KG}, since the fluxes of 
Pirozzoli and Kennedy \& Gruber are KEP according to the numerical experiments in 
Section \ref{EC_KEP_Test} Figure 4.2. Therefore, the results in this work show 
that the treatment of the pressure in the numerical two point volume and 
surface fluxes is the most significant contribution to the balance between 
internal and kinetic energy at least for low Mach number ($\textrm{Ma}
\lesssim 0.3$) turbulent vortex dominated flows, e.g. the TGV test. Overall, it 
should be mentioned that the fluxes of Kuya et al. and Ranocha 
show a very similar behavior, in particular both fluxes are nearly KEP and EC in the 
numerical experiments in Section \ref{EC_KEP_Test} Figure 4.2. This observation is 
very important for questions related to the efficiency of the numerical 
approximation. The flux by Ranocha includes terms with logarithmic averages 
\eqref{LogMean}. These terms are computationally expensive compared with arithmetic 
averages \eqref{VolumeAverages} (see the PID in Table \ref{tbl:PIDSD7003}). Thus, it 
seems to be more convenient to use the flux of Kuya et al. However, a more precise 
investigation showed that the flux of Ranocha gives the expected EC behavior (see in 
Section \ref{EC_KEP_Test} Figure 4.3.). On the other hand, the flux of Kuya et al. 
is not EC.    

Afterward, the KED and ES ALE DGSEM for the NSE were used for a real-world problem 
with moving boundaries. We considered the transitional flow around a plunging SD7003 
airfoil at Reynolds number $\textrm{Re}=40,000$ and Mach number $\textrm{Ma}=0.1$. 
An iLES split form ALE DGSEM approach with the KEP flux of Pirozzoli and the EC flux 
of Chandrashekar as well as Roe type dissipation was used. For a comparison 
experimental measurements and numerical simulation results from literature 
\cite{mcgowan2008computation, visbal2009high} were considered and showed very good 
agreement (see Section \ref{Sim_P_Air} Figures 5.12. and 5.13).

Then the computational performance of the split form ALE DGSEM was compared with the 
overintegrated DG variant. We increased the quadrature nodes by a factor of $
\frac{3}{3}N$ per spatial direction. In comparison with the most efficient, but 
unstable, standard DGSEM \eqref{StandardDGSEM}, the DGSEM using KEP fluxes is only $
\approx 30 \%$ and the DGSEM using ES fluxes is  $\approx 50 \%$  more expensive. 
The overintegration approach increases the computational costs nearly by a factor of 
$4$ compared to the standard DGSEM. Thus, concluding, the novel split form ALE DGSEM 
is an accurate and efficient framework.  

\section*{Acknowledgement} 
Gero Schnücke and Gregor Gassner are supported by the European Research Council (ERC) under the European Union's Eights Framework Program Horizon 2020 with the research project Extreme, ERC grant agreement no. 714487. The authors gratefully acknowledge the support and the computing time on "Hazel Hen" provided by the HLRS through the project "hpcdg".

\appendix
\section*{Appendix} 
\renewcommand\thefigure{\thesection.\arabic{figure}}    
\setcounter{figure}{0} 

\section{Block vector nomenclature from \cite{Gassner2017}}\label{Nomenclature:BlockVector} 
 A block vector is highlighted by the double arrow    
\begin{equation}\label{BlockFlux}
\blockvec{\textbf{f}}:=\begin{bmatrix}\textbf{f}_{1}^{T}, \textbf{f}_{2}^{T}, \textbf{f}_{3}^{T}\end{bmatrix}^{T}, \qquad \textbf{f}_{i} \in \R^{5}, \qquad i=1,2,3.   
\end{equation}
The dot product of two block vectors is given by   
\begin{equation}\label{ProductBlock}
\blockvec{\textbf{f}}\cdot\blockvec{\textbf{g}}:=\sum_{i=1}^{3}\textbf{f}_{i}^{T}\textbf{g}_{i}.
\end{equation}
Furthermore, the dot product of a vector $\vec{v}$ in the three dimensional spatial space and a block vector is defined by    
\begin{equation}\label{ProductVecBlock}
\vec{v}\cdot\blockvec{\textbf{f}}:=\sum_{i=1}^{3}v_{i}\textbf{f}_{i}.
\end{equation}
We note that the dot product \eqref{ProductBlock} is a scalar quantity and the dot product \eqref{ProductVecBlock} is a vector in a $5$ dimensional space, where the number $5$ corresponds to the number of conserved variables in the NSE. The interaction between a vector $\vec{v}$ and the conserved variables is defined as the block vector   
\begin{equation}\label{BlockVector}
\vec{v}\,\textbf{u}:=\begin{bmatrix}v_{1}\textbf{u}\\
v_{2}\textbf{u}\\
v_{3}\textbf{u}
\end{bmatrix}.
\end{equation}
Thus, in particular, the spatial gradient of the conserved variables is defined by     
\begin{equation}\label{BlockGradient}
\vec{\nabla}_{x}\textbf{u}:=\begin{bmatrix}\pderivative{\textbf{u}}{x_{1}}\\
\pderivative{\textbf{u}}{x_{2}}\\
\pderivative{\textbf{u}}{x_{3}}
\end{bmatrix}.
\end{equation}

\section{Transformation of differential operator}\label{TransformationDifferentialOperator} 
The covariant and the contravariant vectors allow to transform differential operators on the time-independent reference element $\mathrm{E}$. In \cite{Kopriva2009}, it has been proven that on the reference element the gradient of a function $f$ is given by 
\begin{equation}\label{Trans:Gradient}
\vec{\nabla}_{x}f=\frac{1}{J}\left(\sum_{i=1}^{3}J\vec{a}^{i}\pderivative{f}{\xi^{i}}\right)
\end{equation}
and the divergence of a vector valued function $\vec{g}$ is given by    
\begin{equation}\label{Trans:Divergence}
\vec{\nabla}_{x}\cdot\vec{g}=\frac{1}{J}\sum_{i=1}^{3}\pderivative{}{\xi^{i}}\left(J\vec{a}^{i}\cdot\vec{g}\right)=\frac{1}{J}\vec{\nabla}_{\xi}\cdot\vec{\tilde{g}}. 
\end{equation}
In \cite{Gassner2017}, the  block matrix \eqref{TransformationBlock} has been introduced to give the following transformation for the gradient and the divergence. Thus, the transformation of the gradient for the state vector $\textbf{u}$ in the NSE becomes     
\begin{equation}\label{BlockGradientReference}
\vec{\nabla}_{x}\textbf{u}=\frac{1}{J}\blockmatx{M}\,\vec{\nabla}_{\xi}\textbf{u} 
\end{equation}
and the the transformation of the divergence for a block vector $\blockvec{\textbf{f}}$ can be written as 
\begin{equation}\label{BlockDivergenceReference}
\vec{\nabla}_{x}\cdot\blockvec{\textbf{f}}=\frac{1}{J}\vec{\nabla}_{\xi}\cdot\blockmatx{M}^{T}\blockvec{\textbf{f}}.
\end{equation}

\section{Numerical flux functions}\label{NumericalFlux} 
In this section is $\iota =1,2,3$, the quantity $\overline{\left|\vec{u}\right|}^{2}$ is given by \eqref{NormAVG} and $\avg{\cdot}^{\text{log}}$ is the logarithmic average  \eqref{LogMean}.
\subsection{The flux from Pirozzoli (PI) \cite{Pirozzoli2010}}\label{PI}
\begin{align}
\begin{split} 
\mathbf{G}{}_{\iota}^{1,\text{PI}\phantom{+1,}}=&\phantom{+} \avg{\rho}\avg{u_{\iota}-\nu_{\iota}}, \\
\mathbf{G}{}_{\iota}^{\upsilon+1,\text{PI}}=& \phantom{+}
\avg{\rho}\avg{u_{\iota}-\nu_{\iota}}\avg{u_{\upsilon}}+\avg{p}\delta_{\iota\upsilon},\qquad\upsilon=1,2,3, \\
\mathbf{G}{}_{\iota}^{5,\text{PI}\phantom{+1,}}=&\phantom{+}  
\avg{\rho}\avg{u_{\iota}-\nu_{\iota}}\left[\avg{e}+\frac{1}{2}\avg{\left|\vec{u}\right|^{2}}\right]+\avg{\rho}\avg{\frac{p}{\rho}}\avg{u_{\iota}}.
\end{split}
\end{align}
\subsection{The flux from Kennedy and Gruber (KG) \cite{Kennedy2008}}\label{KG}  
\begin{align}
\begin{split} 
\mathbf{G}{}_{\iota}^{1,\text{KG}\phantom{+1,}}=&\phantom{+} \avg{\rho}\avg{u_{\iota}-\nu_{\iota}}, \\
\mathbf{G}{}_{\iota}^{\upsilon+1,\text{KG}}=& \phantom{+}
\avg{\rho}\avg{u_{\iota}-\nu_{\iota}}\avg{u_{\upsilon}}+\avg{p}\delta_{\iota\upsilon},\qquad\upsilon=1,2,3, \\
\mathbf{G}{}_{\iota}^{5,\text{KG}\phantom{+1,}}=&\phantom{+}  
\avg{\rho}\avg{u_{\iota}-\nu_{\iota}}\left[\avg{e}+\frac{1}{2}\avg{\left|\vec{u}\right|^{2}}\right]+\avg{p}\avg{u_{\iota}}.
\end{split}
\end{align}
\subsection{The flux from Kuya, Totani and Kawai (KTK) \cite{kuya2018kinetic}}\label{Kuya}  
\begin{align}\label{KTK_Flux}
\begin{split} 
\mathbf{G}{}_{\iota}^{1,\text{KTK}\phantom{+1,}}=&\phantom{+}\avg{\rho}\avg{u_{\iota}-\nu_{\iota}},\\
\mathbf{G}{}_{\iota}^{\upsilon+1,\text{KTK}}=&\phantom{+}\avg{\rho}\avg{u_{\iota}-\nu_{\iota}}\avg{u_{\upsilon}}+\avg{p}\delta_{\iota\upsilon},\qquad\upsilon=1,2,3, \\
\mathbf{G}{}_{\iota}^{5,\text{KTK}\phantom{+1,}}=&\phantom{+}\avg{\rho}\avg{u_{\iota}-\nu_{\iota}}\left[\avg{e}+\frac{1}{2}\overline{\left|\vec{u}\right|}^{2}\right]+2\avg{p}\avg{u_{\iota}}-\avg{pu_{\iota}}.
\end{split}
\end{align} 
\subsection{The flux from Ranocha (RA) \cite{Ranocha2018}}\label{RA} 
\begin{align}
\begin{split} 
\textbf{G}{}_{\iota}^{1,\text{RA}\phantom{+1,}}=& \phantom{+}
\avg{\rho}^{\text{log}}\avg{u_{\iota}-\nu_{\iota}}, \\
\textbf{G}{}_{\iota}^{\upsilon+1,\text{RA}}=& \phantom{+}
\avg{\rho}^{\text{log}}\avg{u_{\iota}-\nu_{\iota}}\avg{u_{\upsilon}}+\avg{p}\delta_{\iota\upsilon}, 
\qquad \upsilon=1,2,3, \\
\textbf{G}{}_{\iota}^{5,\text{RA}\phantom{+1,}}=& \phantom{+}
\avg{\rho}^{\text{log}}\left[\frac{1}{\avg{\frac{1}{e}}^{\text{log}}}+\frac{1}{2}\overline{\left|\vec{u}\right|}^{2}\right]\avg{u_{\iota}-\nu_{\iota}}+2\avg{p}\avg{u_{\iota}}-\avg{pu_{\iota}}.
\end{split}
\end{align}
\subsection{The flux from Chandrashekar (CH) \cite{Chandrashekar2013}}\label{CH} 
\begin{align}
\begin{split} 
\textbf{G}{}_{\iota}^{1,\text{CH}\phantom{+1,}}=& \phantom{+}
\avg{\rho}^{\text{log}}\avg{u_{\iota}-\nu_{\iota}}, \\
\textbf{G}{}_{\iota}^{\upsilon+1,\text{CH}}=& \phantom{+}
\avg{\rho}^{\text{log}}\avg{u_{\iota}-\nu_{\iota}}\avg{u_{\upsilon}}+\frac{\avg{\rho}}{\avg{\frac{\rho}{p}}}\delta_{\iota\upsilon}, 
\qquad \upsilon=1,2,3, \\
\textbf{G}{}_{\iota}^{5,\text{CH}\phantom{+1,}}=& \phantom{+}
\avg{\rho}^{\text{log}}\left[\frac{1}{\avg{\frac{1}{e}}^{\text{log}}}+\frac{1}{2}\overline{\left|\vec{u}\right|}^{2}\right]\avg{u_{\iota}-\nu_{\iota}}+\frac{\avg{\rho}}{\avg{\frac{\rho}{p}}}\avg{u_{\iota}}.
\end{split}
\end{align}
\section{Proofs for the discrete kinetic energy analysis in Section \ref{Sec:KESEuler}}
\subsection{Algebraic tools}
For two states $a$, $b$ with right limits  $a^{+}$, $b^{+}$ and left limits $a^{-}$, $b^{-}$, we have the algebraic relation 
\begin{align}
\jump{ab}=&\avg{a}\jump{b}+\jump{a}\avg{b}, \label{Tool1} 
\end{align}
with respect to the orientated jump average operator \eqref{SurfaceJumpMean}. We note that  the relation is also true for the
volume average operator \eqref{VolumeAverages} and the volume jump
\begin{align}\label{VolumeJump}
\jump{\ast}_{(i,m)jk}:=\left(\ast\right)_{ ijk}+\left(\ast\right)_{ mjk},
\end{align}
if we consider generic nodal values $\left\{ a\right\} _{i=0}^{N}$, $\left\{ b\right\} _{i=0}^{N}$ and $\left\{ c\right\} _{i=0}^{N}$. In addition, the SBP property \eqref{SBP} of the matrix $\matx{Q}$ provides     
\begin{align}\label{SBPRelation1}
\begin{split}
\sum_{i,j=0}^{N}\matx{Q}_{ij}\jump{a}_{\left(i,j\right)}\avg{b}_{\left(i,j\right)}\avg{c}_{\left(i,j\right)}
=&  \phantom{-} \frac{1}{2}\sum_{i,j=0}^{N}\matx{Q}_{ij}a_{i}b_{i}c_{j}+\frac{1}{2}\sum_{i,j=0}^{N}\matx{Q}_{ij}a_{i}b_{j}c_{i} \\
&+\frac{1}{2}\sum_{i,j=0}^{N}\matx{Q}_{ij}a_{i}b_{j}c_{j}-\left[a_{N}b_{N}c_{N}-a_{0}b_{0}c_{0}\right]. 
\end{split}
\end{align}
This identity can be proven in a similar way as the discrete split forms in \cite[Lemma 1]{Gassner2016}. Thus, we skip a proof in this paper.

\subsection{Proof of Theorem 3.1}\label{ProofTheorem3_1} 
We define for $i,j,k=0,\dots,N$ 
\begin{equation}\label{NodalKinEnergyVar}
\mathbf{V}_{ijk}:=\left[\interpolation{N}{\left(\mathbf{v}\right)}\right]_{ijk}=\left[-\frac{1}{2}\left|\vec{u}_{ijk}\right|^{2},\vec{u}_{ijk},0\right]^{T}.
\end{equation}
Then it follows 
\begin{align}\label{NodalKinEnergyVar1}
\begin{split}
\mathbf{V}_{ijk}^{T}\mathbf{U}_{ijk}=\Bbbk_{ijk},\qquad\mathbf{V}_{ijk}^{T}\left(\mathbf{F}_{\iota}\right)_{ijk}=\left(u_{\iota}\right)_{ijk}\left(\Bbbk_{ijk}+p_{ijk}\right), 
\qquad \iota=1,2,3.
\end{split}
\end{align}
We use $\textbf{V}$ (see \eqref{NodalKinEnergyVar}) as test function in equation \eqref{DGSEM:CL} and obtain 
\begin{align}\label{LocalDisKEP1}
\begin{split}
\left\langle 
\pderivative{\left(\mathcal{J}\textbf{U}\right)}{\tau},\mathbf{V}\right\rangle _{N}=&-\left\langle \Dprojection{N}\cdot\blockvec{\tilde{\textbf{G}}}{}^{\#},\mathbf{V}\right\rangle _{N}-\int\limits _{\partial\mathrm{E},N}\mathbf{V}^{T}\left(\tilde{\textbf{G}}_{\hat{n}}^{*}-\tilde{\textbf{G}}_{\hat{n}}\right)\,dS.
\end{split} 
\end{align} 
\textbf{Step 1:} Consider the first discrete volume integral in equation \eqref{LocalDisKEP1}: Suppose the time integration method is exact. Then it is possible to apply the chain rule in time and we obtain for $i,j,k=0,\dots,N$ by \eqref{LocalDisKEP1}
\begin{align}\label{TemporalKinetic1}
\begin{split}
\mathbf{V}_{ijk}^{T}\left(\pderivative{\left(\mathcal{J}\mathbf{U}\right)}{\tau}\right)_{ijk}
=& \left(\pderivative{\mathcal{J}}{\tau}\right)_{ijk}\mathbf{V}_{ijk}^{T}\mathbf{U}_{ijk}+\mathcal{J}_{ijk}\mathbf{V}_{ijk}^{T}\left(\pderivative{\mathbf{U}}{\tau}\right)_{ijk} \\ 
=& \left(\pderivative{\mathcal{J}}{\tau}\right)_{ijk}\Bbbk_{ijk}+\mathcal{J}_{ijk}\left(\pderivative{\Bbbk}{\tau}\right)_{ijk}=\left(\pderivative{\left( \mathcal{J}\Bbbk\right)}{\tau}\right)_{ijk}. 
\end{split}
\end{align}
We multiply the equation \eqref{TemporalKinetic1} by $\omega_{ijk}$ and sum over all quadrature points. This results in the identity 
\begin{equation}
\left\langle \pderivative{\left(\mathcal{J}\mathbf{U}\right)}{\tau},\mathbf{V}\right\rangle _{N}= \pderivative{}{\tau}\left\langle \interpolation{N}{\left(\Bbbk\right)},\mathcal{J}\right\rangle _{N}.
\end{equation}
\noindent
\textbf{Step 2:} Consider the second discrete volume integral in equation \eqref{LocalDisKEP1}: First of all, we obtain by the definition of the derivative projection operator \eqref{SpatialDerivativeProjectionOperator2} 
\begin{align}\label{SpatialKinetic1}
\begin{split}
\left\langle \Dprojection{N}\cdot\blockvec{\tilde{\mathbf{G}}}{}^{\#},\mathbf{V}\right\rangle _{N} 
=& \phantom{{+}}\sum_{j,k=0}^{N}\omega_{j}\omega_{k}\sum_{i,m=0}^{N}2\,\mathcal{Q}_{im}\mathbf{V}_{ijk}^{T}\left(\blockvec{\textbf{G}}{}^{\#}\left(\vec{\nu}_{ijk},\vec{\nu}_{mjk},\textbf{U}_{ijk},\textbf{U}_{mjk}\right)\cdot\avg{\mathbb{J}\vec{a}^{1}}_{\left(i,m\right)jk}\right)\\
& +\sum_{i,k=0}^{N}\omega_{i}\omega_{k}\sum_{j,m=0}^{N}2\,\mathcal{Q}_{jm}\mathbf{V}_{ijk}^{T}
\left(\blockvec{\textbf{G}}{}^{\#}\left(\vec{\nu}_{ijk},\vec{\nu}_{imk},\textbf{U}_{ijk},\textbf{U}_{imk}\right)\cdot\avg{\mathbb{J}\vec{a}^{2}}_{i\left(j,m\right)k}\right)\\
&+\sum_{i,j=0}^{N}\omega_{i}\omega_{j}\sum_{k,m=0}^{N}2\,\mathcal{Q}_{km}\mathbf{V}_{ijk}^{T}\left(\blockvec{\textbf{G}}{}^{\#}\left(\vec{\nu}_{ijk},\vec{\nu}_{ijm},\textbf{U}_{ijk},\textbf{U}_{ijm}\right)\cdot\avg{\mathbb{J}\vec{a}^{3}}_{ij\left(k,m\right)}\right).
\end{split}
\end{align}
In the following, we investigate the first sum on the right hand side in \eqref{SpatialKinetic1}. The SBP property \eqref{SBP} of the matrix $\matx{Q}$ provides $2\matx{Q}_{ij}=\matx{Q}_{ij}-\matx{Q}_{ji}+\matx{B}_{ij}$. Hence, it follows  
\begin{align}\label{SpatialKinetic2}
\begin{split}
& \phantom{+}\sum_{j,k=0}^{N}\omega_{j}\omega_{k}\sum_{i,m=0}^{N}2\,\mathcal{Q}_{im}\mathbf{V}_{ijk}^{T}\left(\blockvec{\textbf{G}}{}^{\#}\left(\vec{\nu}_{ijk},\vec{\nu}_{mjk},\textbf{U}_{ijk},\textbf{U}_{mjk}\right)\cdot\avg{\mathbb{J}\vec{a}^{1}}_{\left(i,m\right)jk}\right)\\
=& \phantom{+} \sum_{j,k=0}^{N}\omega_{j}\omega_{k}\sum_{i,m=0}^{N}\mathcal{Q}_{im}\jump{\mathbf{V}}_{\left(i,m\right)jk}^{T}
\left(\blockvec{\textbf{G}}{}^{\#}\left(\vec{\nu}_{ijk},\vec{\nu}_{mjk},\textbf{U}_{ijk},\textbf{U}_{mjk}\right)\cdot\avg{\mathbb{J}\vec{a}^{1}}_{\left(i,m\right)jk}\right)\\
&+\sum_{j,k=0}^{N}\omega_{j}\omega_{k}\sum_{i,m=0}^{N}\mathcal{B}_{im}\mathbf{V}_{ijk}^{T}
\left(\blockvec{\textbf{G}}{}^{\#}\left(\vec{\nu}_{ijk},\vec{\nu}_{mjk},\textbf{U}_{ijk},\textbf{U}_{mjk}\right)\cdot\avg{\mathbb{J}\vec{a}^{1}}_{\left(i,m\right)jk}\right),
\end{split}
\end{align}
since the flux functions $\mathbf{G}{}_{\iota}^{\#}$, $\iota=1,2,3$, are symmetric. Next, we apply the condition \eqref{Jameson}, \eqref{Tool1} and \eqref{SBPRelation1} to evaluate the first sum on the right hand side in \eqref{SpatialKinetic2}. This results in   
\begin{align}\label{SpatialKinetic3}
\begin{split}
& \phantom{+} \sum_{j,k=0}^{N}\omega_{j}\omega_{k}\sum_{i,m=0}^{N}\mathcal{Q}_{im}\jump{\mathbf{V}}_{\left(i,m\right)jk}^{T}\left(\blockvec{\textbf{G}}{}^{\#}\left(\vec{\nu}_{ijk},\vec{\nu}_{mjk},\textbf{U}_{ijk},\textbf{U}_{mjk}\right)\cdot\avg{\mathbb{J}\vec{a}^{1}}_{\left(i,m\right)jk}\right) \\ 
=& \phantom{+} \sum_{\beta=1}^{3}\sum_{j,k=0}^{N}\omega_{j}\omega_{k}\sum_{i,m=0}^{N}\mathcal{Q}_{im}\left(\jump{\mathbf{V}}_{\left(i,m\right)jk}^{T}\textbf{G}{}_{\beta}^{\#}\left(\vec{\nu}_{ijk},\vec{\nu}_{imk},\textbf{U}_{ijk},\textbf{U}_{imk}\right)\right)\avg{\mathbb{J}\vec{a}_{\beta}^{1}}_{\left(i,m\right)jk}\\
=& -\sum_{\beta=1}^{3}\sum_{j,k=0}^{N}\omega_{j}\omega_{k}\sum_{i,m=0}^{N}\mathcal{Q}_{im}\frac{1}{2}\jump{\left|\vec{u}\right|^{2}}_{\left(i,m\right)jk}\mathbf{G}{}_{\beta}^{1,\#}\left(\vec{\nu}_{ijk},\vec{\nu}_{mjk},\mathbf{U}_{ijk},\mathbf{U}_{mjk}\right)\avg{\mathbb{J}\vec{a}_{\beta}^{1}}_{\left(i,m\right)jk} \\
&+\sum_{\beta=1}^{3}\sum_{j,k=0}^{N}\omega_{j}\omega_{k}\sum_{i,m=0}^{N}\mathcal{Q}_{im}\jump{\vec{u}}_{\left(i,m\right)jk}\cdot\avg{\vec{u}}_{\left(i,m\right)jk}\mathbf{G}{}_{1}^{1,\#}\left(\vec{\nu}_{ijk},\vec{\nu}_{mjk},\mathbf{U}_{ijk},\mathbf{U}_{mjk}\right)\avg{\mathbb{J}\vec{a}_{\beta}^{1}}_{\left(i,m\right)jk} \\
& + \sum_{\beta=1}^{3}\sum_{j,k=0}^{N}\omega_{j}\omega_{k}\sum_{i,m=0}^{N}\mathcal{Q}_{im}\jump{u_{\beta}}_{\left(i,m\right)jk}\avg{p}_{\left(i,m\right)jk}\avg{\mathbb{J}\vec{a}_{\beta}^{1}}_{\left(i,m\right)jk} \\
=& \phantom{+} \frac{1}{2}\sum_{\beta=1}^{3}\sum_{j,k=0}^{N}\omega_{j}\omega_{k}\sum_{i,m=0}^{N}\mathcal{Q}_{im}\left(u_{\beta}\right)_{ijk}p_{ijk}\left(\mathbb{J}\vec{a}_{\beta}^{1}\right)_{mjk} \\
&+\frac{1}{2}\sum_{\beta=1}^{3}\sum_{j,k=0}^{N}\omega_{j}\omega_{k}\sum_{i,m=0}^{N}\mathcal{Q}_{im}\left(u_{\beta}\right)_{ijk}p_{mjk}\left(\mathbb{J}\vec{a}_{\beta}^{1}\right)_{ijk} \\
&+\frac{1}{2}\sum_{\beta=1}^{3}\sum_{j,k=0}^{N}\omega_{j}\omega_{k}\sum_{i,m=0}^{N}\mathcal{Q}_{im}\left(u_{\beta}\right)_{ijk}p_{mjk}\left(\mathbb{J}\vec{a}_{\beta}^{1}\right)_{mjk} \\
&-\sum_{\beta=1}^{3}\sum_{j,k=0}^{N}\omega_{j}\omega_{k}\left(\left(u_{\beta}\right)_{Njk}p_{Njk}\left(\mathbb{J}\vec{a}_{\beta}^{1}\right)_{Njk}-\left(u_{\beta}\right)_{0jk}p_{0jk}\left(\mathbb{J}\vec{a}_{\beta}^{1}\right)_{0jk}\right)\\
=& \phantom{+}\frac{1}{2}\left\langle p\left(\pderivative{\interpolation{N}{\left(\mathbb{J}\vec{a}^{1}\right)}}{\xi^{1}}\right)+\left(\pderivative{\interpolation{N}{\left(p\right)}}{\xi^{1}}\right)\mathbb{J}\vec{a}^{1}+\left(\pderivative{\interpolation{N}{\left(p\mathbb{J}\vec{a}^{1}\right)}}{\xi^{1}}\right),\vec{u}\right\rangle _{N}-\int\limits _{\partial\mathrm{E},N}\hat{n}^{1}\left\{ p\mathbb{J}\vec{a}^{1}\cdot\vec{u}\right\} \,dS.
\end{split}
\end{align}
Since the flux functions $\mathbf{G}{}_{\iota}^{\#}$, $\iota=1,2,3$ are consistent with the contravariant flux vectors
$\mathbf{G}_{\iota}=\interpolation{N}{\left(\mathbf{f}_{\iota}-\nu_{\iota}\mathbf{u}\right)}$, where $\mathbf{f}_{\iota}$ are the
flux functions given by \eqref{AdvectiveFluxes} and $\mathbf{u}$ is the state vector given by \eqref{Statevector}, it follows by the definition of the matrix $\matx{B}$ 
\begin{align}\label{SpatialKinetic4}
\begin{split}
& \phantom{+} \sum_{j,k=0}^{N}\omega_{j}\omega_{k}\sum_{i,m=0}^{N}\mathcal{B}_{im}\mathbf{V}_{ijk}^{T}
\left(\blockvec{\textbf{G}}{}^{\#}\left(\vec{\nu}_{ijk},\vec{\nu}_{mjk},\textbf{U}_{ijk},\textbf{U}_{mjk}\right)\cdot\avg{\mathbb{J}\vec{a}^{1}}_{\left(i,m\right)jk}\right) \\ 
=& \phantom{+} \sum_{j,k=0}^{N}\omega_{j}\omega_{k}\left[\mathbf{V}_{Njk}^{T}\left\{ \left(\mathbb{J}\vec{a}^{1}\right)_{Njk}\cdot\blockvec{\mathbf{G}}_{Njk}\right\} -\mathbf{V}_{0jk}^{T}\left\{ \left(\mathbb{J}\vec{a}^{1}\right)_{0jk}\cdot\blockvec{\mathbf{G}}_{0jk}\right\} \right] \\
=& \int\limits _{\partial\mathrm{E},N}\hat{n}^{1}\left\{ \mathbf{V}^{T}\left[\mathbb{J}\vec{a}^{1}\cdot\blockvec{\mathbf{G}}\right]\right\} \,dS.
\end{split}
\end{align}
Hence, by \eqref{SpatialKinetic2}, \eqref{SpatialKinetic3}, \eqref{SpatialKinetic4}, the first sum on the right hand side in \eqref{SpatialKinetic1} can be written as follows  
\begin{align}\label{SpatialKinetic5}
\begin{split}
& \phantom{+}\sum_{j,k=0}^{N}\omega_{j}\omega_{k}\sum_{i,m=0}^{N}2\,\mathcal{Q}_{im}\mathbf{V}_{ijk}^{T}\left(\blockvec{\textbf{G}}{}^{\#}\left(\vec{\nu}_{ijk},\vec{\nu}_{mjk},\textbf{U}_{ijk},\textbf{U}_{mjk}\right)\cdot\avg{\mathbb{J}\vec{a}^{1}}_{\left(i,m\right)jk}\right) \\
=&  \phantom{+}
\frac{1}{2}\left\langle p\left(\pderivative{\interpolation{N}{\left(\mathbb{J}\vec{a}^{1}\right)}}{\xi^{1}}\right)+\left(\pderivative{\interpolation{N}{\left(p\right)}}{\xi^{1}}\right)\mathbb{J}\vec{a}^{1}+\left(\pderivative{\interpolation{N}{\left(p\mathbb{J}\vec{a}^{1}\right)}}{\xi^{1}}\right),\vec{u}\right\rangle _{N}\\
&+\int\limits _{\partial\mathrm{E},N}\hat{n}^{1}\left\{ \mathbf{V}^{T}\left[\mathbb{J}\vec{a}^{1}\cdot\blockvec{\mathbf{G}}\right]-p\mathbb{J}\vec{a}^{1}\cdot\vec{u}\right\} \,dS. 
\end{split}
\end{align}
The second and third sum on the right hand side in \eqref{SpatialKinetic1} can be analyzed in the same way. This provides the identities  
\begin{align}\label{SpatialKinetic6}
\begin{split}
& \sum_{i,k=0}^{N}\omega_{i}\omega_{k}\sum_{j,m=0}^{N}2\,\mathcal{Q}_{im}\mathbf{V}_{ijk}^{T}\left(\blockvec{\textbf{G}}{}^{\#}\left(\vec{\nu}_{ijk},\vec{\nu}_{imk},\textbf{U}_{ijk},\textbf{U}_{imk}\right)\cdot\avg{\mathbb{J}\vec{a}^{2}}_{i\left(j,m\right)k}\right) \\
=&\phantom{+}  
\frac{1}{2}\left\langle p\left(\pderivative{\interpolation{N}{\left(\mathbb{J}\vec{a}^{2}\right)}}{\xi^{2}}\right)+\left(\pderivative{\interpolation{N}{\left(p\right)}}{\xi^{2}}\right)\mathbb{J}\vec{a}^{2}+\left(\pderivative{\interpolation{N}{\left(p\mathbb{J}\vec{a}^{2}\right)}}{\xi^{2}}\right),\vec{u}\right\rangle _{N}\\
&+\int\limits _{\partial\mathrm{E},N}\hat{n}^{2}\left\{ \mathbf{V}^{T}\left[\mathbb{J}\vec{a}^{2}\cdot\blockvec{\mathbf{G}}\right]-p\mathbb{J}\vec{a}^{2}\cdot\vec{u}\right\} \,dS, 
\end{split}
\end{align}
\begin{align}\label{SpatialKinetic7}
\begin{split}
& \phantom{+}\sum_{i,j=0}^{N}\omega_{i}\omega_{j}\sum_{k,m=0}^{N}2\,\mathcal{Q}_{km}\mathbf{V}_{ijk}^{T}\left(\blockvec{\textbf{G}}{}^{\#}\left(\vec{\nu}_{ijk},\vec{\nu}_{ijm},\textbf{U}_{ijk},\textbf{U}_{ijm}\right)\cdot\avg{\mathbb{J}\vec{a}^{3}}_{ij\left(k,m\right)}\right) \\
=& \phantom{+}
\frac{1}{2}\left\langle p\left(\pderivative{\interpolation{N}{\left(\mathbb{J}\vec{a}^{3}\right)}}{\xi^{3}}\right)+\left(\pderivative{\interpolation{N}{\left(p\right)}}{\xi^{3}}\right)\mathbb{J}\vec{a}^{3}+\left(\pderivative{\interpolation{N}{\left(p\mathbb{J}\vec{a}^{3}\right)}}{\xi^{3}}\right),\vec{u}\right\rangle _{N} \\
&+\int\limits _{\partial\mathrm{E},N}\hat{n}^{3}\left\{ \mathbf{V}^{T}\left[\mathbb{J}\vec{a}^{3}\cdot\blockvec{\mathbf{G}}\right]-p\mathbb{J}\vec{a}^{3}\cdot\vec{u}\right\} \,dS.
\end{split}
\end{align}
Finally, we combine \eqref{SpatialKinetic1}, \eqref{SpatialKinetic5}, \eqref{SpatialKinetic6}, \eqref{SpatialKinetic7}, and obtain 
\begin{align}\label{SpatialKinetic9}
\begin{split}
\left\langle \Dprojection{N}\cdot\blockvec{\tilde{\mathbf{G}}}{}^{\#},\mathbf{V}\right\rangle _{N} 
=& \phantom{+} \frac{1}{2}\sum_{\iota=1}^{3}\left\langle p\left(\pderivative{\interpolation{N}{\left(\mathbb{J}\vec{a}^{\iota}\right)}}{\xi^{\iota}}\right)+\left(\pderivative{\interpolation{N}{\left(p\right)}}{\xi^{\iota}}\right)\mathbb{J}\vec{a}^{\iota}+\left(\pderivative{\interpolation{N}{\left(p\mathbb{J}\vec{a}^{\iota}\right)}}{\xi^{\iota}}\right),\vec{u}\right\rangle _{N} \\
&+\phantom{\frac{1}{2}} \sum_{\iota=1}^{3}\ \int\limits _{\partial\mathrm{E},N}\hat{n}^{\iota}\left\{ \mathbf{V}^{T}\left[\mathbb{J}\vec{a}^{\iota}\cdot\blockvec{\mathbf{G}}\right]-p\mathbb{J}\vec{a}^{\iota}\cdot\vec{u}\right\} \,dS \\
=&\phantom{+}
\frac{1}{2}\sum_{\iota=1}^{3}\left\langle p\left(\pderivative{\interpolation{N}{\left(\mathbb{J}\vec{a}^{\iota}\right)}}{\xi^{\iota}}\right)+\left(\pderivative{\interpolation{N}{\left(p\right)}}{\xi^{\iota}}\right)\mathbb{J}\vec{a}^{\iota}+\left(\pderivative{\interpolation{N}{\left(p\mathbb{J}\vec{a}^{\iota}\right)}}{\xi^{\iota}}\right),\vec{u}\right\rangle _{N}\\
&+\phantom{\frac{1}{2}}\sum_{\iota=1}^{3}\ \int\limits _{\partial\mathrm{E},N}\hat{s}n_{\iota}\left\{ \mathbf{V}^{T}\mathbf{G}_{\iota}-pu_{\iota}\right\} \,dS,
\end{split}
\end{align}
\noindent 
\textbf{Step 3:} Consider the discrete surface integral in equation \eqref{LocalDisKEP1}:
First of all, it follows
\begin{equation}
-\frac{1}{2}\left|\vec{u}\right|^{2}-\sum_{\upsilon=1}^{3}u_{\upsilon}\avg{u_{\upsilon}}=\sum_{\upsilon=1}^{3}\avg{u_{\upsilon}}^{2}-\frac{1}{2}\avg{u_{\upsilon}^{2}}=\frac{1}{2}\overline{\left|\vec{u}\right|}^{2}.
\end{equation}
Thus, Jameson's \cite{Jameson2008} conditions \eqref{Jameson} provide  
\begin{align}\label{SurfacKinetic1}
\begin{split}
\int\limits _{\partial\mathrm{E},N}\mathbf{V}^{T}\left(\tilde{\textbf{G}}_{\hat{n}}^{*}-\tilde{\textbf{G}}_{\hat{n}}\right)\,dS
=& \phantom{+}\sum_{\iota=1}^{3}\ \int\limits _{\partial\mathrm{E},N}\hat{s}n_{\iota}\mathbf{V}^{T}\left(\mathbf{G}_{\iota}^{*}-\mathbf{G}_{\iota}\right)\,dS \\
=&\phantom{+} \sum_{\iota=1}^{3}\ \int\limits _{\partial\mathrm{E},N}\hat{s}n_{\iota}\left[-\frac{1}{2}\left|\vec{u}\right|^{2}\mathbf{G}_{\iota}^{1,*}+\sum_{\upsilon=1}^{3}u_{\upsilon}\mathbf{G}_{\iota}^{\upsilon,*}-\mathbf{V}^{T}\mathbf{G}_{\iota}\right]\,dS \\
=& \phantom{+}\sum_{\iota=1}^{3}\ \int\limits _{\partial\mathrm{E},N}\hat{s}n_{\iota}\left[\frac{1}{2}\overline{\left|\vec{u}\right|}^{2}\mathbf{G}_{\iota}^{1,*}+\avg{p}^{\star}u_{\iota}-\mathbf{V}^{T}\mathbf{G}_{\iota}\right]\,dS,
\end{split}
\end{align}
where $\avg{p}^{\star}$ is a consistent numerical trace approximation of the pressure. 
\noindent
Finally, we combine \eqref{LocalDisKEP1} with the equations \eqref{TemporalKinetic1}, \eqref{SpatialKinetic9}, \eqref{SurfacKinetic1} and obtain the identity \eqref{LocalDisKEP}   
\begin{align}
\begin{split}
\pderivative{}{\tau}\left\langle \interpolation{N}{\left(\Bbbk\right)},\mathcal{J}\right\rangle _{N}=&
-\frac{1}{2}\sum_{\iota=1}^{3}\left\langle p\left(\pderivative{\interpolation{N}{\left(\mathbb{J}\vec{a}^{\iota}\right)}}{\xi^{\iota}}\right)+\left(\pderivative{\interpolation{N}{\left(p\right)}}{\xi^{\iota}}\right)\mathbb{J}\vec{a}^{\iota}+\left(\pderivative{\interpolation{N}{\left(p\mathbb{J}\vec{a}^{\iota}\right)}}{\xi^{\iota}}\right),\vec{u}\right\rangle _{N}\\
&-\sum_{\iota=1}^{3}\ \int\limits _{\partial\mathrm{E},N}\hat{s}n_{\iota}\left[\frac{1}{2}\overline{\left|\vec{u}\right|}^{2}\mathbf{G}_{\iota}^{1,*}+\left(\avg{p}^{\star}-p\right)u_{\iota}\right]\,dS.
\end{split}
\end{align}

\bibliography{References.bib}

\begin{thebibliography}{70}
\providecommand{\natexlab}[1]{#1}
\providecommand{\url}[1]{\texttt{#1}}
\expandafter\ifx\csname urlstyle\endcsname\relax
  \providecommand{\doi}[1]{doi: #1}\else
  \providecommand{\doi}{doi: \begingroup \urlstyle{rm}\Url}\fi

\bibitem[Barth(1999)]{Barth2018}
T.~J. Barth.
\newblock Numerical methods for gasdynamic systems on unstructured meshes.
\newblock In \emph{An introduction to recent developments in theory and
  numerics for conservation laws}, pages 195--285. Springer, 1999.

\bibitem[Bassi and Rebay(1997)]{Bassi1997}
F.~Bassi and S.~Rebay.
\newblock A high-order accurate discontinuous finite element method for the
  numerical solution of the compressible {N}avier--{S}tokes equations.
\newblock \emph{Journal of Computational Physics}, 131\penalty0 (2):\penalty0
  267--279, 1997.

\bibitem[Beck et~al.(2014)Beck, Bolemann, Flad, Frank, Gassner, Hindenlang, and
  Munz]{beck2014high}
A.~Beck, T.~Bolemann, D.~Flad, H.~Frank, G.~Gassner, F.~Hindenlang, and C.-D.
  Munz.
\newblock High-order discontinuous {G}alerkin spectral element methods for
  transitional and turbulent flow simulations.
\newblock \emph{International Journal for Numerical Methods in Fluids},
  76\penalty0 (8):\penalty0 522--548, 2014.

\bibitem[Blaisdell et~al.(1996)Blaisdell, Spyropoulos, and
  Qin]{blaisdell1996effect}
G.~A. Blaisdell, E.~T. Spyropoulos, and J.~H. Qin.
\newblock The effect of the formulation of nonlinear terms on aliasing errors
  in spectral methods.
\newblock \emph{Applied Numerical Mathematics}, 21\penalty0 (3):\penalty0
  207--219, 1996.

\bibitem[Canuto et~al.(2006)Canuto, Hussaini, Quarteroni, and Zang]{CHQZ:2006}
C.~Canuto, M.~Y. Hussaini, A.~Quarteroni, and T.~A. Zang.
\newblock \emph{Spectral methods}.
\newblock Springer, 2006.

\bibitem[Carton~de Wiart and Hillewaert(2012)]{de2012dns}
C.~Carton~de Wiart and K.~Hillewaert.
\newblock {DNS} and {ILES} of transitional flows around a {SD7003} using a high
  order discontinuous {G}alerkin method.
\newblock In \emph{Seventh International Conference on Computational Fluid
  Dynamics (ICCFD7), Big Island, Hawaii}, 2012.

\bibitem[Chandrashekar(2013)]{Chandrashekar2013}
P.~Chandrashekar.
\newblock Kinetic energy preserving and entropy stable finite volume schemes
  for compressible {E}uler and {N}avier--{S}tokes equations.
\newblock \emph{Communications in Computational Physics}, 14\penalty0
  (5):\penalty0 1252--1286, 2013.

\bibitem[Cockburn and Shu(2001)]{Cockburn2001}
B.~Cockburn and C.-W. Shu.
\newblock {R}unge--{K}utta discontinuous {G}alerkin methods for
  convection-dominated problems.
\newblock \emph{Journal of Scientific Computing}, 16\penalty0 (3):\penalty0
  173--261, 2001.

\bibitem[Dalcin et~al.(2019)Dalcin, Rojas, Zampini, Fern{\'a}ndez, Carpenter,
  and Parsani]{dalcin2019conservative}
L.~Dalcin, D.~Rojas, S.~Zampini, D.~C. Del~Rey Fern{\'a}ndez, M.~H. Carpenter,
  and M.~Parsani.
\newblock Conservative and entropy stable solid wall boundary conditions for
  the compressible {N}avier--{S}tokes equations: Adiabatic wall and heat
  entropy transfer.
\newblock \emph{Journal of Computational Physics}, 397:\penalty0 108775, 2019.

\bibitem[Donea et~al.(2017)Donea, Huerta, Ponthot, and
  Rodr{\'\i}guez-Ferran]{donea2017arbitrary}
J.~Donea, A.~Huerta, J.-P. Ponthot, and A.~Rodr{\'\i}guez-Ferran.
\newblock Arbitrary {L}agrangian--{E}ulerian methods.
\newblock \emph{Encyclopedia of Computational Mechanics Second Edition}, pages
  1--23, 2017.

\bibitem[Ducros et~al.(2000)Ducros, Laporte, Soul{\`e}res, Guinot, Moinat, and
  Caruelle]{ducros2000high}
F.~Ducros, F.~Laporte, T.~Soul{\`e}res, V.~Guinot, Ph. Moinat, and B.~Caruelle.
\newblock High-order fluxes for conservative skew-symmetric-like schemes in
  structured meshes: application to compressible flows.
\newblock \emph{Journal of Computational Physics}, 161\penalty0 (1):\penalty0
  114--139, 2000.

\bibitem[Fern{\'a}ndez et~al.(2014)Fern{\'a}ndez, Hicken, and
  Zingg]{Fernandez2014}
D.~C. Del~Rey Fern{\'a}ndez, J.~E. Hicken, and D.~W. Zingg.
\newblock Review of summation-by-parts operators with simultaneous
  approximation terms for the numerical solution of partial differential
  equations.
\newblock \emph{Computers \& Fluids}, 95:\penalty0 171--196, 2014.

\bibitem[Fisher and Carpenter(2013)]{Fisher2013}
T.~C. Fisher and M.~H. Carpenter.
\newblock High-order entropy stable finite difference schemes for nonlinear
  conservation laws: Finite domains.
\newblock \emph{Journal of Computational Physics}, 252:\penalty0 518--557,
  2013.

\bibitem[Flad and Gassner(2017)]{flad2017use}
D.~Flad and G.~J. Gassner.
\newblock On the use of kinetic energy preserving {DG}-schemes for large eddy
  simulation.
\newblock \emph{Journal of Computational Physics}, 350:\penalty0 782--795,
  2017.

\bibitem[Flad et~al.(2016)Flad, Beck, and Munz]{flad2016simulation}
D.~Flad, A.~D. Beck, and C.-D. Munz.
\newblock Simulation of underresolved turbulent flows by adaptive filtering
  using the high order discontinuous {G}alerkin spectral element method.
\newblock \emph{Journal of Computational Physics}, 313:\penalty0 1--12, 2016.

\bibitem[Galbraith and Visbal(2010)]{galbraith2010implicit}
M.~Galbraith and M.~Visbal.
\newblock Implicit large eddy simulation of low-{R}eynolds-number transitional
  flow past the {SD7003} airfoil.
\newblock In \emph{40th Fluid Dynamics Conference and Exhibit}, page 4737,
  2010.

\bibitem[Gassner(2013)]{Gassner2013}
G.~J. Gassner.
\newblock A skew-symmetric discontinuous {G}alerkin spectral element
  discretization and its relation to {SBP}-{SAT} finite difference methods.
\newblock \emph{SIAM Journal on Scientific Computing}, 35\penalty0
  (3):\penalty0 A1233--A1253, 2013.

\bibitem[Gassner(2014)]{gassner2014kinetic}
G.~J Gassner.
\newblock A kinetic energy preserving nodal discontinuous {G}alerkin spectral
  element method.
\newblock \emph{International Journal for Numerical Methods in Fluids},
  76\penalty0 (1):\penalty0 28--50, 2014.

\bibitem[Gassner and Beck(2013)]{gassner2013accuracy}
G.~J. Gassner and A.~D. Beck.
\newblock On the accuracy of high-order discretizations for underresolved
  turbulence simulations.
\newblock \emph{Theoretical and Computational Fluid Dynamics}, 27\penalty0
  (3-4):\penalty0 221--237, 2013.

\bibitem[Gassner and Kopriva(2011)]{gassner2011comparison}
G.~J. Gassner and D.~A. Kopriva.
\newblock A comparison of the dispersion and dissipation errors of {G}auss and
  {G}auss--{L}obatto discontinuous {G}alerkin spectral element methods.
\newblock \emph{SIAM Journal on Scientific Computing}, 33\penalty0
  (5):\penalty0 2560--2579, 2011.

\bibitem[Gassner et~al.(2016)Gassner, Winters, and Kopriva]{Gassner2016}
G.~J. Gassner, A.~R. Winters, and D.~A. Kopriva.
\newblock Split form nodal discontinuous {G}alerkin schemes with
  summation-by-parts property for the compressible {E}uler equations.
\newblock \emph{Journal of Computational Physics}, 327:\penalty0 39--66, 2016.

\bibitem[Gassner et~al.(2018)Gassner, Winters, Hindenlang, and
  Kopriva]{Gassner2017}
G.~J. Gassner, A.~R. Winters, F.~J. Hindenlang, and D.~A. Kopriva.
\newblock The {BR1} scheme is stable for the compressible {N}avier--{S}tokes
  equations.
\newblock \emph{Journal of Scientific Computing}, 77\penalty0 (1):\penalty0
  154--200, 2018.

\bibitem[Harten(1983)]{Harten1983}
A.~Harten.
\newblock On the symmetric form of systems of conservation laws with entropy.
\newblock \emph{Journal of Computational Physics}, 49:\penalty0 151--164, 1983.

\bibitem[Hesthaven and Warburton(2007)]{hesthaven2007nodal}
J.~S. Hesthaven and T.~Warburton.
\newblock \emph{Nodal discontinuous {G}alerkin methods: algorithms, analysis,
  and applications}.
\newblock Springer Science \& Business Media, 2007.

\bibitem[Hughes et~al.(1986)Hughes, Franca, and Mallet]{Hughes1986}
T.~J. Hughes, L.~P. Franca, and M.~Mallet.
\newblock A new finite element formulation for computational fluid dynamics: I.
  symmetric forms of the compressible {E}uler and {N}avier--{S}tokes equations
  and the second law of thermodynamics.
\newblock \emph{Computer Methods in Applied Mechanics and Engineering},
  54\penalty0 (2):\penalty0 223--234, 1986.

\bibitem[Ismail and Roe(2009)]{ismail2009affordable}
F.~Ismail and P.~L. Roe.
\newblock Affordable, entropy-consistent {E}uler flux functions ii: Entropy
  production at shocks.
\newblock \emph{Journal of Computational Physics}, 228\penalty0 (15):\penalty0
  5410--5436, 2009.

\bibitem[Jameson(2008)]{Jameson2008}
A.~Jameson.
\newblock Formulation of kinetic energy preserving conservative schemes for gas
  dynamics and direct numerical simulation of one-dimensional viscous
  compressible flow in a shock tube using entropy and kinetic energy preserving
  schemes.
\newblock \emph{Journal of Scientific Computing}, 34\penalty0 (2):\penalty0
  188--208, 2008.

\bibitem[Jones et~al.(1998)Jones, Dohring, and Platzer]{jones1998experimental}
K.~D. Jones, C.~M. Dohring, and M.~F. Platzer.
\newblock Experimental and computational investigation of the {K}noller--{B}etz
  effect.
\newblock \emph{AIAA journal}, 36\penalty0 (7):\penalty0 1240--1246, 1998.

\bibitem[Kennedy and Gruber(2008)]{Kennedy2008}
C.~A. Kennedy and A.~Gruber.
\newblock Reduced aliasing formulations of the convective terms within the
  {N}avier--{S}tokes equations for a compressible fluid.
\newblock \emph{Journal of Computational Physics}, 227\penalty0 (3):\penalty0
  1676--1700, 2008.

\bibitem[Kennedy et~al.(2000)Kennedy, Carpenter, and Lewis]{kennedy2000low}
C.~A. Kennedy, M.~H. Carpenter, and R.~M. Lewis.
\newblock Low-storage, explicit {R}unge--{K}utta schemes for the compressible
  {N}avier--{S}tokes equations.
\newblock \emph{Applied Numerical Mathematics}, 35\penalty0 (3):\penalty0
  177--219, 2000.

\bibitem[Kirby and Karniadakis(2003)]{kirby2003aliasing}
R.~M. Kirby and G.~E. Karniadakis.
\newblock De-aliasing on non-uniform grids: algorithms and applications.
\newblock \emph{Journal of Computational Physics}, 191\penalty0 (1):\penalty0
  249--264, 2003.

\bibitem[Kopriva(2006)]{Kopriva2006}
D.~A. Kopriva.
\newblock Metric identities and the discontinuous spectral element method on
  curvilinear meshes.
\newblock \emph{Journal of Scientific Computing}, 26\penalty0 (3):\penalty0
  301, 2006.

\bibitem[Kopriva(2009)]{Kopriva2009}
D.~A. Kopriva.
\newblock \emph{Implementing spectral methods for partial differential
  equations: Algorithms for scientists and engineers}.
\newblock Springer Science \& Business Media, 2009.

\bibitem[Kopriva(2018)]{kopriva2018stability}
D.~A. Kopriva.
\newblock Stability of overintegration methods for nodal discontinuous
  {G}alerkin spectral element methods.
\newblock \emph{Journal of Scientific Computing}, 76\penalty0 (1):\penalty0
  426--442, 2018.

\bibitem[Kopriva et~al.(2016)Kopriva, Winters, Bohm, and
  Gassner]{kopriva2016provably}
D.~A. Kopriva, A.~R. Winters, M.~Bohm, and G.~J. Gassner.
\newblock A provably stable discontinuous {G}alerkin spectral element
  approximation for moving hexahedral meshes.
\newblock \emph{Computers \& Fluids}, 139:\penalty0 148--160, 2016.

\bibitem[Krais et~al.(2019)Krais, Beck, Bolemann, Frank, Flad, Gassner,
  Hindenlang, Hoffmann, Kuhn, Sonntag, and Munz]{krais2019flexi}
N.~Krais, A.~Beck, T.~Bolemann, H.~Frank, D.~Flad, G.~Gassner, F.~Hindenlang,
  M.~Hoffmann, T.~Kuhn, M.~Sonntag, and C.-D. Munz.
\newblock {FLEXI}: A high order discontinuous {G}alerkin framework for
  hyperbolic-parabolic conservation laws.
\newblock \emph{arXiv preprint arXiv:1910.02858}, 2019.

\bibitem[Kreiss and Oliger(1972)]{kreiss1}
H.-O. Kreiss and J.~Oliger.
\newblock Comparison of accurate methods for the integration of hyperbolic
  equations.
\newblock \emph{Tellus}, 24\penalty0 (3):\penalty0 199--215, 1972.

\bibitem[Kuya et~al.(2018)Kuya, Totani, and Kawai]{kuya2018kinetic}
Y.~Kuya, K.~Totani, and S.~Kawai.
\newblock Kinetic energy and entropy preserving schemes for compressible flows
  by split convective forms.
\newblock \emph{Journal of Computational Physics}, 375:\penalty0 823--853,
  2018.

\bibitem[Lefloch et~al.(2002)Lefloch, Mercier, and Rohde]{lefloch2002fully}
P.~G. Lefloch, J.-M. Mercier, and C.~Rohde.
\newblock Fully discrete, entropy conservative schemes of arbitrary order.
\newblock \emph{SIAM Journal on Numerical Analysis}, 40\penalty0 (5):\penalty0
  1968--1992, 2002.

\bibitem[Lele(1992)]{lele1992compact}
S.~K. Lele.
\newblock Compact finite difference schemes with spectral-like resolution.
\newblock \emph{Journal of Computational Physics}, 103\penalty0 (1):\penalty0
  16--42, 1992.

\bibitem[Lomtev et~al.(1999)Lomtev, Kirby, and
  Karniadakis]{lomtev1999discontinuous}
I.~Lomtev, R.~M. Kirby, and G.~E. Karniadakis.
\newblock A discontinuous {G}alerkin {ALE} method for compressible viscous
  flows in moving domains.
\newblock \emph{Journal of Computational Physics}, 155\penalty0 (1):\penalty0
  128--159, 1999.

\bibitem[McGowan et~al.(2008)McGowan, Gopalarathnam, Ol, Edwards, and
  Fredberg]{mcgowan2008computation}
G.~McGowan, A.~Gopalarathnam, M.~Ol, J.~Edwards, and D.~Fredberg.
\newblock Computation vs. experiment for high-frequency low-{R}eynolds number
  airfoil pitch and plunge.
\newblock In \emph{46th AIAA Aerospace Sciences Meeting and Exhibit}, page 653,
  2008.

\bibitem[Mengaldo et~al.(2015)Mengaldo, De~Grazia, Moxey, Vincent, and
  Sherwin]{mengaldo2015dealiasing}
G.~Mengaldo, D.~De~Grazia, D.~Moxey, P.~E. Vincent, and S.~J. Sherwin.
\newblock Dealiasing techniques for high-order spectral element methods on
  regular and irregular grids.
\newblock \emph{Journal of Computational Physics}, 299:\penalty0 56--81, 2015.

\bibitem[Merriam(1989)]{Merriam1989}
M.~L. Merriam.
\newblock Towards a rigorous approach to artificial dissipation.
\newblock Technical report, National Aeronautics and Space Administration,
  Moffett Field, CA (USA). Ames Research Center, 1989.

\bibitem[Minoli and Kopriva(2011)]{minoli2011discontinuous}
C.~A.~A. Minoli and D.~A. Kopriva.
\newblock Discontinuous {G}alerkin spectral element approximations on moving
  meshes.
\newblock \emph{Journal of Computational Physics}, 230\penalty0 (5):\penalty0
  1876--1902, 2011.

\bibitem[Morinishi(2010)]{morinishi2010skew}
Y.~Morinishi.
\newblock Skew-symmetric form of convective terms and fully conservative finite
  difference schemes for variable density low-{M}ach number flows.
\newblock \emph{Journal of Computational Physics}, 229\penalty0 (2):\penalty0
  276--300, 2010.

\bibitem[Moura et~al.(2015)Moura, Sherwin, and Peir{\'o}]{moura2015linear}
R.~C. Moura, S.~J. Sherwin, and J.~Peir{\'o}.
\newblock Linear dispersion--diffusion analysis and its application to
  under-resolved turbulence simulations using discontinuous {G}alerkin
  spectral/hp methods.
\newblock \emph{Journal of Computational Physics}, 298:\penalty0 695--710,
  2015.

\bibitem[Moura et~al.(2017)Moura, Mengaldo, Peir{\'o}, and
  Sherwin]{moura2017eddy}
R.~C. Moura, G.~Mengaldo, J.~Peir{\'o}, and S.~J. Sherwin.
\newblock On the eddy-resolving capability of high-order discontinuous
  {G}alerkin approaches to implicit {LES}/under-resolved {DNS} of {E}uler
  turbulence.
\newblock \emph{Journal of Computational Physics}, 330:\penalty0 615--623,
  2017.

\bibitem[Nguyen(2010)]{nguyen2010arbitrary}
V.-T. Nguyen.
\newblock An arbitrary {L}agrangian--{E}ulerian discontinuous {G}alerkin method
  for simulations of flows over variable geometries.
\newblock \emph{Journal of Fluids and Structures}, 26\penalty0 (2):\penalty0
  312--329, 2010.

\bibitem[Nigro et~al.(2019)Nigro, De~Bartolo, Crivellini, Franciolini, Colombo,
  and Bassi]{nigro2019low}
A.~Nigro, C.~De~Bartolo, A.~Crivellini, M.~Franciolini, A.~Colombo, and
  F.~Bassi.
\newblock A low-dissipation {DG} method for the under-resolved simulation of
  low {M}ach number turbulent flows.
\newblock \emph{Computers \& Mathematics with Applications}, 77\penalty0
  (6):\penalty0 1739--1755, 2019.

\bibitem[Nordstr{\"o}m and Carpenter(1999)]{nordstrom1999boundary}
J.~Nordstr{\"o}m and M.~H. Carpenter.
\newblock Boundary and interface conditions for high-order finite-difference
  methods applied to the {E}uler and {N}avier--{S}tokes equations.
\newblock \emph{Journal of Computational Physics}, 148\penalty0 (2):\penalty0
  621--645, 1999.

\bibitem[Ol et~al.(2009)Ol, Reeder, Fredberg, McGowan, Gopalarathnam, and
  Edwards]{ol2009computation}
M.~V. Ol, M.~Reeder, D.~Fredberg, G.~Z. McGowan, A.~Gopalarathnam, and J.~R.
  Edwards.
\newblock Computation vs. experiment for high-frequency low-{R}eynolds number
  airfoil plunge.
\newblock \emph{International Journal of Micro Air Vehicles}, 1\penalty0
  (2):\penalty0 99--119, 2009.

\bibitem[Parsani et~al.(2015)Parsani, Carpenter, and
  Nielsen]{parsani2015entropy}
M.~Parsani, M.~H. Carpenter, and E.~J. Nielsen.
\newblock Entropy stable wall boundary conditions for the three-dimensional
  compressible {N}avier--s{t}okes equations.
\newblock \emph{Journal of Computational Physics}, 292:\penalty0 88--113, 2015.

\bibitem[Persson et~al.(2009)Persson, Bonet, and
  Peraire]{persson2009discontinuous}
P.-O. Persson, J.~Bonet, and J.~Peraire.
\newblock Discontinuous {G}alerkin solution of the {N}avier--{S}tokes equations
  on deformable domains.
\newblock \emph{Computer Methods in Applied Mechanics and Engineering},
  198\penalty0 (17-20):\penalty0 1585--1595, 2009.

\bibitem[Pirozzoli(2010)]{Pirozzoli2010}
S.~Pirozzoli.
\newblock Generalized conservative approximations of split convective
  derivative operators.
\newblock \emph{Journal of Computational Physics}, 229\penalty0 (19):\penalty0
  7180--7190, 2010.

\bibitem[Radespiel et~al.(2007)Radespiel, Windte, and
  Scholz]{radespiel2007numerical}
R.~Radespiel, J.~Windte, and U.~Scholz.
\newblock Numerical and experimental flow analysis of moving airfoils with
  laminar separation bubbles.
\newblock \emph{AIAA journal}, 45\penalty0 (6):\penalty0 1346--1356, 2007.

\bibitem[Ranocha(2018)]{Ranocha2018}
H.~Ranocha.
\newblock \emph{Generalised summation-by-parts operators and entropy stability
  of numerical methods for hyperbolic balance laws}.
\newblock Cuvillier Verlag, 2018.

\bibitem[Ray(2017)]{Ray2017}
D.~Ray.
\newblock \emph{Entropy-stable finite difference and finite volume schemes for
  compressible flows}.
\newblock PhD thesis, Tata Institute of Fundamental Research, Mumbai, 2017.

\bibitem[Schn{\"u}cke et~al.(2020)Schn{\"u}cke, Krais, Bolemann, and
  Gassner]{Schnuecke2018}
G.~Schn{\"u}cke, N.~Krais, T.~Bolemann, and G.~J. Gassner.
\newblock Entropy stable discontinuous galerkin schemes on moving meshes for
  hyperbolic conservation laws.
\newblock \emph{Journal of Scientific Computing}, 2020.
\newblock \doi{10.1007/s10915-020-01171-7}.

\bibitem[Shu et~al.(2005)Shu, Don, Gottlieb, Schilling, and
  Jameson]{shu2005numerical}
C.-W. Shu, W.-S. Don, D.~Gottlieb, O.~Schilling, and L.~Jameson.
\newblock Numerical convergence study of nearly incompressible, inviscid
  {T}aylor--{G}reen vortex flow.
\newblock \emph{Journal of Scientific Computing}, 24\penalty0 (1):\penalty0
  1--27, 2005.

\bibitem[Stoellinger et~al.(2019)Stoellinger, Edmonds, Kirby, Mavriplis, and
  Heinz]{stoellinger2019dynamic}
M.~K. Stoellinger, A.~P. Edmonds, A.~C. Kirby, D.~J. Mavriplis, and S.~Heinz.
\newblock Dynamic {SGS} modeling in {LES} using {DG} with kinetic energy
  preserving flux schemes.
\newblock In \emph{AIAA Scitech 2019 Forum}, page 1648, 2019.

\bibitem[Sv{\"a}rd and {\"O}zcan(2014)]{svard2014entropy}
M.~Sv{\"a}rd and H.~{\"O}zcan.
\newblock Entropy-stable schemes for the {E}uler equations with far-field and
  wall boundary conditions.
\newblock \emph{Journal of Scientific Computing}, 58\penalty0 (1):\penalty0
  61--89, 2014.

\bibitem[Tadmor(2003)]{Tadmor2003}
E.~Tadmor.
\newblock Entropy stability theory for difference approximations of nonlinear
  conservation laws and related time-dependent problems.
\newblock \emph{Acta Numerica}, 12:\penalty0 451--512, 2003.

\bibitem[Tang(2005)]{tang2005moving}
T.~Tang.
\newblock Moving mesh methods for computational fluid dynamics.
\newblock \emph{Contemporary mathematics}, 383\penalty0 (8):\penalty0 141--173,
  2005.

\bibitem[Uranga et~al.(2011)Uranga, Persson, Drela, and
  Peraire]{uranga2011implicit}
A.~Uranga, P.-O. Persson, M.~Drela, and J.~Peraire.
\newblock Implicit large eddy simulation of transition to turbulence at low
  {R}eynolds numbers using a discontinuous {G}alerkin method.
\newblock \emph{International Journal for Numerical Methods in Engineering},
  87\penalty0 (1-5):\penalty0 232--261, 2011.

\bibitem[Visbal(2009)]{visbal2009high}
M.~R. Visbal.
\newblock High-fidelity simulation of transitional flows past a plunging
  airfoil.
\newblock \emph{AIAA journal}, 47\penalty0 (11):\penalty0 2685--2697, 2009.

\bibitem[Winters and Kopriva(2014)]{winters2014ale}
A.~R. Winters and D.~A. Kopriva.
\newblock {ALE}--{DGSEM} approximation of wave reflection and transmission from
  a moving medium.
\newblock \emph{Journal of Computational Physics}, 263:\penalty0 233--267,
  2014.

\bibitem[Winters et~al.(2018)Winters, Moura, Mengaldo, Gassner, Walch, Peiro,
  and Sherwin]{winters2018comparative}
A.~R. Winters, R.~C. Moura, G.~Mengaldo, G.~J. Gassner, S.~Walch, J.~Peiro, and
  S.~J. Sherwin.
\newblock A comparative study on polynomial dealiasing and split form
  discontinuous {G}alerkin schemes for under-resolved turbulence computations.
\newblock \emph{Journal of Computational Physics}, 372:\penalty0 1--21, 2018.

\bibitem[Yamaleev et~al.(2019)Yamaleev, Fernandez, Lou, and
  Carpenter]{yamaleev2019entropy}
N.~K. Yamaleev, D.~C. Del~Rey Fernandez, J.~Lou, and M.~H. Carpenter.
\newblock Entropy stable spectral collocation schemes for the {3-D}
  {N}avier--{S}tokes equations on dynamic unstructured grids.
\newblock \emph{Journal of Computational Physics}, 399:\penalty0 108897, 2019.

\bibitem[Yuan et~al.(2005)Yuan, Khalid, Windte, Scholz, and
  Radespiel]{yuan2005investigation}
W.~Yuan, M.~Khalid, J.~Windte, U.~Scholz, and R.~Radespiel.
\newblock An investigation of low-{R}eynolds-number flows past airfoils.
\newblock In \emph{23rd AIAA Applied Aerodynamics Conference}, page 4607, 2005.

\end{thebibliography}
\bibliographystyle{plainnat}
\end{document}